\documentclass[12pt,a4paper]{article}
\usepackage[utf8]{inputenc}
\usepackage[total={6in, 8in}]{geometry}
\usepackage[english]{babel}
\usepackage[round]{natbib}
\usepackage{amsmath, amssymb,amsthm,bm,dsfont,tabularx}
\usepackage{spverbatim,hhline}
\usepackage{graphicx,accents,url}
\usepackage{blkarray}
\usepackage{makecell}
\newcolumntype{Y}{>{\centering\arraybackslash}X}
\newcolumntype{S}{>{\hsize=.1\hsize}Y}
\makeatletter
\newcommand{\pushright}[1]{\ifmeasuring@#1\else\omit\hfill$\displaystyle#1$\fi\ignorespaces}
\newcommand{\pushleft}[1]{\ifmeasuring@#1\else\omit$\displaystyle#1$\hfill\fi\ignorespaces}
\makeatother
\newcommand\smallO{
  \mathchoice
    {{\scriptstyle\mathcal{O}}}% \displaystyle
    {{\scriptstyle\mathcal{O}}}% \textstyle
    {{\scriptscriptstyle\mathcal{O}}}% \scriptstyle
    {\scalebox{.7}{$\scriptscriptstyle\mathcal{O}$}}%\scriptscriptstyle
  }
\emergencystretch=\maxdimen
\title{\textbf{\Large Parametric dependence between random vectors via copula-based divergence measures}}
\author{Steven De Keyser, Ir\`{e}ne Gijbels}
\date{\today}
\begin{document}
\maketitle 
\vspace*{2cm}
\noindent 
\textbf{Abstract.} This article proposes copula-based dependence quantification between multiple groups of random variables of possibly different sizes via the family of $\Phi$-divergences. An axiomatic framework for this purpose is provided, after which we focus on the absolutely continuous setting assuming copula densities exist. We consider parametric and semi-parametric frameworks, discuss estimation procedures, and report on asymptotic properties of the proposed estimators. In particular, we first concentrate on a Gaussian copula approach yielding explicit and attractive dependence coefficients for specific choices of $\Phi$, which are more amenable for estimation. Next, general parametric copula families are considered, with special attention to nested Archimedean copulas, being a natural choice for dependence modelling of random vectors. The results are illustrated by means of examples. Simulations and a real-world application on financial data are provided as well. \newline  
\par\noindent\rule{\textwidth}{0.4pt}
\clearpage \noindent 
\textbf{1. Introduction} 
\\

The fundamental problem of measuring dependence between two random variables is customary in the analysis of bivariate data. Linear relationships are embodied in the Pearson correlation coefficient and concordance measures like Kendall's tau or Spearman's rho, among many others, extend to incorporate monotone dependence. The interest of generalizations of such concordance measures to more than two univariate random variables is also widely recognized, see e.g. \cite{Nelsen1996}, \cite{Schmid2007} and \cite{Gijbels2021}. Important is the compliance with certain postulated axioms, starting with \cite{Renyi1959} and followed by e.g. \cite{Lancaster1963}, \cite{Schweizer1981} and \cite{Embrechts2002} for the case of two univariate random variables, and e.g. \cite{Wolff1980}, \cite{Nelsen1996} and \cite{Gijbels2021} when the interest is in more than two variables.

Another extension consists of looking at two groups of random variables. In this context, one is typically aware of the statistical analysis of canonical correlations, as in \cite{Hotelling1936}. \cite{Grothe2014} suggest using concordance measures and \cite{Mordant2022} measure dependence between two random vectors via optimal transport, making a Gaussian assumption when going to statistical inference. Often, dependence capturing is restricted to monotone relationships, either due to making rather stringent assumptions (e.g. Gaussianity), or because the measures in question have limited detection ability (e.g. concordance measures do not measure tail dependence). 

\cite{Gijbels2023} work in the broader setting of $k$ random vectors, think of e.g. answers to $k$ different questionnaires or $k$ groups of financial assets like shares from different stock indexes, and define the general family of $\Phi$-dependence measures, which complies with their postulated properties that are driven by the objective of quantifying any deviation from independence. In this article, we elaborate more on these dependence measures. After proving several desirable properties, we focus on some examples of maximal dependence in singular copula distributions. Thereafter, we assume absolute continuity and concentrate on parametric and semi-parametric modelling and estimation from a copula density point of view. 

Unlike numerous dependence measures (concordance measures, tail dependence coefficients, $L_{p}$-copula distances like the Hoeffding's $\Phi^{2}$ of \cite{Ruppert2010}, $\dots$) inquiring about the (bounded) copula cdf, we now have functionals of the density. Copula densities may have rather cumbersome mathematical expressions, but numerical approximations are at hand if the dependence coefficient has no explicit analytical form in terms of the copula parameters. We also take extra care at boundaries, where copula densities commonly tend to infinity or zero.

The family of $\Phi$-dependence measures includes many popular measures that have a strong ability to detect deviations from independence, and there exists a great deal of statistics providing inference procedures and pursuing practical usefulness. The outline of this paper is the following.

Section 2 discusses possible axioms for dependence measures in the general context of $k$ random vectors, and they are verified for the proposed dependence measures is Section 3. A Gaussian copula approach and corresponding statistical inference is considered in Section 4, after which general parametric copulas are dealt with in Section 5, where the focus will be on maximum likelihood estimation. Some simulation studies are presented in Section 6, and a real life application to financial data is to be found in Section 7. We end this paper with a brief discussion in Section 8. For proofs related to the asymptotic properties of the proposed estimators (Theorems 1 and 2), we refer to the Appendix. 
\\ \newline \noindent 
\textbf{2. Notation and axioms} 
\\

The general setting is the same is in \cite{Gijbels2023}, i.e. we consider a $q$-dimensional random vector $\mathbf{X} = (\mathbf{X}_{1},\dots,\mathbf{X}_{k})$ having $k$ marginal random vectors $\mathbf{X}_{i} = (X_{i1},\dots,X_{id_{i}})$ for $i = 1,\dots,k$ composed of $d_{i}$ marginal univariate random variables $X_{ij}$ for $j = 1,\dots,d_{i}$, with $q = d_{1} + \cdots + d_{k}$, which are assumed to be continuous. The interest is in dependence measures $\mathcal{D}^{d_{1},\dots,d_{k}}(\mathbf{X}) = \mathcal{D}(\mathbf{X}_{1},\dots,\mathbf{X}_{k})$. We have $q$ continuous marginal cdf's, say $F_{ij}$, of $X_{ij}$ for $i = 1,\dots,k$ and $j = 1,\dots,d_{i}$. Sklar's theorem (\cite{Sklar1959}) guarantees the existence of a unique ($q$-dimensional) copula $C$ of $\mathbf{X}$ and marginal ($d_{i}$-dimensional) copulas $C_{i}$ of $\mathbf{X}_{i}$ for $i = 1,\dots,k$. They bring forth respective probability measures $\mu_{C}$ and $\mu_{C_{i}}$. Plausible axioms for a valid dependence measure are as follows.
\begin{enumerate} 
\addtolength{\itemindent}{0.4cm}
\item[(A1)]{For every permutation $\pi$ of $\mathbf{X}_{1}, \dots, \mathbf{X}_{k}$: $\mathcal{D}^{d_{1},\dots,d_{k}}(\mathbf{X}) = \mathcal{D}^{d_{1},\dots,d_{k}} \big (\pi(\mathbf{X}) \big )$; and \hspace*{0.2cm} for every permutation $\pi_{i}$ of $X_{i1},\dots,X_{id_{i}}$, for $i \in \{1, \dots, k \}$, it holds: \hspace*{0.3cm}
 $\mathcal{D}^{d_{1},\dots,d_{k}}(\mathbf{X}) = \mathcal{D} \big (\mathbf{X}_{1}, \dots, \pi_{i}(\mathbf{X}_{i}), \dots, \mathbf{X}_{k} \big )$.}
 \vspace{0.12cm}
\item[(A2)]{$0 \leq \mathcal{D}^{d_{1},\dots,d_{k}}(\mathbf{X}) \leq 1$.}
\vspace{0.12cm}
\item[(A3)]{$\mathcal{D}^{d_{1},\dots,d_{k}}(\mathbf{X}) = 0$ if and only if $\mathbf{X}_{1}, \dots, \mathbf{X}_{k}$ are mutually independent.}
\vspace{0.12cm}
\item[(A4)]{$\mathcal{D} (\mathbf{X}_{1},\dots,\mathbf{X}_{k}, \mathbf{X}_{k+1}) \geq \mathcal{D}(\mathbf{X}_{1},\dots,\mathbf{X}_{k})$ with equality if and only if $\mathbf{X}_{k+1}$ is \hspace*{0.3cm} independent of $(\mathbf{X}_{1},\dots,\mathbf{X}_{k})$.}
\item[(A5)]{$\mathcal{D}^{d_{1},\dots,d_{k}}(\mathbf{X})$ is well defined for any $q$-dimensional random vector $\mathbf{X}$ and is a \hspace*{0.3cm} functional of solely the copula $C$ of $\mathbf{X}$.}
\vspace{0.12cm}
\item[(A6)]{Let $T_{ij}$ for $i = 1, \dots, k$ and $j = 1, \dots, d_{i}$ be strictly increasing, continuous \hspace*{0.4cm}transformations. Then $$\mathcal{D} \big (T_{1}(\mathbf{X}_{1}),\dots,T_{k}(\mathbf{X}_{k}) \big ) = \mathcal{D}(\mathbf{X}_{1},\dots,\mathbf{X}_{k}),$$
\hspace*{0.3cm} where $T_{i}(\mathbf{X}_{i}) = (T_{i1}(X_{i1}),\dots,T_{id_{i}}(X_{id_{i}}))$ for $i = 1,\dots, k$}.
\vspace{0.12cm}
\item[(A7)]{Let $T_{ij}$ be a strictly decreasing, continuous transformation for a fixed $i \in \hspace*{0.3cm} \{1,\dots,k\}$ and a fixed $j \in \{1,\dots,d_{i} \}$. Then $$\mathcal{D} \big ( \mathbf{X}_{1},\dots, T_{i}(\mathbf{X}_{i}), \dots,\mathbf{X}_{k} \big ) = \mathcal{D}(\mathbf{X}_{1},\dots,\mathbf{X}_{k}),$$ \hspace*{0.3cm} where $T_{i}(\mathbf{X}_{i}) = (X_{i1},\dots,T_{ij}(X_{ij}),\dots,X_{id_{i}})$}.
\vspace{0.12cm}
\item[(A8)] Let $(\mathbf{X}_n)_{n \in \mathbb{N}}$ be a sequence of $q$-dimensional random vectors having copulas \hspace*{0.3cm} $(C_{n})_{n \in \mathbb{N}}$, then 
$$\lim_{n \to \infty} \mathcal{D}^{d_{1},\dots,d_{k}}(\mathbf{X}_{n}) = \mathcal{D}^{d_{1},\dots,d_{k}}(\mathbf{X})$$
\hspace*{0.40cm} if $C_{n} \to C$ uniformly, where $C$ denotes the copula of $\mathbf{X}$.
\end{enumerate}

We now bring forward the family of $\Phi$-dependence measures and show its compliance with the above properties. For (A8), we will restrict ourselves to uniform convergence of the copula densities. 
\newline \\ \noindent
\textbf{3. $\Phi$-dependence measures}
\\

Write $\mu_{C} = \mu_{C}^{\text{ac}} + \mu_{C}^{\text{s}}$ for the Lebesgue decomposition of $\mu_{C}$ with respect to the product measure $\mu_{C_{1}} \times \cdots \times \mu_{C_{k}}$, i.e.
$\mu_{C}^{\text{ac}}$ is absolutely continuous with respect to $\mu_{C_{1}} \times \cdots \times \mu_{C_{k}}$ (denoted as $\mu_{C}^{\text{ac}} \ll   \mu_{C_{1}} \times \cdots \times \mu_{C_{k}}$) and $\mu_{C}^{\text{s}}$ is singular with respect to $ \mu_{C_{1}} \times \cdots \times \mu_{C_{k}}$ (denoted as $\mu_{C}^{\text{s}} \perp   \mu_{C_{1}} \times \cdots \times \mu_{C_{k}}$).
\newline \\ \noindent
\textit{3.1. Definition and properties}

The family of $\Phi$-dependence measures between $k$ random vectors is defined in \cite{Gijbels2023} as follows.
\newline \\ \noindent 
\textbf{Definition 1. ($\Phi$-dependence measures)}
Consider a continuous, convex function $\Phi: (0,\infty) \rightarrow \mathbb{R}$ with $\Phi(1) = 0$. Extend $\Phi$ by defining
\begin{equation*}
    \Phi(0)  = \lim_{\substack{t \to 0 \\ >}} \Phi(t),  
\qquad \qquad 
\mbox{and} \qquad \qquad  
    \Phi^{*}(0)  = \lim_{t \to \infty} \frac{\Phi(t)}{t}.
\end{equation*}
The \textit{$\Phi$-dependence} between $\mathbf{X}_{1},\dots,\mathbf{X}_{k}$ is the quantity $\mathcal{D}_{\Phi} = \mathcal{D}_{\Phi}(\mathbf{X}_{1},\dots,\mathbf{X}_{k}) \in [0,\infty]$ defined by
\begin{equation*}
\mathcal{D}_{\Phi} =
 \int  \Phi \bigg ( \frac{d\mu_{C}^{\text{ac}}}{d(\mu_{C_{1}} \times \cdots \times \mu_{C_{k}})} \bigg) d(\mu_{C_{1}} \times \cdots \times \mu_{C_{k}})
 + \Phi^{*}(0) \mu_{C}^{\text{s}}(B),
\end{equation*}
with $B$ the set on which $\mu_{C}^{s}$ is concentrated. 
We use the convention $0 \cdot \infty = 0$. 
\newline 

The maximum value of $\mathcal{D}_{\Phi}$ is $\Phi(0) + \Phi^{*}(0)$, and attained when $\mu_{C}=\mu_{C}^{\text{s}}$. After looking at some examples of maximal $\Phi$-dependence in singular copulas, we restrict ourselves in this paper to $\mu_{C} \ll \lambda^{q}$ with $\lambda^{q}$ the Lebesgue measure (hence, $\mu_{C} = \mu_{C}^{\text{ac}}$ and $\mu_{C}^{\text{s}} = 0$), implying that 
\begin{equation}\label{eq: phidiv absc}
    \mathcal{D}_{\Phi} \left (\mathbf{X}_{1},\dots,\mathbf{X}_{k} \right ) = \int_{\mathbb{I}^{q}} \prod_{i=1}^{k} c_{i}(\mathbf{u}_{i}) \Phi \left (\frac{c(\mathbf{u})}{\prod_{i=1}^{k} c_{i}(\mathbf{u}_{i})} \right )d\mathbf{u},
\end{equation},
    where $\mathbb{I} = [0,1]$, $c$ and $c_{i}$ the copula densities (w.r.t. $\lambda^{q}$ and $\lambda^{d_{i}}$) corresponding to $C$ and $C_{i}$ for $i = 1,\dots,k$ and $\mathbf{u} = (\mathbf{u}_{1},\dots,\mathbf{u}_{k})$ with $\mathbf{u}_{i} = (u_{i1},\dots,u_{id_{i}})$ for $i = 1,\dots,k$. Before showing compliance of the $\Phi$-dependence measures with our stated axioms, we define an artificial normalization $N$ to be a continuous, strictly increasing mapping $N : [0,\infty] \rightarrow \mathbb{I}$ satisfying $N(0) = 0$ and $N(\infty) = 1$. As an example, \cite{Joe1989} suggests $N(t) = \sqrt{1-e^{-2t}}$ in case $\Phi(t) = t \log(t)$, because then the normalized dependence coefficient reduces to $|\rho|$ in case of a bivariate Gaussian distribution with correlation $\rho$.
\newline \\ \noindent 
\textbf{Proposition 1.} \textit{Let $\mathcal{D}_{\Phi}$ be defined by \eqref{eq: phidiv absc} and normalized to \begin{equation*}
    N_{\Phi} \circ \mathcal{D}_{\Phi} = \begin{cases}
    \mathcal{D}_{\Phi} / (\Phi(0) + \Phi^{*}(0)) \hspace{0.2cm} \text{if} \hspace{0.2cm} \Phi(0) + \Phi^{*}(0) < \infty \\
    N(\mathcal{D}_{\Phi}) \hspace{2.4cm} \text{if} \hspace{0.2cm} \Phi(0) + \Phi^{*}(0) = \infty,
    \end{cases}
\end{equation*}
where $N$ is an artificial normalization as explained above. Then, $N_{\Phi} \circ \mathcal{D}_{\Phi}$ satisfies (A1),(A2),(A5),(A6) and (A7). If $\Phi$ is strictly convex at $1$, property (A3) is satisfied, and (A4) holds if $\Phi$ is strictly convex on $(0,\infty)$. Axiom (A8) is fulfilled if we replace $C$ and $C_{n}$ by the existing densities $c$ and $c_{n}$, and when $c$ is uniformly bounded from below and above by a strictly positive constant.}
\begin{proof}
Property (A1) holds because of Fubini's theorem and knowing that permuting the components of a random vector, results in permuting the copula components accordingly. The results stated about (A2) and (A3) follow from Theorem 1 of \cite{Gijbels2023} and applying the normalization. 

For Property (A4), suppose that $(\mathbf{X}_{1},\dots,\mathbf{X}_{k},\mathbf{X}_{k+1})$ has copula density $\widetilde{c}$ with marginal copula density $c_{k+1}$ of $\mathbf{X}_{k+1}$. Put $\widetilde{\mathbf{u}} = (\mathbf{u},\mathbf{u}_{k+1})$ Then,
\begin{equation*}
    \begin{split}
        & \hspace{-1cm} \mathcal{D}_{\Phi}(\mathbf{X}_{1},\dots,\mathbf{X}_{k},\mathbf{X}_{k+1})  \\ & \hspace{1cm} = \int_{\mathbb{I}^{q+d_{k+1}}} \prod_{i=1}^{k+1} c_{i}(\mathbf{u}_{i}) \Phi \Bigg (\frac{\widetilde{c}(\widetilde{\mathbf{u}})}{\prod_{i=1}^{k+1} c_{i}(\mathbf{u}_{i})} \Bigg ) d\widetilde{\mathbf{u}} \\
        & \hspace{1cm} = \int_{\mathbb{I}^{q}} \prod_{i=1}^{k} c_{i}(\mathbf{u}_{i})  \int_{\mathbb{I}^{d_{k+1}}} c_{k+1}(\mathbf{u}_{k+1}) \Phi \Bigg (\frac{\widetilde{c}(\mathbf{u},\mathbf{u}_{k+1})}{c_{k+1}(\mathbf{u}_{k+1}) \prod_{i=1}^{k} c_{i}(\mathbf{u}_{i}) } \Bigg ) d\mathbf{u}_{k+1}d\mathbf{u} \\
        & \hspace{1cm} \geq \int_{\mathbb{I}^{q}} \prod_{i=1}^{k} c_{i}(\mathbf{u}_{i}) \Phi \Bigg (\int_{\mathbb{I}^{d_{k+1}}} \frac{\widetilde{c}(\mathbf{u},\mathbf{u}_{k+1})}{\prod_{i=1}^{k} c_{i}(\mathbf{u}_{i})}d\mathbf{u}_{k+1}  \Bigg ) d\mathbf{u} \\ & \hspace{1cm} = \mathcal{D}_{\Phi}(\mathbf{X}_{1},\dots,\mathbf{X}_{k}),
    \end{split}
\end{equation*}
where we used Jensen's inequality. If $\mathbf{X}_{k+1}$ is independent from $(\mathbf{X}_{1},\dots,\mathbf{X}_{k})$, we have $\widetilde{c}(\mathbf{u},\mathbf{u}_{k+1}) = c(\mathbf{u})  c_{k+1}(\mathbf{u}_{k+1})$ and the equality holds. If $\Phi$ is strictly convex, the equality holds if and only if
$\widetilde{c}(\mathbf{u},\mathbf{u}_{k+1}) = A(\mathbf{u}) c_{k+1}(\mathbf{u}_{k+1}) \prod_{i=1}^{k}c_{i}(\mathbf{u}_{i})$, where $A(\mathbf{u})$ is a function of $\mathbf{u}$ not depending on $\mathbf{u}_{k+1}$ almost surely for almost every $\mathbf{u}$. Integrating this equality with respect to $\mathbf{u}_{k+1}$ gives $A(\mathbf{u}) = c(\mathbf{u})/\prod_{i=1}^{k}c_{i}(\mathbf{u}_{i})$, i.e. $\widetilde{c}(\mathbf{u},\mathbf{u}_{k+1}) = c(\mathbf{u}) c_{k+1}(\mathbf{u}_{k+1})$.

Obviously, by definition, Property (A5) and hence also (A6) are fulfilled. Next, in the context of Property (A7), assume without loss of generality that $X_{11}$ gets transformed to $T_{11}(X_{11})$ for a strictly decreasing transformation $T_{11}$ and let $\widetilde{c}$ be the copula density of $(T_{1}(\mathbf{X}_{1}),\mathbf{X}_{2},\dots, \mathbf{X}_{k})$ with $T_{1}(\mathbf{X}_{1}) = (T_{11}(X_{11}),X_{12},\dots,X_{1d_{1}})$. Then, $\widetilde{c}(\mathbf{u}_{1},\mathbf{u}_{2},\dots,\mathbf{u}_{k}) = c(\widetilde{\mathbf{u}}_{1},\mathbf{u}_{2},\dots,\mathbf{u}_{k})$ with $\widetilde{\mathbf{u}}_{1} = (1-u_{11},u_{12},\dots,u_{1d_{1}})$ and $\widetilde{c}_{1}(\mathbf{u}_{1}) = c_{1}(\widetilde{\mathbf{u}}_{1})$ with $\widetilde{c}_{1}$ the copula density of $T_{1}(\mathbf{X}_{1})$. Hence, 
\begin{equation*}
    \begin{split}
        \text{\scalebox{0.95}{$\mathcal{D}_{\Phi}(T_{1}(\mathbf{X}_{1}),\dots,\mathbf{X}_{k}) = \int_{\mathbb{I}^{q}} c_{1}(\widetilde{\mathbf{u}}_{1}) \prod_{i=2}^{k} c_{i}(\mathbf{u}_{i}) \Phi \Bigg (\frac{c(\widetilde{\mathbf{u}}_{1},\mathbf{u}_{2},\dots,\mathbf{u}_{k})}{c_{1}(\widetilde{\mathbf{u}}_{1}) \prod_{i=2}^{k} c_{i}(\mathbf{u}_{i})} \Bigg ) d\widetilde{\mathbf{u}} = \mathcal{D}_{\Phi}(\mathbf{X}_{1},\dots,\mathbf{X}_{k})$}},
    \end{split}
\end{equation*}
by simply doing a substitution $t_{11} = 1-u_{11}$. 

Finally, given there exist $m, M > 0$ such that $0 < m \leq c(\mathbf{u}) \leq M$ for all $\mathbf{u} \in \mathbb{I}^{q}$, and $c_{n} \to c$ uniformly on $\mathbb{I}^{q}$, some basic analysis implies that
\begin{equation*}
    \left ( \prod_{i=1}^{k} c_{ni} \right ) \Phi \left (\frac{c_{n}}{\prod_{i=1}^{k}c_{ni}} \right ) \to \left ( \prod_{i=1}^{k} c_{i} \right ) \Phi \left (\frac{c}{\prod_{i=1}^{k}c_{i}} \right )
\end{equation*}
uniformly on $\mathbb{I}^{q}$ as $n \to \infty$, where $c_{ni}$ are the marginal copula densities of $c_{n}$ for $i = 1,\dots,k$. Property (A8) in terms of copula densities is then satisfied by the Lebesgue dominated convergence theorem.
\end{proof}
\noindent 
\textit{Remark 1.}
For showing (A3), we assumed that $\Phi$ is strictly convex at $1$. While convexity is usually defined as a global property of a function, we use Definition 1 of local strict convexity of \cite{Liese2006}, i.e. $\Phi$ is strictly convex at $1$ if it is convex and not locally linear at $1$. 
\newline \\ \noindent
\textit{Remark 2.} 
Consistency results for dependence measures based on the copula cdf are typically based on the weak uniform convergence of the empirical copula process. In copula density terms, the stated conditions in Proposition 1 for fulfilling (A8) are rather stringent. Indeed, it is known that many of the common copula families (e.g. normal, Student, Clayton, Gumbel) have densities that explode to infinity near some boundaries points, see e.g. \cite{Omelka2009}. This means that consistency arguments relying on uniform convergence are typically limited to compact subsets of $\mathbb{I}^{q}$. However, there are theoretical properties that favour copula density based dependence measures, see Remark 3.
\newline \\ \noindent
\textit{Remark 3.}
When using a dependence measure that compares the true copula cdf to the one under independence, like the Hoeffding's $\Phi^{2}$ of \cite{Ruppert2010} (with $d_{1} = \cdots = d_{k} = 1$) using the $L_{2}$-distance, the independence characterization Axiom (A3) is still satisfied, but stays rather ambiguous. The reason is that such dependence measures can be made arbitrarily small, while maintaining an exact deterministic relationship (singularity) between all the variables. An explicit proof of this follows from Theorem 3.2.2 of \cite{Nelsen2006}, telling us that we can approximate the independence copula arbitrarily and uniformly closely by copulas (shuffles of Min) exhibiting a perfect deterministic relationship (‘complete dependence' in the sense of \cite{Lancaster1963}). The copula density based $\Phi$-dependence measures are more alert to such singularities. 
\newline 

It is also interesting to think about the meaning of maximal $\Phi$-dependence. Such maximal dependence occurs if there is a certain singularity (and hence the copula density does not exist everywhere). First, we give an overview of popular choices for $\Phi$, the corresponding name, and its maximum value $\Phi(0) + \Phi^{*}(0)$, see Table \ref{tab: table1}. 

Note that all the $\Phi$-functions in Table \ref{tab: table1} are strictly convex on $(0,\infty)$, except for the total variation distance, which is only strictly convex at $1$. The mutual information is a prominent quantity in information theory, see e.g. \cite{Cover2006}. Differential Shannon entropy quantifies the average amount of uncertainty and mutual information equals
the difference between the differential entropy under independence and under the true model. A general family is $\Phi(t) = |t-1|^{\alpha}$ with $\alpha \geq 1$, for which $\Phi(0) = 1$ and $\Phi^{*}(0) = \infty$ for $\alpha > 1$ (for $\alpha = 1$, this is the total variation distance). 
\newpage \noindent 
We refer to \cite{Liese2006} and references therein for further statistical applications, as well as for other choices of $\Phi$. For the Jensen-Shannon divergence, we refer to \cite{Osterreicher2003}.
\begin{table}[h]
\begin{tabularx}{\textwidth}{||Y|Y|Y||}
 \hhline{|===|}
 $\Phi(t)$ & Name & $\Phi(0) + \Phi^{*}(0)$ \\ 
 \Xhline{4\arrayrulewidth}
 $t \log(t)$ & mutual information & $0 + \infty$ \\ 
 $(t-1)^{2}$ & Pearson distance & $1 + \infty$ \\
 $(\sqrt{t}-1)^{2}$ & Hellinger distance & $1+1$ \\ 
 $|t-1|$ & total variation distance & $1+1$ \\
 $-(t+1)\log \left (\frac{t+1}{2} \right ) + t \log(t)$ & Jensen-Shannon distance & $\log(2) + \log(2)$ \\
 \hhline{|===|}
\end{tabularx}
\caption{Common choices for the function $\Phi$.}
\label{tab: table1}
\end{table} 
\newline \\ \noindent 
\textit{3.2. Perfect dependence} 

Next to independence, there is some kind of maximal dependence at the opposite end of the spectrum, to which we will refer as perfect dependence. Perfect dependence is inherent to the dependence measure and occurs if and only if the measure in question reaches its maximum value. 

Two random variables $X_{1},X_{2}$ are often seen as maximally dependent if their copula is the Fr\'{e}chet upper or lower bound copula, that is if $F_{1}(X_{1}) = U_{1} = U_{2} = F_{2}(X_{2})$ almost surely, with $F_{i}$ the cdf of $X_{i}$ for $i = 1,2$, or if $U_{1} = 1 - U_{2}$ almost surely. In this case, concordance measures like Kendall's tau and Spearman's rho are maximal. This however focuses on monotonic dependence, and perfect co- or counter-monotonicity are only particular cases of ‘strict dependence' as in \cite{Renyi1959}, telling that $X_{2} = \Psi_{1}(X_{1})$ for some function $\Psi_{1}$, or $X_{1} = \Psi_{2}(X_{2})$ for some function $\Psi_{2}$ (deterministic predictability of one variable through the other). When $\Psi_{1}$ (or $\Psi_{2}$) is invertible, we get the ‘complete dependence' of \cite{Lancaster1963}. 

\begin{figure}[h!] \centering
\includegraphics[scale = 0.59]{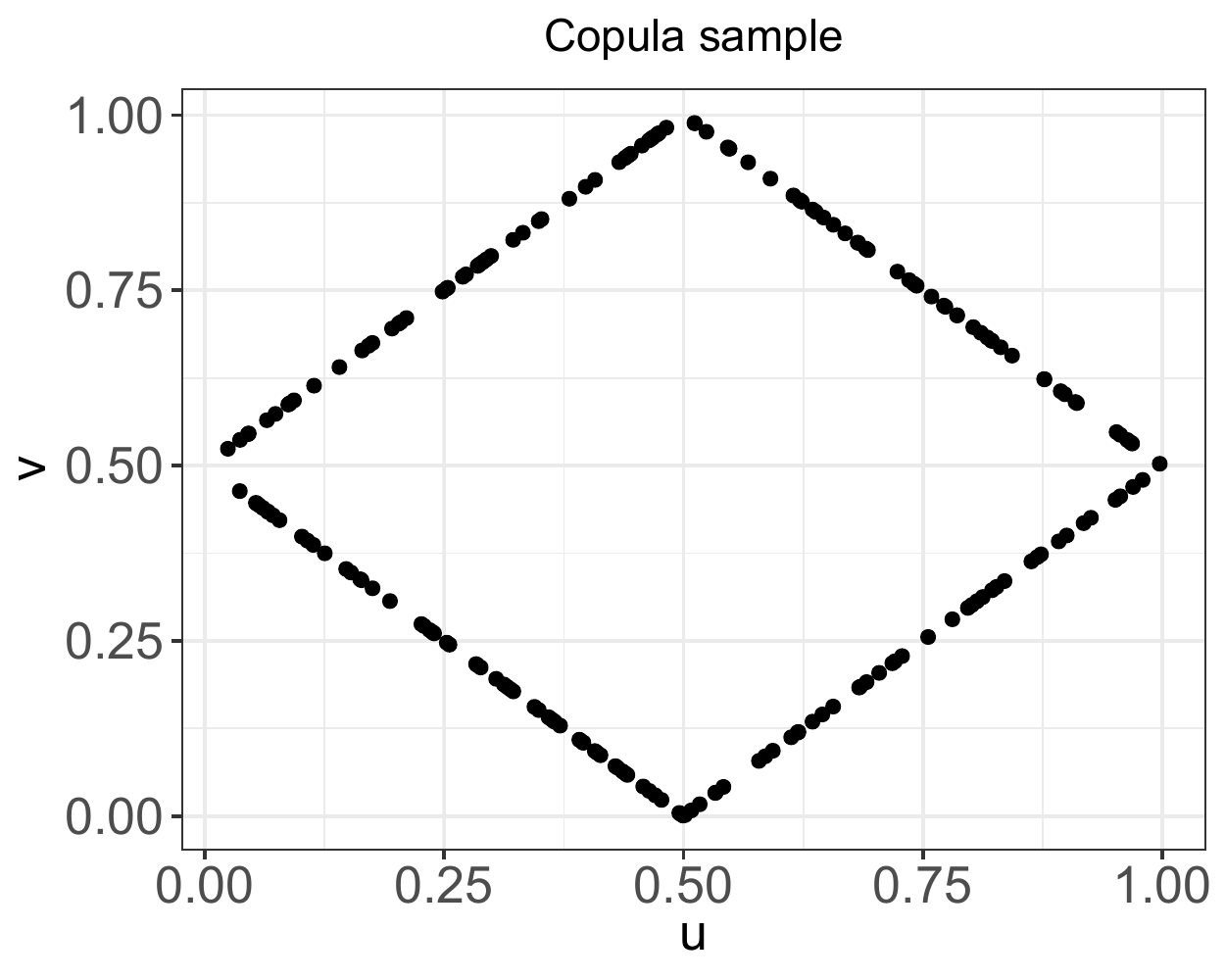}
\includegraphics[scale = 0.59]{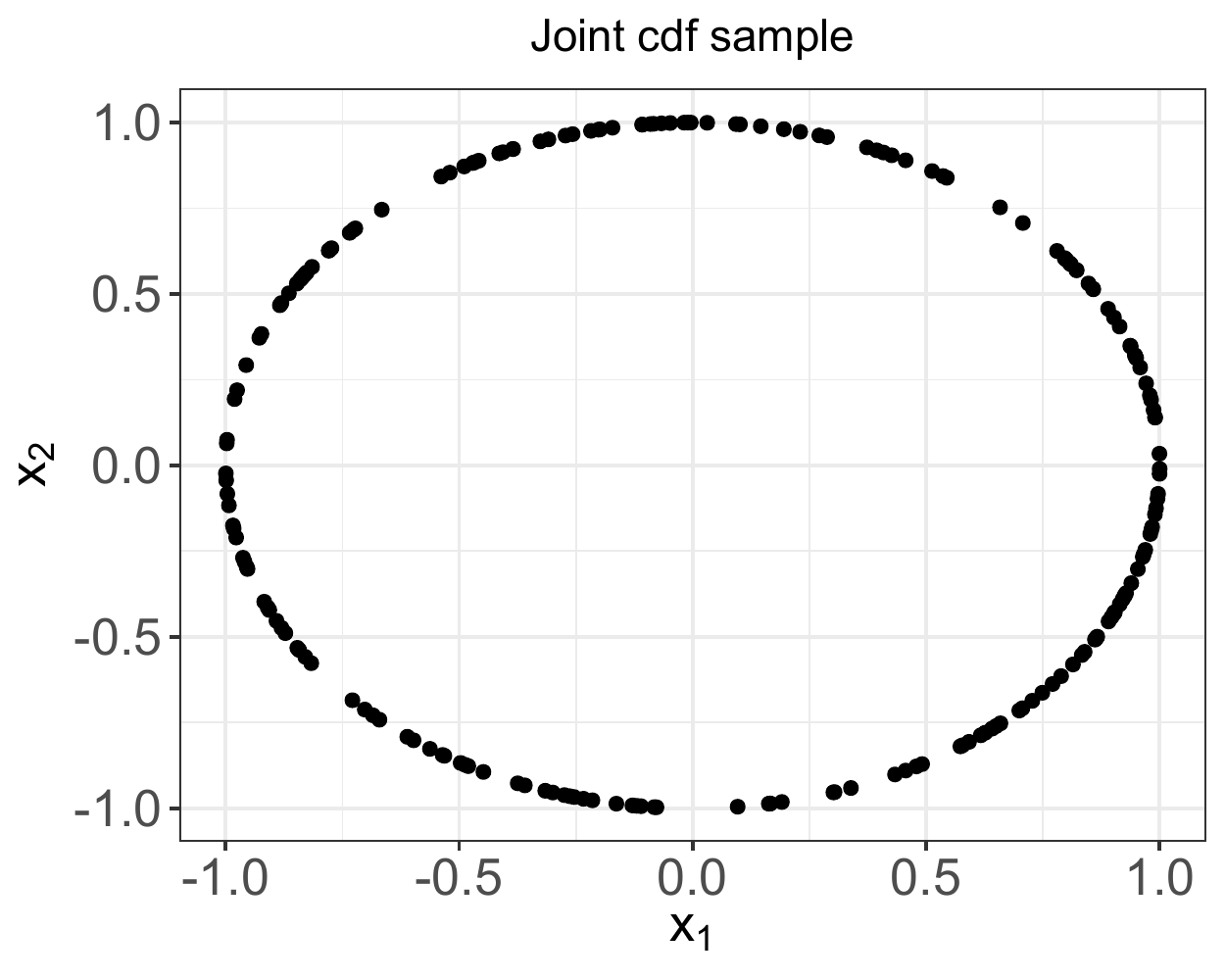}
\caption{Scatter plot (sample size $200$) of bivariate copula \eqref{eq: cop1} (left) and cdf (right) when \eqref{eq: cdf1} are the marginals.}
\label{fig: figure1}
\end{figure}

More general is the ‘pure dependence' of \cite{Geenens2022}, being the existence of a function $\Psi : \mathbb{R} \rightarrow \mathbb{R}^{2}$ such that $(X_{1},X_{2}) = \Psi(U)$ for a certain $U \sim \mathcal{U}[0,1]$, and with $\mathbb{P}_{X_{1} \perp X_{2}}((X_{1},X_{2}) \in \Psi(\mathbb{I})) = 0$, where $\Psi(\mathbb{I})$ is the image of $\mathbb{I}$ under $\Psi$. We provide an example of pure dependence. \newline \\ \noindent
\textbf{Example 1.} Consider $X_{1} \stackrel{d}{=} X_{2}$ (equality in distribution) with
\begin{equation}\label{eq: cdf1}
    \mathbb{P}(X_{1} \leq x_{1}) = \frac{2\arcsin(x_{1})+\pi}{2\pi} \hspace{0.2cm} \text{for} \hspace{0.2cm} x_{1} \in [-1,1],
\end{equation}
and interconnected by the copula 
\begin{equation}\label{eq: cop1}
    \text{\scalebox{0.85}{$C(u_{1},u_{2}) = \begin{cases} \frac{2(u_{1}+u_{2})-1}{4} + \frac{2(u_{1}+u_{2})-3}{4}\mathds{1}\Big (u_{1}+u_{2} > \frac{3}{2}\Big ) & \mbox{if } (u_{1},u_{2}) \in \Big [\frac{1}{2},1 \Big]^{2} \\ \frac{2(u_{1}+u_{2})-1}{4} \mathds{1}\Big (u_{1}+u_{2} > \frac{1}{2}\Big ) & \mbox{if } (u_{1},u_{2}) \in \Big [0,\frac{1}{2} \Big]^{2} \\
    \frac{1+u_{1}-\max\{1-u_{1},\frac{3}{2}-u_{2}\}}{2 } & \mbox{if } (u_{1},u_{2}) \in \Big [0,\frac{1}{2} \Big] \times  \Big [\frac{1}{2},1 \Big] \\
    \frac{\min\{1+u_{1},\frac{3}{2}+u_{2}\}+u_{2}-\frac{3}{2}}{2} & \mbox{if } (u_{1},u_{2}) \in \Big [\frac{1}{2},1 \Big] \times \Big [0,\frac{1}{2} \Big],\end{cases}$}}
\end{equation}
\hspace{-0.3cm}
where $\mathds{1}$ is the indicator function. Figure \ref{fig: figure1} shows scatter plots of $(U_{1},U_{2})$ and $(X_{1},X_{2})$ based on a random sample of size $200$. This is clearly an example of pure dependence, with function $\Psi$ given by $\Psi(U) = $ $(\cos(2\pi U),\sin(2\pi U))$. Note however that this is not an example of strict dependence nor of complete dependence. 
\newline 

We can easily extent this notion of pure dependence to $q$ random variables $X_{1},\dots,X_{q}$, and formulate it in terms of $F_{1}(X_{1}) = U_{1},\dots,F_{q}(X_{q}) = U_{q}$, as the existence of a $\Psi : \mathbb{I}^{q-1} \rightarrow \mathbb{I}^{q}$ such that $(U_{1},\dots,U_{q}) = \Psi(V_{1},\dots,V_{q-1})$ for certain $V_{1},\dots,V_{q-1} \sim \mathcal{U}[0,1]$, and with $\lambda^{q}(\Psi(\mathbb{I}^{q-1})) = 0$. The intuition associated with pure dependence is akin to understanding perfect dependence inherent in the $\Phi$-dependence measures. 

Theorem 1 of \cite{Gijbels2023} suggests that perfect $\Phi$-dependence ($\mathcal{D}_{\Phi}$ maximal) between random vectors $\mathbf{X}_{1},\dots,\mathbf{X}_{k}$ arises when there exists a $B \in \mathcal{B}(\mathbb{I}^{q})$ such that $\mu_{C}(B) = 1$ and $(\mu_{C_{1}} \times \cdots \times \mu_{C_{k}})(B) = 0$ (i.e. $\mu_{C} \perp \mu_{C_{1}} \times \cdots \times \mu_{C_{k}}$, hence $\mu_{C} = \mu_{C}^{\text{s}})$, where $\mathcal{B}(\mathbb{I}^{q})$ is the Borel sigma-algebra on $\mathbb{I}^{q}$. Restricting to $\Phi$-functions that are strictly convex at $1$, and satisfy $\Phi(0) + \Phi^{*}(0) < \infty$, this is the only case in which $\mathcal{D}_{\Phi}$ is maximal (and hence we obtain a characterization of perfect dependence). We illustrate this singularity between random vectors by means of another example.
\newline \\ \noindent
\textbf{Example 2.} Consider $(X_{1},X_{2},X_{3},X_{4})$ having copula
\begin{equation*}
    C(u_{1},u_{2},u_{3},u_{4}) = \exp \bigg (- \Big ( \big ( -\log(\min\{u_{1},u_{3}\}) \big )^{\theta} + \big (-\log(\min\{u_{2},u_{4}\}) \big )^{\theta} \Big )^{\frac{1}{\theta}} \bigg )
\end{equation*}
for $\theta \in [1,\infty)$. Then $(X_{1},X_{2})$ and $(X_{3},X_{4})$ have a Gumbel$(\theta)$ copula and $(X_{1},X_{3})$ and $(X_{2},X_{4})$ the comonotonicity copula. We have $\mu_{C} \perp \lambda^{4}$, since $\mu_{C}$ is concentrated on $\{(u_{1},u_{2},u_{3},u_{4}) \in \mathbb{I}^{4} : u_{1} = u_{3}, u_{2} = u_{4} \}$. Also, $\mathcal{D}_{\Phi}\big ((X_{1},X_{2});(X_{3},X_{4}) \big ) = \infty$, since $\mu_{C} \perp \mu_{C_{1}} \times \mu_{C_{2}}$ ($C_{1}$ the copula of $(X_{1},X_{2})$ and $C_{2}$ the copula of $(X_{3},X_{4}))$.

On the other hand, we have  $\mathcal{D}_{\Phi}\big ((X_{1},X_{3});(X_{2},X_{4}) \big ) < \infty$ since 
$\mu_{C_{1}} \perp \lambda^{2}$ and $\mu_{C_{2}} \perp \lambda^{2}$ (now defining $C_{1}$ as the copula of $(X_{1},X_{3})$ and $C_{2}$ as the copula of $(X_{2},X_{4})$). In fact, this is an example of $\mu_{C} \perp \lambda^{4}$, $\mu_{C_{1}} \times \mu_{C_{2}} \perp \lambda^{4}$, but $\mu_{C} \ll \mu_{C_{1}} \times \mu_{C_{2}}$.
In Example 6, we illustrate that $\mathcal{D}_{\Phi}\big ((X_{1},X_{3});(X_{2},X_{4}) \big )$ (for one particular choice of $\Phi$) is actually nothing more than $\mathcal{D}_{\Phi}$ of a two dimensional Gumbel$(\theta)$ copula, which is pretty intuitive, since $(X_{1},X_{3})$ and $(X_{2},X_{4})$ both have the comonotonicity copula and thus $X_{1}$ and $X_{3}$ can be seen as one, as well as $X_{2}$ and $X_{4}$, and $(X_{1},X_{2})$ (just as $(X_{3},X_{4})$) has a Gumbel$(\theta)$ copula.
\newline 

In short, $\Phi$-dependence is maximal when there is a singularity among random variables belonging to different random vectors, indifferent to plausible singularities between random variables within one and the same random vector. 
\newline \\ \noindent
\textbf{4. A Gaussian copula approach}
\\

In this section, we assume a Gaussian copula model for $\mathbf{X} = (\mathbf{X}_{1},\dots,\mathbf{X}_{k})$. Although a restricted framework, it leads to interesting statistical inference results. In fact, we only assume finite second moments and the existence of the covariance matrix of $\mathbf{X}$, say
\begin{equation*}
\mathbf{\Sigma} = \begin{pmatrix}
\mathbf{\Sigma}_{11} & \mathbf{\Sigma}_{12} & \dots & \mathbf{\Sigma}_{1k} \\
\mathbf{\Sigma}_{12}^{\text{T}} & \mathbf{\Sigma}_{22} & \dots & \mathbf{\Sigma}_{2k} \\
\vdots & \vdots & \ddots & \vdots \\
\mathbf{\Sigma}_{1k}^{\text{T}} & \mathbf{\Sigma}_{2k}^{\text{T}} & \dots & \mathbf{\Sigma}_{kk}
\end{pmatrix} \in \mathbb{R}^{q \times q},
\end{equation*}
where $\mathbf{\Sigma}_{ij} \in \mathbb{R}^{d_{i} \times d_{j}}$ contains the covariances between the components of $\mathbf{X}_{i}$ and $\mathbf{X}_{j}$ for $i \neq j$, and $\mathbf{\Sigma}_{ii} \in \mathbb{R}^{d_{i} \times d_{i}}$ the within covariances of $\mathbf{X}_{i}$ for $i = 1,\dots,k$. If $\mathbf{X}$ has a normal $\mathcal{N}(\mathbf{0}_{q},\mathbf{\Sigma})$ distribution, a straightforward calculation shows that \eqref{eq: phidiv absc} reduces to
\begin{equation}\label{eq: phiN}
    \mathcal{D}_{\Phi}^{\mathcal{N}}(\mathbf{\Sigma}) = \int_{\mathbb{R}^{q}} \frac{\exp \left (-\frac{1}{2} \mathbf{x}^{\text{T}} \mathbf{\Sigma}_{0}^{-1} \mathbf{x} \right )}{(2\pi)^{q/2} \prod_{i=1}^{k} \left |\mathbf{\Sigma}_{ii} \right |^{1/2}} \Phi \left (\frac{\prod_{i=1}^{k}\left |\mathbf{\Sigma}_{ii} \right |^{1/2}}{\left |\mathbf{\Sigma} \right |^{1/2}} \exp \left (-\frac{1}{2} \mathbf{x}^{\text{T}}\mathbf{B}^{-1}\mathbf{x} \right ) \right ) d \mathbf{x},
\end{equation}
with
\begin{equation*}
\mathbf{\Sigma}_{0} = \begin{pmatrix}
\mathbf{\Sigma}_{11} & \mathbf{0}_{12} & \dots & \mathbf{0}_{1k} \\
\mathbf{0}_{12}^{\text{T}} & \mathbf{\Sigma}_{22} & \dots & \mathbf{0}_{2k} \\
\vdots & \vdots & \ddots & \vdots \\
\mathbf{0}_{1k}^{\text{T}} & \mathbf{0}_{2k}^{\text{T}} & \dots & \mathbf{\Sigma}_{kk}
\end{pmatrix}
\end{equation*}
the covariance matrix under mutual independence of $\mathbf{X}_{1},\dots,\mathbf{X}_{k}$, and $\mathbf{B}^{-1} = \mathbf{\Sigma}^{-1} - \mathbf{\Sigma}_{0}^{-1}$. In \eqref{eq: phiN}, we used a superscript $\mathcal{N}$ referring to the Gaussian assumption and emphasize the dependence on merely the covariance matrix $\mathbf{\Sigma}$. For certain specific choices of the function $\Phi$, the above integral allows a closed form solution in terms of the covariance matrix $\mathbf{\Sigma}$. We shall go deeper into the cases $\Phi(t) = t \log(t)$ and $\Phi(t) = (\sqrt{t}-1)^{2}$, and denote the resulting dependence measures as $\mathcal{D}_{t\log(t)}^{\mathcal{N}}$ and $\mathcal{D}_{(\sqrt{t}-1)^{2}}^{\mathcal{N}}$ (for the other choices of $\Phi$ listed in Table \ref{tab: table1}, there is no such elegant closed form solution to \eqref{eq: phiN}). It is a straightforward calculation to show that
\begin{equation}\label{eq: mutN}
        \mathcal{D}_{t\log(t)}^{\mathcal{N}}(\mathbf{\Sigma}) = -\frac{1}{2} \log \left (\frac{\left |\mathbf{\Sigma} \right |}{\prod_{i=1}^{k} \left |\mathbf{\Sigma}_{ii} \right |} \right ), \hspace{0.2cm} 
\end{equation}
and, denoting $\mathbb{I}_{q}$ for the $q \times q$ identity matrix,
\begin{equation}\label{eq: helN}
  \mathcal{D}_{(\sqrt{t}-1)^{2}}^{\mathcal{N}}(\mathbf{\Sigma}) =  2 - 2 \frac{2^{q/2}\left |\mathbf{\Sigma} \right |^{1/4}}{\left |\mathbb{I}_{q}+\mathbf{\Sigma}_{0}^{-1} \mathbf{\Sigma} \right |^{1/2}\prod_{i=1}^{k}\left |\mathbf{\Sigma}_{ii} \right |^{1/4}}.
\end{equation}
It is also easily checked that, if $k = 2$ and $d_{1} = d_{2} = 2$, i.e. in the case of two univariate random variables, the expression \eqref{eq: helN} reduces to what was found in Section 4 of \cite{Geenens2022}.

Notice that \eqref{eq: mutN} and \eqref{eq: helN} are formulated in terms of covariance matrices. However, since our dependence measures are copula-based, they should be scale-invariant. Let $\mathbf{D}_{\mathbf{\Sigma}}$ be the diagonal matrix with the variances on the diagonal. Then, the correlation matrix $\mathbf{R}$ corresponding to $\mathbf{\Sigma}$ satisfies $\mathbf{R} = \mathbf{D}_{\mathbf{\Sigma}}^{-1/2} \mathbf{\Sigma} \mathbf{D}_{\mathbf{\Sigma}}^{-1/2}$, such that $|\mathbf{\Sigma}|$ is just a rescaling of $|\mathbf{R}|$ by multiplying the latter with the product of all variances. From this, we easily see that $\mathcal{D}_{t\log(t)}^{\mathcal{N}}(\mathbf{\Sigma}) = \mathcal{D}_{t\log(t)}^{\mathcal{N}}(\mathbf{R})$. Moreover,
\begin{equation*}
\begin{split}
    \left |\mathbb{I}_{q} + \mathbf{\Sigma}_{0}^{-1} \mathbf{\Sigma} \right | & = \left |\mathbf{D}_{\mathbf{\Sigma}}^{1/2} \right | \left |\mathbb{I}_{q} + \mathbf{\Sigma}_{0}^{-1} \mathbf{\Sigma} \right | \left |\mathbf{D}_{\mathbf{\Sigma}}^{-1/2} \right | = \left |\mathbb{I}_{q} + \mathbf{D}_{\mathbf{\Sigma}}^{1/2} \mathbf{\Sigma}_{0}^{-1}\mathbf{D}_{\mathbf{\Sigma}}^{1/2}\mathbf{R}\mathbf{D}_{\mathbf{\Sigma}}^{1/2}\mathbf{D}_{\mathbf{\Sigma}}^{-1/2} \right | \\ &  = \left |\mathbb{I}_{q} + \mathbf{D}_{\mathbf{\Sigma}}^{1/2} \mathbf{\Sigma}_{0}^{-1} \mathbf{D}_{\mathbf{\Sigma}}^{1/2} \mathbf{R}  \right | = \left |\mathbb{I}_{q} + \mathbf{R}_{0}^{-1} \mathbf{R} \right |,
\end{split}
\end{equation*}
with $\mathbf{R}_{0}$ the correlation matrix corresponding to $\mathbf{\Sigma}_{0}$, implying that $\mathcal{D}_{(\sqrt{t}-1)^{2}}^{\mathcal{N}}(\mathbf{\Sigma}) = \mathcal{D}_{(\sqrt{t}-1)^{2}}^{\mathcal{N}}(\mathbf{R})$. For general $\Phi$, we obtain this relation by performing the substitution $\mathbf{y} = \mathbf{D}_{\boldsymbol{\Sigma}}^{-1/2} \mathbf{x}$ in \eqref{eq: phiN}.

Also, since we do not care about univariate marginals, $\mathbf{R}$ needs not to be the traditional Pearson correlation matrix of a multivariate normal distribution, but can be the margin free Gaussian copula correlation matrix of normal scores, that is 
\begin{equation}\label{eq: Normal scores matrix}
(\mathbf{R}_{im})_{jt} = \rho_{ij,mt} = \text{Corr}\left((\phi^{-1} \circ F_{ij})(X_{ij}),(\phi^{-1} \circ F_{mt})(X_{mt})\right),
\end{equation}
for $i,m = 1,\dots,k$, $j = 1,\dots,d_{i}$ and $t = 1,\dots,d_{m}$, where $\text{Corr}$ stands for the traditional Pearson correlation, and $\phi^{-1}$ for the standard normal quantile function. 
\newline \\ \noindent
\textbf{Statistical inference} \newline 

Based on a sample $\mathbf{X}^{(\ell)} = (\mathbf{X}_{1}^{(\ell)},\dots,\mathbf{X}_{k}^{(\ell)})$ for $\ell = 1,\dots,n$ from $\mathbf{X}$, with $\mathbf{X}_{i}^{(\ell)} = (X_{i1}^{(\ell)},\dots,X_{id_{i}}^{(\ell)})$ for $\ell = 1, \dots, n$ a sample from $\mathbf{X}_{i}$ for $i = 1,\dots,k$, the sample version of \eqref{eq: Normal scores matrix} is known as the matrix of normal scores rank correlation coefficients (\cite{Hajek1967}), 
\begin{equation}\label{eq: est cor matrix}
    \widehat{\mathbf{R}}_{n} = \begin{pmatrix}
\widehat{\mathbf{R}}_{11} & \widehat{\mathbf{R}}_{12} & \cdots & \widehat{\mathbf{R}}_{1k} \\
\widehat{\mathbf{R}}_{12}^{\text{T}} & \widehat{\mathbf{R}}_{22} & \cdots & \widehat{\mathbf{R}}_{2k} \\
\vdots & \vdots & \ddots & \vdots \\
\widehat{\mathbf{R}}_{1k}^{\text{T}} & \widehat{\mathbf{R}}_{2k}^{\text{T}} & \cdots & \widehat{\mathbf{R}}_{kk}
\end{pmatrix} \hspace{0.2cm} \text{with} \hspace{0.2cm} (\widehat{\mathbf{R}}_{im})_{jt} =  \widehat{\rho}_{ij,mt} =  \frac{\frac{1}{n}\sum_{\ell=1}^{n} \widehat{Z}_{ij}^{(\ell)} \widehat{Z}_{mt}^{(\ell)}}{\frac{1}{n}\sum_{\ell=1}^{n}\Big (\phi^{-1} \big (\frac{\ell}{n+1} \big ) \Big )^{2}},
\end{equation}
defined through normal scores 
\begin{equation*}
    \widehat{Z}_{ij}^{(\ell)} = \phi^{-1} \left (\frac{n}{n+1} \widehat{F}_{ij}\Big (X_{ij}^{(\ell)} \Big )\right )
\end{equation*}
obtained from the univariate empirical cdf $\widehat{F}_{ij}(x_{ij}) = \frac{1}{n} \sum_{\ell=1}^{n} \mathds{1} \{X_{ij}^{(\ell)} \leq x_{ij}\}$ for $i = 1,\dots,k$ and $j = 1,\dots,d_{i}$. 
The quantity $\widehat{\rho}_{ij,mt}$ is computed as the conventional Pearson correlation of the bivariate sample of scores $\big ((\widehat{Z}_{ij}^{(1)},\widehat{Z}_{mt}^{(1)}),\dots,(\widehat{Z}_{ij}^{(n)},\widehat{Z}_{mt}^{(n)}) \big )$ and by noting that 
\begin{equation*}
\begin{split}
    \frac{1}{n} \sum_{\ell=1}^{n} \widehat{Z}_{ij}^{(\ell)} & = \frac{1}{n} \sum_{\ell=1}^{n} \phi^{-1} \Big (\frac{\ell}{n+1}\Big ) = 0 \\
    \frac{1}{n} \sum_{\ell=1}^{n} \Big ( \widehat{Z}_{ij}^{(\ell)} \Big )^{2} & = \frac{1}{n} \sum_{\ell=1}^{n} \bigg ( \phi^{-1} \Big (\frac{\ell}{n+1}\Big ) \bigg )^{2},
\end{split}
\end{equation*}
which follows from the fact that $\phi^{-1}(\alpha) = -\phi^{-1}(1-\alpha)$ for $\alpha \in [0,1]$ and $n \widehat{F}_{ij}(X_{ij}^{(\ell)})$ is the rank of $X_{ij}^{(\ell)}$ in the sample $X_{ij}^{(1)},\dots,X_{ij}^{(n)}$. Being rank-based, makes the variance of the normal scores independent of the data. 

A next natural step in estimating $\mathcal{D}_{\Phi}^{\mathcal{N}}(\mathbf{R})$ is to just plug in $\widehat{\mathbf{R}}_{n}$ instead of the unknown matrix $\mathbf{R}$. Let $\mathbb{S}^{q}$ be the set of all $q \times q$ covariance matrices and $\mathbb{S}^{q}_{>} \subset \mathbb{S}^{q}$ the set of all positive definite ones. Let $\varphi$ be the map defined by $\varphi(\mathbf{\Sigma}) = \mathbf{D}_{\mathbf{\Sigma}}^{-1/2}\mathbf{\Sigma}\mathbf{D}_{\mathbf{\Sigma}}^{-1/2}$ for $\mathbf{\Sigma} \in \mathbb{S}^{q}$, and $||\cdot||_{\text{F}}$ the Frobenius matrix norm, i.e. $||\mathbf{\Sigma}||_{\text{F}}^{2} = \text{Tr}(\mathbf{\Sigma}^{\text{T}}\mathbf{\Sigma})$. If we can show the Fr\'echet differentiability of the mapping 
\begin{equation}\label{eq: Fremap}
   (\mathbb{S}^{s},||\cdot||_{\text{F}}) \rightarrow (\mathbb{R},|\cdot|) : \mathbf{\Sigma} \mapsto (\mathcal{D}_{\Phi}^{\mathcal{N}} \circ \varphi)(\mathbf{\Sigma})
\end{equation}
on $\mathbb{S}^{q}_{>}$, then the delta method turns an asymptotic normality result for $\widehat{\mathbf{R}}_{n}$ into an asymptotic normality result for $\mathcal{D}_{\Phi}^{\mathcal{N}}(\widehat{\mathbf{R}}_{n})$. For general $\Phi$, interchanging Fr\'{e}chet differentiation and Lebesgue integration in \eqref{eq: phiN} would be useful, and holds if there exists an integrable function on $\mathbb{R}^{q}$ dominating the Fr\'{e}chet derivative of the integrand of \eqref{eq: phiN} uniformly on $\mathbb{S}^{q}_{>}$. Here, special attention is again devoted to $\Phi(t)$ either $t \log(t)$ or $(\sqrt{t}-1)^{2}$, because we have the more explicit expressions \eqref{eq: mutN} and \eqref{eq: helN}. Theorem 1 states formally the asymptotic normality result for $\mathcal{D}_{\Phi}^{\mathcal{N}}(\widehat{\mathbf{R}}_{n})$.
\newline \\ \noindent 
\textbf{Theorem 1.} \textit{Let $\mathbf{X}$ have a Gaussian copula with correlation matrix $\mathbf{R} \in \mathbb{S}^{q}_{>}$, and let $\widehat{\mathbf{R}}_{n}$ be given by \eqref{eq: est cor matrix}, based on which the plug-in estimator $\widehat{\mathcal{D}}^{\mathcal{N}}_{\Phi,n} = \mathcal{D}_{\Phi}^{\mathcal{N}}(\widehat{\mathbf{R}}_{n})$ is constructed. If the mapping defined in \eqref{eq: Fremap} is Fr\'{e}chet differentiable, then the estimator $\widehat{\mathcal{D}}_{\Phi,n}^{\mathcal{N}}$ is asymptotically normal. Moreover, if differentiation and integration in \eqref{eq: phiN} can be interchanged, we have in addition the expression for the  asymptotic variance: 
\begin{equation*}
    \sqrt{n} \left (\widehat{\mathcal{D}}^{\mathcal{N}}_{\Phi,n} - \mathcal{D}_{\Phi}^{\mathcal{N}}(\mathbf{R}) \right ) \xrightarrow{d} \mathcal{N}(0,\zeta_{\Phi}^{2})
\end{equation*}
with asymptotic variance
\begin{equation*}
    \zeta_{\Phi}^{2} = 2 \text{Tr} \Big ( \big (\mathbf{R}  (\mathbf{M}_{\Phi}-\mathbf{D}_{\mathbf{M}_{\Phi}\mathbf{R}}) \big )^{2} \Big ),
\end{equation*}
where $\mathbf{D}_{\mathbf{M}_{\Phi}\mathbf{R}}$ is the diagonal matrix consisting of the diagonal of $\mathbf{M}_{\Phi}\mathbf{R}$, and with
\begin{equation*}
\begin{split}
    \mathbf{M}_{\Phi} = \frac{1}{2}  \Big ( \mathbf{F}_{1} - \mathcal{D}_{\Phi}^{\mathcal{N}}(\mathbf{R}) \mathbf{R}_{0}^{-1} & - \mathbb{E}_{\mathcal{N}(\mathbf{0},\mathbf{R})} \left [\alpha^{\prime}(\mathbf{X}) \right ] \left (\mathbf{R}^{-1} - \mathbf{R}_{0}^{-1} \right )  \\ & + \mathbf{R}^{-1} \mathbb{E}_{\mathcal{N}(\mathbf{0},\mathbf{R})} \left [\alpha^{\prime}(\mathbf{X}) \mathbf{X}\mathbf{X}^{\text{T}} \right ] \mathbf{R}^{-1} - \mathbf{F}_{2} \Big ),
\end{split}
\end{equation*}
defining $\alpha(\mathbf{X}) = \Phi(k(\mathbf{X}))$, $\alpha^{\prime}(\mathbf{X}) = \Phi^{\prime}(k(\mathbf{X}))$, and
\begin{equation*}
\begin{split}
    \mathbf{F}_{1} & = \mathbf{R}_{0}^{-1} \text{diag} \left (\mathbb{E}_{\mathcal{N}(\mathbf{0},\mathbf{R}_{0})} \left [\alpha(\mathbf{X})\mathbf{X}_{1}\mathbf{X}_{1}^{\text{T}} \right ],\cdots, \mathbb{E}_{\mathcal{N}(\mathbf{0},\mathbf{R}_{0})} \left [\alpha(\mathbf{X})\mathbf{X}_{k}\mathbf{X}_{k}^{\text{T}} \right ] \right ) \mathbf{R}_{0}^{-1}, \\ \mathbf{F}_{2} & = \mathbf{R}_{0}^{-1} \text{diag} \left (\mathbb{E}_{\mathcal{N}(\mathbf{0},\mathbf{R})} \left [\alpha^{\prime}(\mathbf{X})\mathbf{X}_{1}\mathbf{X}_{1}^{\text{T}} \right ],\cdots, \mathbb{E}_{\mathcal{N}(\mathbf{0},\mathbf{R})} \left [\alpha^{\prime}(\mathbf{X})\mathbf{X}_{k}\mathbf{X}_{k}^{\text{T}} \right ] \right ) \mathbf{R}_{0}^{-1},
\end{split}
\end{equation*}
for
\begin{equation*}
     k(\mathbf{X}) = \frac{\prod_{i=1}^{k} \left |\boldsymbol{\Sigma}_{ii} \right |^{1/2}}{\left |\boldsymbol{\Sigma} \right |^{1/2}} \exp \left (-\frac{1}{2} \mathbf{X}^{\text{T}} \left ( \mathbf{R}^{-1} - \mathbf{R}_{0}^{-1} \right ) \mathbf{X} \right ).
\end{equation*}
In particular, one can show that 
\begin{equation*}
    \mathbf{M}_{t \log(t)} = -\frac{1}{2} \left (\mathbf{R}^{-1}- \mathbf{R}_{0}^{-1} \right )
\end{equation*}
and 
\begin{equation*}
    \mathbf{M}_{(\sqrt{t}-1)^{2}} = \frac{2^{q/2}\exp \left (-\frac{1}{2}\mathcal{D}_{t \log(t)}^{\mathcal{N}}(\mathbf{R}) \right )}{ \left |\mathbb{I}_{q} + \mathbf{R}_{0}^{-1} \mathbf{R} \right |^{1/2}} \left [-\frac{1}{2} \left (\mathbf{R}^{-1} - \mathbf{R}_{0}^{-1} \right ) + (\mathbb{I}_{q} + \mathbf{R}_{0}^{-1}\mathbf{R})^{-1}\mathbf{R}_{0}^{-1} - \boldsymbol{\gamma} \right ],
\end{equation*}
defining 
\begin{equation*}
    \text{\scalebox{0.95}{$\boldsymbol{\gamma} = \text{diag} \left (\mathbf{R}_{11}^{-1}\mathbf{J}_{11}\mathbf{R}_{11}^{-1},\dots,\mathbf{R}_{kk}^{-1}\mathbf{J}_{kk}\mathbf{R}_{kk}^{-1} \right ), \text{with} \hspace{0.2cm} \mathbf{R}(\mathbb{I}_{q}+\mathbf{R}_{0}^{-1}\mathbf{R})^{-1} = \begin{pmatrix}
\mathbf{J}_{11} & \mathbf{J}_{12} & \dots & \mathbf{J}_{1k} \\
\mathbf{J}_{21} & \mathbf{J}_{22} & \dots & \mathbf{J}_{2k} \\
\vdots & \vdots & \ddots & \vdots \\
\mathbf{J}_{k1} & \mathbf{J}_{k2} & \dots & \mathbf{J}_{kk}
\end{pmatrix}$}}.
\end{equation*}}
\\

The expression for $\mathbf{M}_{t\log(t)}$ in Theorem 1 was obtained in the proof using the explicit formula \eqref{eq: mutN} that we have for the mutual information (similarly for the Hellinger distance). In the following example, we verify, for $\Phi(t) = t \log(t)$, the general formula for $\mathbf{M}_{\Phi}$ relying on interchanging the Fr\'{e}chet derivative and Lebesgue integral. 
\newline \\ \noindent
\textbf{Example 3.}
In case $\Phi(t) = t \log(t)$, the different terms in the general formula for $\text{Tr}(\mathbf{M}_{\Phi}\mathbf{H})$, for a certain $\mathbf{H} \in \mathbb{S}^{q}_{>}$, can be calculated as follows. First, using (A1) in the proof of Theorem 1, we have
\begin{alignat*}{2} \text{Tr} \left (\mathbf{F}_{1} \mathbf{H} \right )  & = && - \mathbb{E}_{\mathcal{N}(\mathbf{0},\mathbf{R}_{0})} \left [\alpha(\mathbf{X}) \mathbf{X}^{\text{T}} \mathbf{D} \mathbf{X} \right ] \\ & =  && - \log \left (\frac{\prod_{i=1}^{k} \left |\mathbf{R}_{ii} \right |^{1/2}}{\left |\mathbf{R} \right |^{1/2}} \right ) \mathbb{E}_{\mathcal{N}(\mathbf{0},\mathbf{R})} \left [\mathbf{X}^{\text{T}} \mathbf{D} \mathbf{X} \right ] \\ & && + \frac{1}{2} \mathbb{E}_{\mathcal{N}(\mathbf{0},\mathbf{R})} \left [\left (\mathbf{X}^{\text{T}} \mathbf{D} \mathbf{X} \right ) \left (\mathbf{X}^{\text{T}}(\mathbf{R}^{-1}-\mathbf{R}_{0}^{-1}) \mathbf{X} \right ) \right ],
\end{alignat*}
where 
\begin{equation*}
    \mathbf{D} = \left (-\mathbf{R}_{11}^{-1}\mathbf{H}_{11}\mathbf{R}_{11}^{-1},\dots, -\mathbf{R}_{kk}^{-1}\mathbf{H}_{kk}\mathbf{R}_{kk}^{-1}\right ).
\end{equation*}
Next, observe that 
\begin{equation*}
    \text{Tr} \left (\mathcal{D}_{t\log(t)}^{\mathcal{N}}(\mathbf{R})\mathbf{R}_{0}^{-1} \mathbf{H} \right ) = \log \left (\frac{\prod_{i=1}^{k} \left |\mathbf{R}_{ii} \right |^{1/2}}{\left |\mathbf{R} \right |^{1/2}} \right )  \text{Tr} \left (\mathbf{R}_{0}^{-1} \mathbf{H} \right ).
\end{equation*}
Using $\Phi^{\prime}(t) = \log(t) + 1$, it is also straightforward to see that
\begin{equation*}
\begin{split}
    & \hspace{-0.5cm} \text{Tr} \left (\mathbb{E}_{\mathcal{N}(0,\mathbf{R})} \left [\alpha^{\prime}(\mathbf{X}) \right ](\mathbf{R}^{-1} - \mathbf{R}_{0}^{-1}) \mathbf{H} \right )  \\ &  = \left ( \log \left (\frac{\prod_{i=1}^{k} \left |\mathbf{R}_{ii} \right |^{1/2}}{\left |\mathbf{R} \right |^{1/2}} \right ) - \frac{1}{2} \mathbb{E} \left [\mathbf{X}^{\text{T}}(\mathbf{R}^{-1} - \mathbf{R}_{0}^{-1})\mathbf{X} \right ] + 1 \right ) \text{Tr} \left ((\mathbf{R}^{-1} - \mathbf{R}_{0}^{-1}) \mathbf{H} \right ),
\end{split}
\end{equation*}
and 
\begin{equation*}
\begin{split}
    & \hspace{-2.5cm} \text{Tr} \left (\mathbf{R}^{-1} \mathbb{E}_{\mathcal{N}(\mathbf{0},\mathbf{R})} \left [\alpha^{\prime}(\mathbf{X}) \mathbf{X}\mathbf{X}^{\text{T}} \right ] \mathbf{R}^{-1} \mathbf{H} \right ) \\ & = \mathbb{E}_{\mathcal{N}(\mathbf{0},\mathbf{R})} \left [\alpha^{\prime}(\mathbf{X}) \mathbf{X}^{\text{T}} \mathbf{R}^{-1} \mathbf{H} \mathbf{R}^{-1} \mathbf{X} \right ] \\
    & = \log \left (\frac{\prod_{i=1}^{k} \left |\mathbf{R}_{ii} \right |^{1/2}}{\left |\mathbf{R} \right |^{1/2}} \right ) \mathbb{E}_{\mathcal{N}(\mathbf{0},\mathbf{R})} \left [\mathbf{X}^{\text{T}}\mathbf{R}^{-1}\mathbf{H}\mathbf{R}^{-1}\mathbf{X} \right ] \\ &  \hspace{0.5cm} - \frac{1}{2} \mathbb{E}_{\mathcal{N}(\mathbf{0},\mathbf{R})} \left [\left (\mathbf{X}^{\text{T}} \mathbf{R}^{-1} \mathbf{H}\mathbf{R}^{-1} \mathbf{X} \right ) \left (\mathbf{X}^{\text{T}}\left (\mathbf{R}^{-1} - \mathbf{R}_{0}^{-1} \right )\mathbf{X} \right ) \right ] \\ &  \hspace{0.5cm} + \mathbb{E}_{\mathcal{N}(\mathbf{0},\mathbf{R})} \left [\mathbf{X}^{\text{T}}\mathbf{R}^{-1}\mathbf{H}\mathbf{R}^{-1}\mathbf{X} \right ],
\end{split}
\end{equation*}
and finally 
\begin{alignat*}{2}
    \text{Tr} \left (\mathbf{F}_{2} \mathbf{H} \right ) & =  && - \mathbb{E}_{\mathcal{N}(\mathbf{0},\mathbf{R})} \left [\alpha^{\prime}(\mathbf{X}) \mathbf{X}^{\text{T}} \mathbf{D} \mathbf{X} \right ] \\ & = && - \log \left (\frac{\prod_{i=1}^{k} \left |\mathbf{R}_{ii} \right |^{1/2}}{\left |\mathbf{R} \right |^{1/2}} \right ) \mathbb{E}_{\mathcal{N}(\mathbf{0},\mathbf{R})} \left [\mathbf{X}^{\text{T}}\mathbf{D}\mathbf{X} \right ] \\ &  && + \frac{1}{2} \mathbb{E}_{\mathcal{N}(\mathbf{0},\mathbf{R})} \left [\left (\mathbf{X}^{\text{T}} \mathbf{D} \mathbf{X} \right ) \left (\mathbf{X}^{\text{T}}\left (\mathbf{R}^{-1} - \mathbf{R}_{0}^{-1} \right )\mathbf{X} \right ) \right ] - \mathbb{E}_{\mathcal{N}(\mathbf{0},\mathbf{R})} \left [\mathbf{X}^{\text{T}}\mathbf{D}\mathbf{X} \right ].
\end{alignat*}
Combining all this, we obtain that
\begin{equation*}
    \begin{split}
        & \hspace{-1.5cm} 2\text{Tr} \left (\mathbf{M}_{t\log(t)} \mathbf{H} \right )  \\ & = -\log \left (\frac{\prod_{i=1}^{k} \left |\mathbf{R}_{ii} \right |^{1/2}}{\left |\mathbf{R} \right |^{1/2}} \right )  \text{Tr} \left (\mathbf{R}_{0}^{-1} \mathbf{H} \right ) \\ &  \hspace{0.5cm} - \left ( \log \left (\frac{\prod_{i=1}^{k} \left |\mathbf{R}_{ii} \right |^{1/2}}{\left |\mathbf{R} \right |^{1/2}} \right ) - \frac{1}{2} \mathbb{E}_{\mathcal{N}(\mathbf{0},\mathbf{R})} \left [\mathbf{X}^{\text{T}}(\mathbf{R}^{-1} - \mathbf{R}_{0}^{-1})\mathbf{X} \right ] + 1 \right ) \\ & \hspace{7.8cm} \cdot \text{Tr} \left ((\mathbf{R}^{-1} - \mathbf{R}_{0}^{-1}) \mathbf{H} \right ) 
    \end{split}
    \end{equation*}

    \begin{equation*}
    \begin{split}
        &  + \log \left (\frac{\prod_{i=1}^{k} \left |\mathbf{R}_{ii} \right |^{1/2}}{\left |\mathbf{R} \right |^{1/2}} \right ) \mathbb{E}_{\mathcal{N}(\mathbf{0},\mathbf{R})} \left [\mathbf{X}^{\text{T}}\mathbf{R}^{-1}\mathbf{H}\mathbf{R}^{-1}\mathbf{X} \right ] \\ &  - \frac{1}{2} \mathbb{E}_{\mathcal{N}(\mathbf{0},\mathbf{R})} \left [\left (\mathbf{X}^{\text{T}} \mathbf{R}^{-1} \mathbf{H}\mathbf{R}^{-1} \mathbf{X} \right ) \left (\mathbf{X}^{\text{T}}\left (\mathbf{R}^{-1} - \mathbf{R}_{0}^{-1} \right )\mathbf{X} \right ) \right ] \\ &  + \mathbb{E}_{\mathcal{N}(\mathbf{0},\mathbf{R})} \left [\mathbf{X}^{\text{T}}\mathbf{R}^{-1}\mathbf{H}\mathbf{R}^{-1}\mathbf{X} \right ] + \mathbb{E}_{\mathcal{N}(\mathbf{0},\mathbf{R})} \left [\mathbf{X}^{\text{T}}\mathbf{D}\mathbf{X} \right ] \\
        & \hspace{-0.5cm}  = - \log \left (\frac{\prod_{i=1}^{k} \left |\mathbf{R}_{ii} \right |^{1/2}}{\left |\mathbf{R} \right |^{1/2}} \right ) \text{Tr} \left (\mathbf{R}^{-1}\mathbf{H} \right ) \\ &  + \frac{1}{2} \mathbb{E}_{\mathcal{N}(\mathbf{0},\mathbf{R})} \left [\mathbf{X}^{\text{T}}(\mathbf{R}^{-1} - \mathbf{R}_{0}^{-1})\mathbf{X} \right ] \text{Tr}\left ((\mathbf{R}^{-1}-\mathbf{R}_{0}^{-1})\mathbf{H} \right )
        \\ &  - \text{Tr}\left ((\mathbf{R}^{-1}-\mathbf{R}_{0}^{-1})\mathbf{H} \right ) + \log \left (\frac{\prod_{i=1}^{k} \left |\mathbf{R}_{ii} \right |^{1/2}}{\left |\mathbf{R} \right |^{1/2}} \right ) \mathbb{E}_{\mathcal{N}(\mathbf{0},\mathbf{R})} \left [\mathbf{X}^{\text{T}}\mathbf{R}^{-1}\mathbf{H}\mathbf{R}^{-1}\mathbf{X} \right ] \\ &  - \frac{1}{2} \mathbb{E}_{\mathcal{N}(\mathbf{0},\mathbf{R})} \left [\left (\mathbf{X}^{\text{T}} \mathbf{R}^{-1} \mathbf{H}\mathbf{R}^{-1} \mathbf{X} \right ) \left (\mathbf{X}^{\text{T}}\left (\mathbf{R}^{-1} - \mathbf{R}_{0}^{-1} \right )\mathbf{X} \right ) \right ] \\ &  + \mathbb{E}_{\mathcal{N}(\mathbf{0},\mathbf{R})} \left [\mathbf{X}^{\text{T}}\mathbf{R}^{-1}\mathbf{H}\mathbf{R}^{-1}\mathbf{X} \right ] + \mathbb{E}_{\mathcal{N}(\mathbf{0},\mathbf{R})} \left [\mathbf{X}^{\text{T}}\mathbf{D}\mathbf{X} \right ]
        \\ & \hspace{-0.5cm}  = \frac{1}{2} \text{Tr} \left ((\mathbf{R}^{-1}-\mathbf{R}_{0}^{-1})\mathbf{R} \right ) \text{Tr} \left ((\mathbf{R}^{-1} - \mathbf{R}_{0}^{-1})\mathbf{H} \right ) - \text{Tr}\left ((\mathbf{R}^{-1} - \mathbf{R}_{0}^{-1})\mathbf{H} \right ) \\ &  - \text{Tr} \left (\mathbf{R}^{-1}\mathbf{H}(\mathbf{R}^{-1}-\mathbf{R}_{0}^{-1})\mathbf{R} \right ) - \frac{1}{2} \text{Tr} \left (\mathbf{R}^{-1}\mathbf{H} \right ) \text{Tr} \left ((\mathbf{R}^{-1}-\mathbf{R}_{0}^{-1})\mathbf{R} \right ) \\ & + \text{Tr} \left (\mathbf{R}^{-1}\mathbf{H} \right )  + \text{Tr} \left (\mathbf{D}\mathbf{R} \right ) \\
        & \hspace{-0.5cm} =  - \text{Tr}\left ((\mathbf{R}^{-1} - \mathbf{R}_{0}^{-1})\mathbf{H} \right ) - \text{Tr} \left (\mathbf{R}^{-1}\mathbf{H}(\mathbf{R}^{-1}-\mathbf{R}_{0}^{-1})\mathbf{R} \right ) + \text{Tr} \left (\mathbf{R}^{-1}\mathbf{H} \right ) + \text{Tr} \left (\mathbf{D}\mathbf{R} \right ) \\ & \hspace{-0.5cm} = - \text{Tr}\left ((\mathbf{R}^{-1} - \mathbf{R}_{0}^{-1})\mathbf{H} \right ),
    \end{split}
\end{equation*}
where we used the following properties of quadratic forms 
\begin{equation*}
    \begin{split}
        \mathbb{E}_{\mathcal{N}(\mathbf{0},\mathbf{R})} \left [\mathbf{X}^{\text{T}} \mathbf{A} \mathbf{X} \right ] & = \text{Tr} \left (\mathbf{A}\mathbf{R} \right ) \\
         \mathbb{E}_{\mathcal{N}(\mathbf{0},\mathbf{R})} \left [\left (\mathbf{X}^{\text{T}}\mathbf{A}_{1}\mathbf{X} \right ) \left (\mathbf{X}^{\text{T}}\mathbf{A}_{2}\mathbf{X} \right ) \right ] & = \text{Tr} \left (\mathbf{A}_{1}\mathbf{R}\left (\mathbf{A}_{2} + \mathbf{A}_{2}^{\text{T}} \right )\mathbf{R} \right ) + \text{Tr} \left (\mathbf{A}_{1}\mathbf{R} \right ) \text{Tr} \left (\mathbf{A}_{2} \mathbf{R} \right ),
    \end{split}
\end{equation*}
for certain compatible matrices $\mathbf{A},\mathbf{A}_{1},\mathbf{A}_{2}$, the cyclic trace property, and the fact that $\text{Tr}((\mathbf{R}^{-1}-\mathbf{R}_{0}^{-1})\mathbf{R}) = 0$ and $\text{Tr} \left (\mathbf{D}\mathbf{R} \right ) = -\text{Tr}(\mathbf{R}_{0}^{-1}\mathbf{H})$. Hence, we have verified that the general formula for $\mathbf{M}_{\Phi}$ also brings us to 
\begin{equation*}
    \mathbf{M}_{t\log(t)} = -\frac{1}{2} \left (\mathbf{R}^{-1}-\mathbf{R}_{0}^{-1} \right ).
\end{equation*}

We end this section by looking at a specific four dimensional Gaussian copula family. \newline \\ \noindent
\textbf{Example 4.} Consider a four dimensional random vector $(X_{1},X_{2},X_{3},X_{4})$ having a Gaussian copula with correlation matrix
\begin{equation*}
    \mathbf{R} = \begin{pmatrix}
    1 & \rho_{1} & \rho_{2} & \rho_{2} \\
    \rho_{1} & 1 & \rho_{2} & \rho_{2} \\
    \rho_{2} & \rho_{2} & 1 & \rho_{1} \\
    \rho_{2} & \rho_{2} & \rho_{1} & 1
    \end{pmatrix}, \hspace{0.3cm} \text{where} \hspace{0.3cm} \rho_{1} \geq 2 |\rho_{2}| - 1.
\end{equation*}
Then one can check that
\small
\begin{equation}\label{eq: xmplmutN}
        \mathcal{D}^{\mathcal{N}}_{t \log(t)}\left ( (X_{1},X_{2}); (X_{3},X_{4}) \right ) = -\frac{1}{2} \log \left ( \frac{(\rho_{1}-2\rho_{2}+1)(\rho_{1}+2\rho_{2}+1)}{(1+\rho_{1})^{2}} \right ),
\end{equation}
and 
\begin{equation}\label{eq: xmplhelN}
    \mathcal{D}^{\mathcal{N}}_{(\sqrt{t}-1)^{2}} \left ( (X_{1},X_{2});(X_{3},X_{4}) \right )/2 = 1 -  \frac{(1+\rho_{1})^{1/2}\left ( (\rho_{1}-2\rho_{2}+1)(\rho_{1}+2\rho_{2}+1) \right)^{1/4}}{\left ( (1+\rho_{1}-\rho_{2})(1+\rho_{1}+\rho_{2}) \right )^{1/2}}.
\end{equation}
\normalsize
\begin{figure}[h!] 
\includegraphics[scale = 0.75]{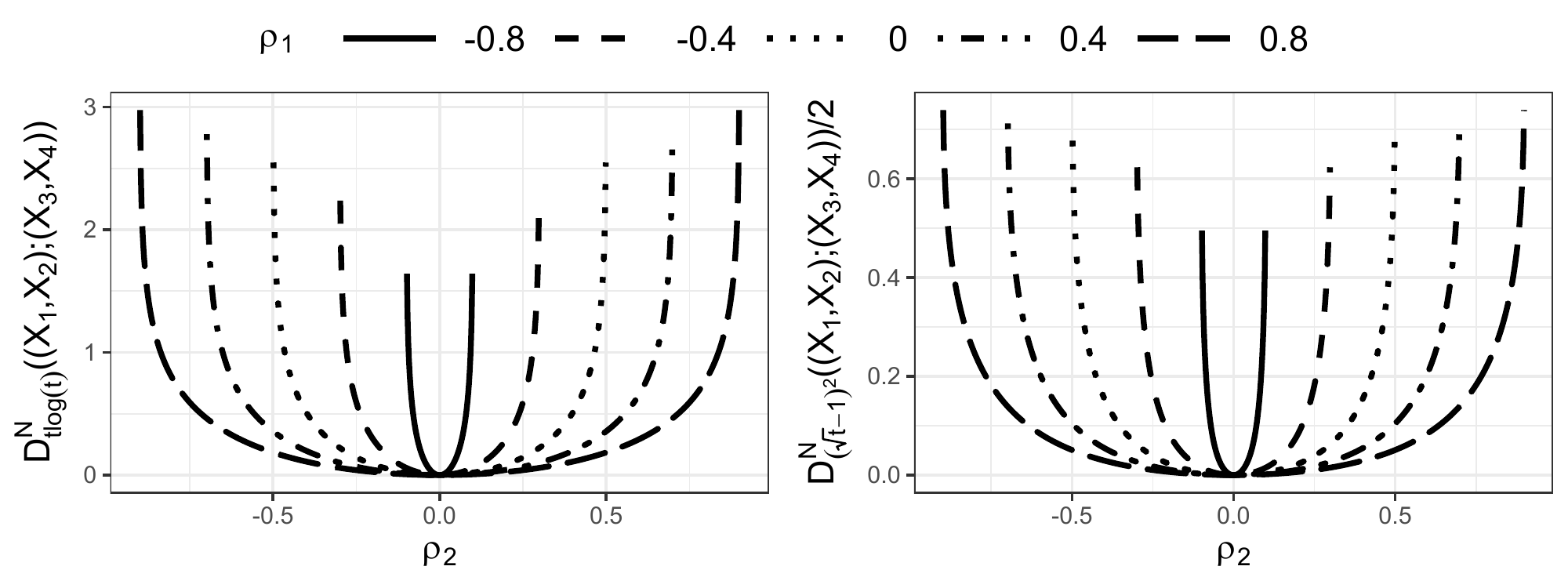}
\caption{Mutual information \eqref{eq: xmplmutN} (left) and half Hellinger distance \eqref{eq: xmplhelN} (right) as a function of $\rho_{2}$ for different values of $\rho_{1}$.}
\label{fig: xmplN}
\end{figure} \newline \noindent 
In \eqref{eq: xmplhelN}, we normalized the Hellinger distance, guaranteeing a dependence measurement in $\mathbb{I}$ (recall Proposition 1). For the mutual information, an artificial normalization is required, but we do not do not implement it here. Figure \ref{fig: xmplN} shows how \eqref{eq: xmplmutN} and \eqref{eq: xmplhelN} depend on $\rho_{2}$ for different values of $\rho_{1}$.
Some observations are:
\begin{itemize}
    \item $\mathcal{D}^{\mathcal{N}}_{t \log(t)} = \mathcal{D}^{\mathcal{N}}_{(\sqrt{t}-1)^{2}} =  0$ iff $\rho_{2} = 0$.
    \item For $\rho_{1} \to -1$ (singularity of $\mu_{C_{1}}$ and $\mu_{C_{2}}$ w.r.t. $\lambda^{2}$), we must have $\rho_{2} \to 0$ and see that $\mathcal{D}^{\mathcal{N}}_{t \log(t)} \to 0$ and $\mathcal{D}^{\mathcal{N}}_{(\sqrt{t}-1)^{2}} \to 0$.
    \item If $\rho_{1} = 2|\rho_{2}|-1$ (absolute continuity of $\mu_{C_{1}}$ and $\mu_{C_{2}}$ w.r.t. $\lambda^{2}$, but singularity of $\mu_{C}$ w.r.t. $\lambda^{4}$), we get $\mathcal{D}^{\mathcal{N}}_{t \log(t)} = \infty$ and $\mathcal{D}^{\mathcal{N}}_{(\sqrt{t}-1)^{2}} = 1$ .
    \item For $\rho_{1} \to 1$ (singularity of $\mu_{C_{1}}$ and $\mu_{C_{2}}$ w.r.t. $\lambda^{2}$), $\mathcal{D}^{\mathcal{N}}_{t \log(t)} \to -1/2 \log(1-\rho_{2}^{2})$, and $\mathcal{D}^{\mathcal{N}}_{(\sqrt{t}-1)^{2}}/2 \to 1-(2(1-\rho_{2}^{2})^{1/4})/(4-\rho_{2}^{2})^{1/2}$, being the mutual information and half Hellinger distance of a bivariate Gaussian copula with correlation $\rho_{2}$ and is maximal iff $|\rho_{2}| = 1$.
\end{itemize} 
Note that the principal components of $(Z_{1},Z_{2}) = ((\phi^{-1} \circ F_{1})(X_{1}),(\phi^{-1} \circ F_{2})(X_{2}))$ are 
\begin{equation*}
    PC_{1} = \frac{1}{\sqrt{2}}(Z_{2} - Z_{1}) \hspace{0.2cm} \text{and} \hspace{0.2cm} PC_{2} = \frac{1}{\sqrt{2}}(Z_{1} + Z_{2}) 
\end{equation*}
and similarly, those of $(Z_{3},Z_{4}) = ((\phi^{-1} \circ F_{3})(X_{3}),(\phi^{-1} \circ F_{4})(X_{4}))$ are
\begin{equation*}
    PC_{3} = \frac{1}{\sqrt{2}}(Z_{4} - Z_{3}) \hspace{0.2cm} \text{and} \hspace{0.2cm} PC_{4} = \frac{1}{\sqrt{2}}(Z_{3} + Z_{4}).
\end{equation*}
Moreover, $\text{Corr}(PC_{1},PC_{3}) = \text{Corr}(PC_{2},PC_{3}) = \text{Corr}(PC_{1},PC_{4}) = 0$, but 
\begin{equation*}
    \text{Corr}(PC_{2},PC_{4}) = \frac{2 \rho_{2}}{1+\rho_{1}}.
\end{equation*}
And so, we see that if $\rho_{1} = 2|\rho_{2}| - 1$, $|\text{Corr}(PC_{2},PC_{4})| = 1$, i.e. the principal components $PC_{2}$ and $PC_{4}$ are perfectly correlated. This means that the four dimensional random vector $(X_{1},X_{2},X_{3},X_{4})$ is propagated in a three dimensional subspace (scatterplot of $(Z_{1},Z_{2},Z_{3},Z_{4})$ constitutes a hyperplane), resulting in the singularity of $\mu_{C}$ with respect to $\mu_{C_{1}} \times \mu_{C_{2}}$ and explaining maximal dependence. 
\begin{figure}[h!] 
\includegraphics[scale = 0.75]{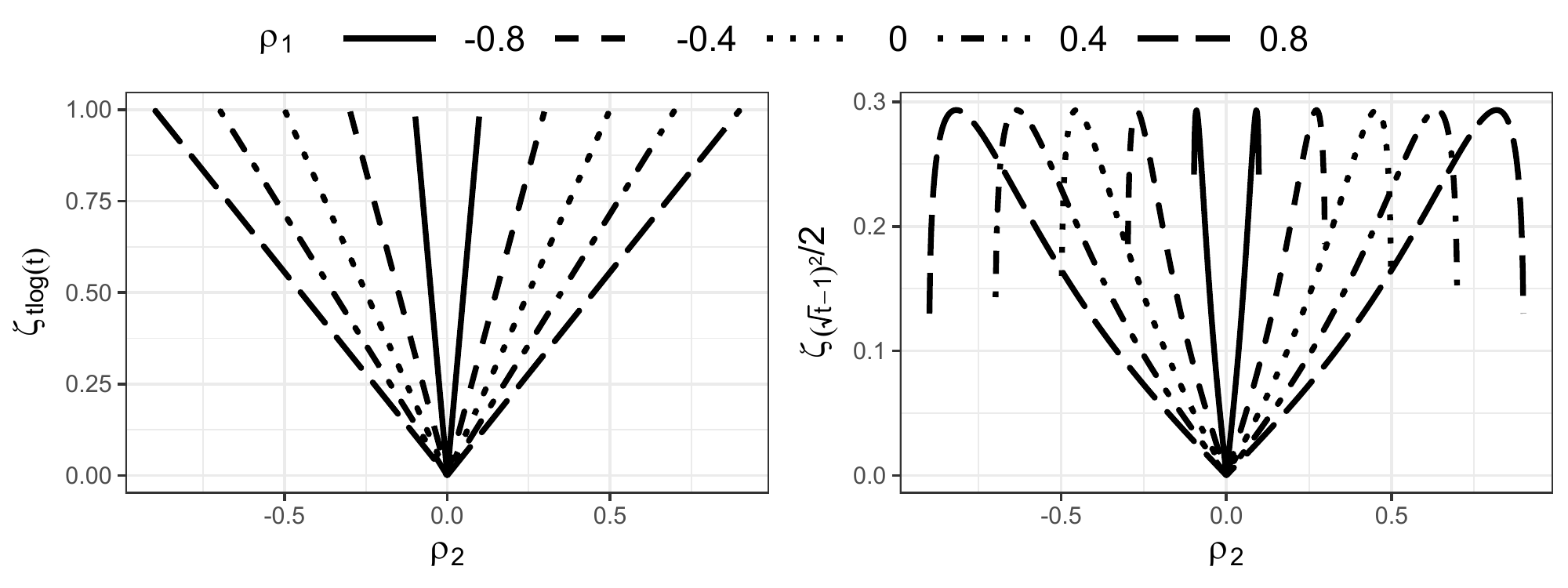}
\caption{Asymptotic standard deviation of mutual information (left) and half Hellinger distance (right) as a function of $\rho_{2}$ for different values of $\rho_{1}$.}
\label{fig: xmplNASD}
\end{figure} 

Figure \ref{fig: xmplNASD} shows the asymptotic standard deviation $\zeta_{t \log(t)}$ of the mutual information and $\zeta_{(\sqrt{t}-1)^{2}}/2$ of the half Hellinger distance (as in Theorem 1) as a function of $\rho_{2}$ for different values of $\rho_{1}$, which can be calculated as
\begin{equation*}
    \begin{split}
        \zeta_{t \log(t)} & = \frac{2\left |\rho_{2} \right |}{1+\rho_{1}}\\
        \zeta_{(\sqrt{t}-1)^{2}}/2 & = \frac{\left ((\rho_{1}-2\rho_{2}+1)(\rho_{1}+2\rho_{2}+1) \right )^{1/4} \left (2\rho_{2}^{2} + (1+\rho_{1})^{2} \right ) \left |\rho_{2} \right |}{2(1+\rho_{1})^{1/2}(\rho_{1}-\rho_{2}+1)^{3/2}(\rho_{1}+\rho_{2}+1)^{3/2}}.
    \end{split}
\end{equation*}
In general, we see that higher degrees of dependence come with higher asymptotic variance. For the Hellinger distance, however, the asymptotic variance goes down to zero when we get close to the previously discussed singularity, i.e. when $\rho_{2}$ gets close to satisfying $|\rho_{2}| = (\rho_{1}+1)/2$. For instance, if $\rho_{1} = 0$, the asymptotic variance is maximal at $\left |\rho_{2} \right | \approx 0.45427$, after which it converges to zero for $\left |\rho_{2} \right | \to 0.5$. From its mathematical expression above, and its general form in Theorem 1, we see the factor $\exp((-1/2)\mathcal{D}_{t \log(t)}^{\mathcal{N}}) \to 0$ in case of singularity makes the asymptotic variance of the Hellinger distance tend to zero.
\newline \\ \noindent
\textbf{5. A maximum likelihood approach} 
\\

A Gaussian copula model is restricted to monotone dependence structures in terms of correlations. However, in many cases other relationships are present. Think for instance of comovements in the tails of stock returns. Thanks to Proposition 1 (in particular the fulfillment of Axiom (A3)), we know that $\Phi$-dependence measures are able to capture such associations as well. In practice, a sufficiently flexible estimation methodology is required. If one is willing to, for example, assume a Clayton copula model for two random variables $X_{1}$ and $X_{2}$, we can use this assumption to estimate $\mathcal{D}_{\Phi}$ and hence, next to a monotone relationship, possible lower tail dependence is guaranteed to be incorporated as well. 
\newline \\ \noindent
\textit{5.1. Maximum likelihood-based inference}

In general, in this section, we assume a specified parametric model for the copula density of $\mathbf{X} = (\mathbf{X}_{1},\dots,\mathbf{X}_{k})$, say $\mathcal{P}_{C} = \{ c(\hspace{0.1cm} \cdot \hspace{0.1cm} ; \boldsymbol{\theta}_{C}) : \boldsymbol{\theta}_{C} \in \boldsymbol{\Theta}_{C} \subset \mathbb{R}^{D}\}$, where $D = \text{dim}(\boldsymbol{\Theta}_{C})$. Notice that a Gaussian copula is also a parametric copula of this form when we stack all upper triangle elements of the correlation matrix in a vector. The univariate marginals $F_{ij}$ for $i = 1,\dots,k$ and $j = 1,\dots,d_{i}$ are approached in either a parametric, or non-parametric way, leading to two different estimation procedures.
\newline \\ \noindent
\underline{Case 1: parametric marginals.}
A first option is to also assume a parametric statistical density model for $F_{ij}$, say $\mathcal{P}_{ij} = \{f_{ij}(\hspace{0.1cm} \cdot \hspace{0.1cm}; \boldsymbol{\theta}_{ij}) : \boldsymbol{\theta}_{ij} \in \mathbf{\Theta}_{ij} \}$. Based on a sample $\mathbf{X}^{(1)}, \dots, \mathbf{X}^{(n)}$ from $\mathbf{X}$, the full log-likelihood is 
\begin{align}
        \ell(\boldsymbol{\theta}_{11},\dots,\boldsymbol{\theta}_{kd_{k}}, \boldsymbol{\theta}_{C}) & = \sum_{\ell=1}^{n} \log \left [ c \left (F_{11}\left (X_{11}^{(\ell)}; \boldsymbol{\theta}_{11} \right ), \dots, F_{kd_{k}}\left (X_{kd_{k}}^{(\ell)};\boldsymbol{\theta}_{kd_{k}} \right ) ; \boldsymbol{\theta}_{C} \right ) \right ] \notag \\
        & \hspace{4.3cm} + \sum_{i=1}^{k} \sum_{j=1}^{d_{i}} \sum_{\ell=1}^{n} \log \left [ f_{ij}\left (X_{ij}^{(\ell)};\boldsymbol{\theta}_{ij} \right ) \right ]. \label{eq: fullll}
\end{align}
Let $\boldsymbol{\eta} = (\boldsymbol{\theta}_{11},\dots,\boldsymbol{\theta}_{kd_{k}},\boldsymbol{\theta}_{C})$ and $\widehat{\boldsymbol{\eta}}_{n}^{\text{MLE}} = (\widehat{\boldsymbol{\theta}}_{11,n}^{\text{MLE}},\dots,\widehat{\boldsymbol{\theta}}_{kd_{k},n}^{\text{MLE}},\widehat{\boldsymbol{\theta}}_{C,n}^{\text{MLE}})$ be the MLE of $\boldsymbol{\eta}$, obtained by maximizing \eqref{eq: fullll}. It is a very well-known result that under certain regularity conditions (true $\boldsymbol{\eta}$ lies in an open subset in which the log-density admits third derivatives w.r.t. the parameters that are uniformly bounded by an integrable function, and interchanging of derivatives and integral for score equations, see e.g. Theorem 4.1 on page 429 in \cite{Lehmann1983}), one has
that $\widehat{\boldsymbol{\eta}}_{n}^{\text{MLE}}$ is an asymptotically normal estimator.
If $\boldsymbol{\theta}_{C}$ is $m$-dimensional and $\boldsymbol{\theta}_{ij}$ all one dimensional, \eqref{eq: fullll} is a $q+m$-dimensional optimization problem. In case of a multivariate normal distribution, one has for $\boldsymbol{\theta}_{C}$ the vector of upper triangle elements of the correlation matrix of the normal copula, and $\boldsymbol{\theta}_{ij} = (\mu_{ij},\sigma_{ij}^{2})$ the mean and variance of $X_{ij}$. 

Given estimators $\widehat{\boldsymbol{\theta}}_{ij,n}^{\text{IFM}}$ of $\boldsymbol{\theta}_{ij}$ (e.g. also based on MLE using the sample $X_{ij}^{(1)},\dots,X_{ij}^{(n)}$), the pseudo log-likelihood for estimating $\boldsymbol{\theta}_{C}$ is 
\begin{equation*}
    \ell(\boldsymbol{\theta}_{C}) = \sum_{\ell=1}^{n} \log \left [ c \left (F_{11} \left (X_{11}^{(\ell)}; \widehat{\boldsymbol{\theta}}_{11,n}^{\text{IFM}} \right ), \dots, F_{kd_{k}} \left (X_{kd_{k}}^{(\ell)};\widehat{\boldsymbol{\theta}}_{kd_{k},n}^{\text{IFM}} \right ) ; \boldsymbol{\theta}_{C} \right ) \right ].
\end{equation*}
If $\boldsymbol{\theta}_{C}$ is $m$-dimensional and $\boldsymbol{\theta}_{ij}$ all one dimensional, this is
an $m$-dimensional optimization problem, after having done $q$ one dimensional optimization problems. This method is known as the inference functions for margins (IFM) method, and extensively studied in Section 10.1 of \cite{Joe1997}. Asymptotic normality is known to hold under the same regularity conditions as for the MLE.
For a multivariate normal distribution, it is known that $\widehat{\boldsymbol{\eta}}_{n}^{\text{IFM}} = \widehat{\boldsymbol{\eta}}_{n}^{\text{MLE}}$.
\newline \\ \noindent
\underline{Case 2: non-parametric marginals.} If no appropriate parametric models can be proposed for the marginals, an option consists of estimating the marginals via the empirical cdf 
\begin{equation*}
    \widehat{F}_{ij}(x_{ij}) = \frac{1}{n+1} \sum_{\ell=1}^{n} \mathds{1} \left \{X_{ij}^{(\ell)} \leq x_{ij}\right \},
\end{equation*}
and the copula parameter $\boldsymbol{\theta}_{C}$ through maximizing
\begin{equation}\label{eq: npll}
    \ell \left (\boldsymbol{\theta}_{C} \right ) = \sum_{\ell=1}^{n} \log \left [c \left (\widehat{F}_{11}\left (X_{11}^{(\ell)} \right ),\dots,\widehat{F}_{kd_{k}} \left (X_{kd_{k}}^{(\ell)} \right ) ; \boldsymbol{\theta}_{C} \right ) \right ],
\end{equation}
resulting in an estimator $\widehat{\boldsymbol{\theta}}_{C,n}^{\text{NP}}$ for $\boldsymbol{\theta}_{C}$. Again under the same regularity conditions, it is shown in \cite{Genest1995} that $\widehat{\boldsymbol{\theta}}^{\text{NP}}_{C,n}$ is asymptotically normal. If $c$ is the Gaussian copula density, there is no known expression for the maximizer of \eqref{eq: npll} over all correlation matrices, and numerical maximization is often used (see Section 5.5.3 in \cite{McNeil2005} for more details), being quite unfeasible in high dimensions. Of course, we have the results from Section 4 dealing with the Gaussian copula case.
\\

In all the above cases, we have an asymptotically normal estimator  $\widehat{\boldsymbol{\theta}}_{n} = \widehat{\boldsymbol{\theta}}_{C,n}$ for $\boldsymbol{\theta}_{C}$. A natural estimator for the copula density $c(\hspace{0.1cm} \cdot \hspace{0.1cm};\boldsymbol{\theta}_{C})$ is then $\widehat{c}(\hspace{0.1cm} \cdot \hspace{0.1cm}) = c(\hspace{0.1cm} \cdot \hspace{0.1cm};\widehat{\boldsymbol{\theta}}_{n})$. If $c$ is unbounded, we will generally not have uniform consistency of $\widehat{c}$ (recall Remark 2). Nevertheless, Fr\'{e}chet differentiability of the mapping
\begin{equation}\label{eq: parmap}
    \boldsymbol{\Theta}_{C} \rightarrow \mathbb{R} : \boldsymbol{\theta}_{C} \mapsto \mathcal{D}_{\Phi} \left (c(\hspace{0.1cm} \cdot \hspace{0.1cm} ; \hspace{0.1cm} \boldsymbol{\theta}_{C}) \right ) = \int_{\mathbb{I}^{q}} \prod_{i=1}^{k} c_{i} \left (\mathbf{u}_{i};\boldsymbol{\theta}_{C} \right ) \Phi \left (\frac{c \left (\mathbf{u} ; \boldsymbol{\theta}_{C} \right )}{\prod_{i=1}^{k} c_{i} \left (\mathbf{u}_{i};\boldsymbol{\theta}_{C} \right )} \right ) d\mathbf{u}
\end{equation}
suffices to turn the asymptotic normality result of $\widehat{\boldsymbol{\theta}}_{n}$ into an asymptotic normality result for the plug-in estimator $\mathcal{D}_{\Phi}(\widehat{c}) = \mathcal{D}_{\Phi}(c(\hspace{0.1cm} \cdot \hspace{0.1cm} ; \hspace{0.1cm} \widehat{\boldsymbol{\theta}}_{n})) = \mathcal{D}_{\Phi}(\widehat{\boldsymbol{\theta}}_{n})$. Still, the estimator $\mathcal{D}_{\Phi}(\widehat{c})$ would often require high-dimensional numerical integration, because usually the integral in \eqref{eq: parmap} does not have a closed-form expression in terms of $\boldsymbol{\theta}_{C}$. That is why we suggest the following general approach. 

Let $M_{n}$ be a user defined sample size for each $n$, and $\widetilde{\mathbf{U}}^{(1)},\dots,\widetilde{\mathbf{U}}^{(M_{n})}$ with $\widetilde{\mathbf{U}}^{(\ell)} = (\widetilde{\mathbf{U}}_{1}^{(\ell)},\dots,\widetilde{\mathbf{U}}_{k}^{(\ell)})$ for $\ell = 1,\dots, M_{n}$, a sample from $\widetilde{\mathbf{U}}$ having distribution $C(\hspace{0.1cm} \cdot \hspace{0.1cm};\widetilde{\boldsymbol{\theta}})$ given that $\widehat{\boldsymbol{\theta}}_{n} = \widetilde{\boldsymbol{\theta}}$, that is
\begin{equation*}
    \mathbb{P} \left (\widetilde{\mathbf{U}} \leq \widetilde{\mathbf{u}} \hspace{0.05cm} | \hspace{0.05cm} \widehat{\boldsymbol{\theta}}_{n} = \widetilde{\boldsymbol{\theta}} \right ) = C(\widetilde{\mathbf{u}};\widetilde{\boldsymbol{\theta}}).
\end{equation*}
We then propose 
\begin{equation}\label{eq: parest}
    \widehat{\mathcal{D}}_{\Phi,n,M_{n}} = \frac{1}{M_{n}} \sum_{\ell = 1}^{M_{n}} \left \{\frac{\prod_{i=1}^{k} c_{i} \left (\widetilde{\mathbf{U}}_{i}^{(\ell)} ; \widehat{\boldsymbol{\theta}}_{n} \right )}{c \left (\widetilde{\mathbf{U}}^{(\ell)} ; \widehat{\boldsymbol{\theta}}_{n} \right )} \Phi \left (\frac{c \left (\widetilde{\mathbf{U}}^{(\ell)} ; \widehat{\boldsymbol{\theta}}_{n} \right )}{\prod_{i=1}^{k} c_{i} \left (\widetilde{\mathbf{U}}_{i}^{(\ell)} ; \widehat{\boldsymbol{\theta}}_{n} \right )} \right ) \right \}
\end{equation}
as estimator for the $\Phi$-dependence $\mathcal{D}_{\Phi}$ given in \eqref{eq: phidiv absc}. The rationale is that $\Phi$-dependence measures can be seen as an expectation and the empirical mean in \eqref{eq: parest} tries to approximate this integral. Since we have an explicit form of the estimated copula, we can take as many samples as we want when $\widehat{\boldsymbol{\theta}}_{n}$ is given, and the larger $M_{n}$, the better the approximation will be. It is intuitively clear (yet not trivial to prove) that, for $M_{n}$ large enough in some sense, the asymptotic properties of $\mathcal{D}_{\Phi}(\widehat{c})$ carry over to $\widehat{\mathcal{D}}_{\Phi,n,M_{n}}$. A conditional law of large numbers will give substance to this, see Remark 4. Theorem 2 states the asymptotic normality result for $\widehat{\mathcal{D}}_{\Phi,n,M_{n}}$.
\newline \\ \noindent
\textbf{Theorem 2.} \textit{Let $\widehat{\boldsymbol{\theta}}_{n}$ be an asymptotically normal estimator for $\boldsymbol{\theta}_{C}$ based on which the estimator \eqref{eq: parest} is constructed, where the user chosen parameter $M_{n}$ is such that
\begin{equation*}
    \left |\widehat{\mathcal{D}}_{\Phi,n,M_{n}} - \mathcal{D}_{\Phi}\left (\widehat{\boldsymbol{\theta}}_{n} \right ) \right | = \smallO_{p}{\left (\frac{1}{\sqrt{n}} \right )} , 
\end{equation*}
for $n \to \infty$.
If the mapping defined in \eqref{eq: parmap} is Fr\'{e}chet differentiable, the estimator $\widehat{\mathcal{D}}_{\Phi,n,M_{n}}$ is asymptotically normal. Moreover, if the derivative can be moved into the integral, we have the expression for the asymptotic variance-covariance:
\begin{equation*}
    \sqrt{n} \left (\widehat{\mathcal{D}}_{\Phi,n,M_{n}} - \mathcal{D}_{\Phi}(c) \right ) \xrightarrow{d} \mathcal{N}(0,\boldsymbol{\beta}^{\text{T}}\mathbf{V}\boldsymbol{\beta})
\end{equation*}
as $n \to \infty$, where $\mathbf{V}$ is the asymptotic variance-covariance matrix of $\widehat{\boldsymbol{\theta}}_{n}$, and
\begin{equation*}
    \boldsymbol{\beta}^{\text{T}} = \begin{pmatrix} \beta_{1} & \beta_{2} & \cdots & \beta_{D} \end{pmatrix} \hspace{0.2cm} \text{with} \hspace{0.2cm} \beta_{i} = \int_{\mathbb{I}^{q}} \frac{\partial f}{\partial \theta_{C,i}}(\mathbf{u};\boldsymbol{\theta}_{C}) d\mathbf{u},
\end{equation*}
defining $D = \text{dim}(\boldsymbol{\Theta}_{C})$, the $i$'th component of $\boldsymbol{\theta}_{C}$ as $\theta_{C,i}$, and
\begin{equation*}
    f(\mathbf{u};\boldsymbol{\theta}_{C}) = \prod_{i=1}^{k} c_{i} \left (\mathbf{u}_{i};\boldsymbol{\theta}_{C} \right ) \Phi \left (\frac{c(\mathbf{u};\boldsymbol{\theta}_{C})}{\prod_{i=1}^{k}c_{i}(\mathbf{u}_{i};\boldsymbol{\theta}_{C})} \right ).
\end{equation*}}
\\
\noindent
\textit{Remark 4.} A conditional version of Kolmogorov’s law of large numbers (see Theorem 4.2 in \cite{Majerek2005}) yields, for all $n \in \mathbb{N}$,
\begin{equation*}
    \sqrt{n} \widehat{\mathcal{D}}_{\Phi,n,M_{n}} \xrightarrow{p} \sqrt{n} \mathcal{D}_{\Phi} \left (\widehat{\boldsymbol{\theta}}_{n} \right ),
\end{equation*}
as $M_{n} \to \infty$. Formally, this means that
\begin{equation*}
    \text{\scalebox{0.95}{$\forall n \in \mathbb{N} : \forall \epsilon, \delta > 0 : \exists M_{n}^{(0)} \in \mathbb{N} : \forall M_{n}^{(1)} \geq M_{n}^{(0)} : \mathbb{P} \left (\sqrt{n} \left |\widehat{\mathcal{D}}_{\Phi,n,M_{n}^{(1)}} - \mathcal{D}_{\Phi} \left (\widehat{\boldsymbol{\theta}}_{n} \right ) \right | > \epsilon \right ) < \delta.$}}
\end{equation*}
This implies that the condition on $M_{n}$ in Theorem 2, stating that
\begin{equation*}
    \forall \epsilon, \delta > 0 : \exists N \in \mathbb{N} : \forall n \geq N : \mathbb{P} \left (\sqrt{n} \left |\widehat{\mathcal{D}}_{\Phi,n,M_{n}} - \mathcal{D}_{\Phi} \left (\widehat{\boldsymbol{\theta}}_{n} \right ) \right | > \epsilon \right ) < \delta,
\end{equation*}
is reasonable. Indeed, take $\epsilon,\delta > 0$ arbitrary. Take for example $N = 1$. Let $n \geq 1$ arbitrary. For these $\epsilon,\delta$ and $n$, we can take $M_{n} = M_{n}^{(0)}$ since $M_{n}$ is user defined, and hence 
\begin{equation*}
    \mathbb{P} \left (\sqrt{n} \left |\widehat{\mathcal{D}}_{\Phi,n,M_{n}} - \mathcal{D}_{\Phi} \left (\widehat{\boldsymbol{\theta}}_{n} \right ) \right | > \epsilon \right ) < \delta.
\end{equation*}
A small value of $N$ corresponds to a good approximation of the integral in $\mathcal{D}_{\Phi}$ already for smaller $n$.
\\

Before looking into another class of parametric copulas, we focus on the Gaussian setting once more.
\newline \\ \noindent
\textbf{Example 5.} 
Consider the simple case of estimating, for two univariate random variables $X_{1}, X_{2}$ with cdf's $F_{1},F_{2}$ and Gaussian copula, the parameter
\begin{equation*}
    \rho = \text{Corr} \left (\left (\phi^{-1} \circ F_{1} \right )(X_{1}), \left (\phi^{-1} \circ F_{2} \right )(X_{2}) \right ),
\end{equation*}
with $\phi^{-1}$ the standard normal quantile function and $\text{Corr}$ the Pearson correlation, and afterwards their mutual information $\mathcal{D}_{t\log(t)}(\rho) = -(1/2) \log(1-\rho^{2})$ via $\mathcal{D}_{t \log(t)}(\widehat{\rho})$ for a certain estimator $\widehat{\rho}$ of $\rho$. The following semi-parametric approaches might be considered.
\begin{itemize}
    \item \textit{Approach 1}: Assume a Gaussian copula model and estimate the marginals non-parametrically. In particular, take the estimator \eqref{eq: est cor matrix}.
    \item \textit{Approach 2}:  Assume a Gaussian copula model, estimate the marginals non-parametrically, and the copula parameter via the pseudo-likelihood \eqref{eq: npll}. 
\end{itemize}
If we look at Approach 1 from a matrix point of view, i.e. estimator $\widehat{\mathbf{R}}_{n}$ in \eqref{eq: est cor matrix}, we can apply Theorem 1 with
\begin{equation*}
\begin{split}
    \hspace{3cm} \mathbf{R} = \begin{pmatrix} 1 & \rho \\ \rho & 1 \end{pmatrix}, \hspace{0.2cm}  &   \mathbf{M}_{t \log(t)} = \begin{pmatrix} \frac{\rho^{2}}{2\rho^{2}-2} & \frac{-\rho}{2\rho^{2}-2} \\ \frac{-\rho}{2\rho^{2}-2} & \frac{\rho^{2}}{2\rho^{2}-2} \end{pmatrix} \\ & \hspace{-4cm} 2 \text{Tr} \left (\left (\mathbf{R}\left (\mathbf{M}_{t \log(t)} - \mathbf{D}_{\mathbf{M}_{t \log(t)} \mathbf{R}} \right ) \right )^{2} \right ) = 2 \text{Tr} \begin{pmatrix} \frac{\rho^{2}}{4} & 0 \\ 0 & \frac{\rho^{2}}{4} \end{pmatrix} = \rho^{2},
\end{split}
\end{equation*}
such that 
\begin{equation*}
    \sqrt{n} \left (\mathcal{D}_{t \log(t)} \left (\widehat{\mathbf{R}}_{n} \right ) - \mathcal{D}_{t \log(t)}(\mathbf{R}) \right ) \xrightarrow{d} \mathcal{N} \left (0,\rho^{2} \right ), 
\end{equation*}
as $n \to \infty$. From a one parameter $\rho$ point of view, with estimator $\widehat{\rho}_{1}$ that is on the off-diagonal of $\widehat{\mathbf{R}}_{n}$, we know that 
\begin{equation*}
\sqrt{n} \left (\widehat{\rho}_{1} - \rho \right ) \xrightarrow{d} \mathcal{N} \left (0, (1-\rho^{2})^{2} \right ),
\end{equation*}
as $n \to \infty$, since, as seen in the proof of Theorem 1, the estimator $\widehat{\rho}_{1}$ has the same asymptotic distribution as in the case where the marginals are known, i.e. as the usual sample Pearson correlation in case of a bivariate normal distribution. Hence, noting that $d/d\rho \mathcal{D}_{t\log(t)}(\rho) = \rho / (1-\rho^{2})$, the univariate delta method implies that 
\begin{equation*}
    \sqrt{n} \left (\mathcal{D}_{t \log(t)} \left (\widehat{\rho}_{1} \right ) - \mathcal{D}_{t \log(t)}(\rho) \right ) \xrightarrow{d} \mathcal{N} \left (0, (1-\rho^{2})^{2} \frac{\rho^{2}}{(1-\rho^{2})^{2}} \right ) = \mathcal{N} \left (0,\rho^{2} \right ),
\end{equation*}
as it should.

As for Approach 2, the resulting estimator, say $\widehat{\rho}_{2}$, has an asymptotic variance that is rather complex to calculate. In \cite{Genest1995}, it is shown that this asymptotic variance cannot be smaller than the one of the maximum likelihood estimator in case the marginals are known. Thus, since $\widehat{\rho}_{1}$ has the same asymptotic distribution as the maximum likelihood estimator in a bivariate Gaussian model, $\widehat{\rho}_{2}$ cannot have a smaller asymptotic variance than $\widehat{\rho}_{1}$, and hence the asymptotic variance of $\mathcal{D}_{t \log(t)}(\widehat{\rho}_{2})$ cannot be smaller than the asymptotic variance of $\mathcal{D}_{t \log(t)}(\widehat{\rho}_{1})$. We conclude that Approach 1, in which we have a nice explicit formula for the estimator, always performs at least equally well as Approach 2 in terms of asymptotic variance. 
\newline \\
\textit{5.2. Nested Archimedean copulas framework}

\begin{figure}[h!] \centering
\includegraphics[width = 0.49\textwidth,height = 6cm]{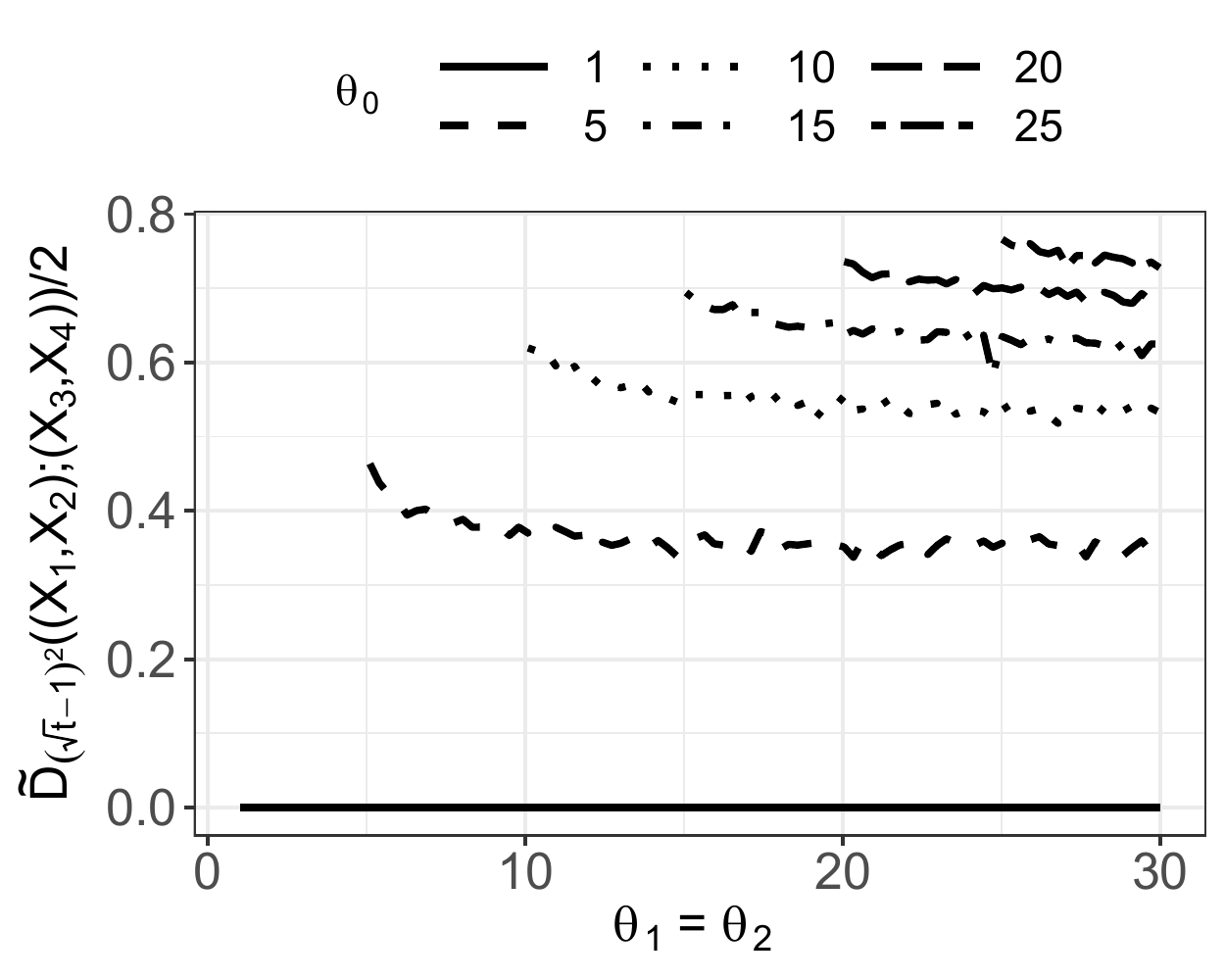}
\includegraphics[width = 0.49\textwidth,height = 6cm]{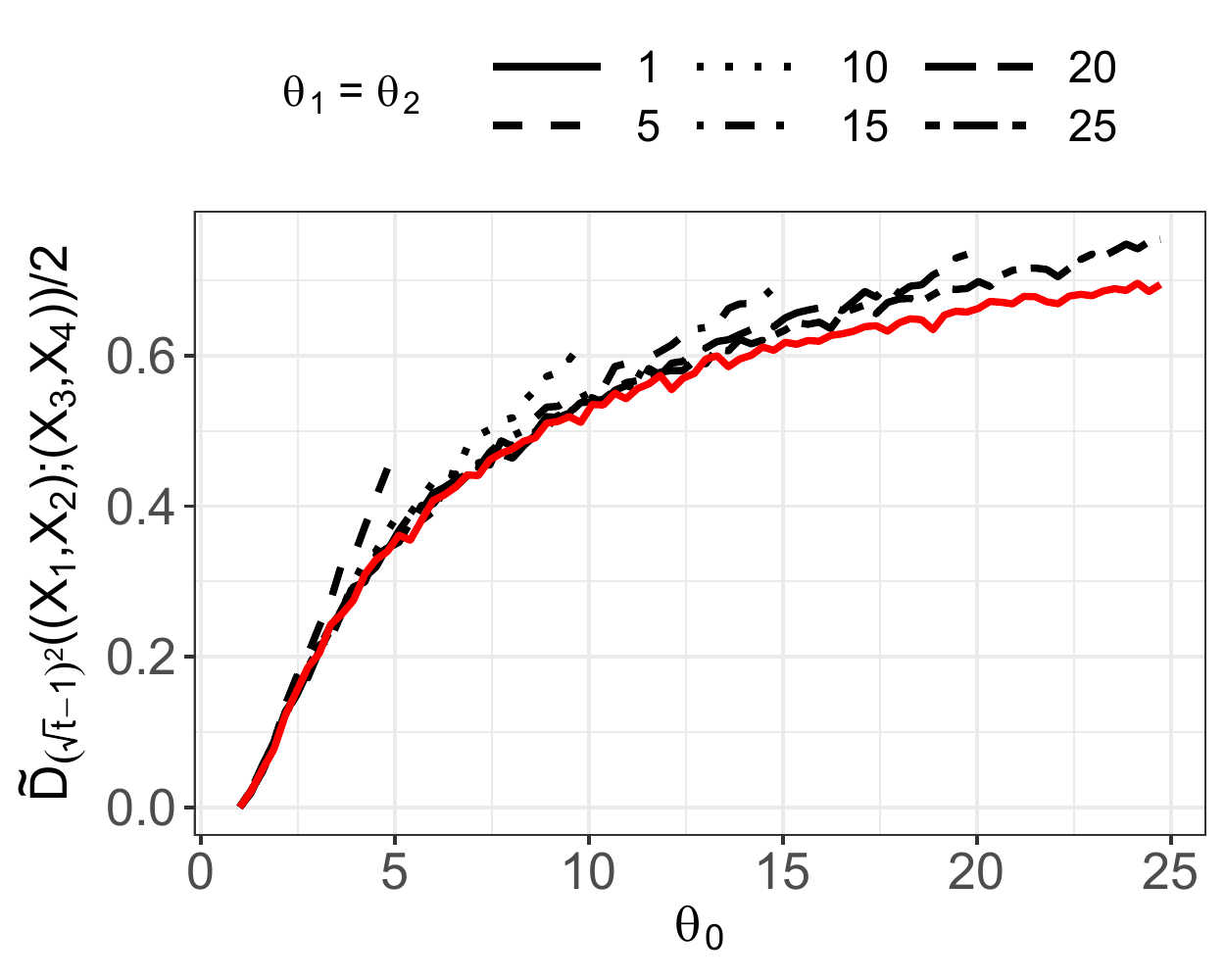}
\caption{Half Hellinger distance of copula \eqref{eq: xmplNAC} as a function of $\theta_{1} = \theta_{2}$ for different values of $\theta_{0}$ (left), and as a function of $\theta_{0}$ for different values of $\theta_{1} = \theta_{2}$ (right). The red line in the right plot shows half the Hellinger distance of a bivariate Gumbel copula with parameter $\theta_{0}$.}
\label{fig: HelNAC}
\end{figure}

We now turn some attention to a specific parametric family of copulas. The hierarchical models of nested Archimedean copulas extend the frequently used Archimedean copulas and are an instinctive choice for modelling dependence between random vectors. See e.g. \cite{Hofert2013} for the following definition.
\newline \\ \noindent
\textbf{Definition 2. (nested Archimedean copulas)} A nested Archimedean copula $C$ with two nesting levels and $k$ child copulas is given by 
\begin{equation}\label{eq: NAC}
    C(\mathbf{u}) = C_{0} \big (C_{1}(\mathbf{u}_{1}),\dots,C_{k}(\mathbf{u}_{k}) \big ),
\end{equation}
where $k$ denotes the dimension of the \textit{root copula} $C_{0}$, and each \textit{child copula} $C_{i}$ for $i \in \{0,\dots, k\}$ is an Archimedean copula with a completely monotone generator $\psi_{i}$, that is
\begin{equation*}
   C_{i}(\mathbf{u}_{i}) = \psi_{i} \big (\psi_{i}^{-1}(u_{i1}),\dots, \psi_{i}^{-1}(u_{id_{i}}) \big ) 
\end{equation*}
and $\psi_{i} : [0,\infty) \to [0,1]$ is continuous with $\psi_{i}(0) = 1, \lim_{t \to \infty} \psi_{i}(t) = 0$ and $(-1)^{\ell}\psi_{i}^{(\ell)}(t) \geq 0$ for all $\ell \in \mathbb{N}, t \in (0,\infty)$.
\newline \\
Note that we can further nest the child copulas in \eqref{eq: NAC}, although this is superfluous for our purposes (in particular because that makes densities excessively complicated). The condition $(-1)^{\ell}\psi_{i}^{(\ell)}(t) \geq 0$ for all $\ell \in \mathbb{N}, t \in (0,\infty)$ is called complete monotonicity of the function $\psi_{i}$ and a sufficient condition to guarantee that \eqref{eq: NAC} indeed is a copula, is that $\psi_{0}^{-1} \circ \psi_{i}$ for all $i \in \{1,\dots,k\}$ have completely monotone first order derivatives. The latter condition is often softened (but definitely not equivalent) to the sufficient nesting condition (e.g. \cite{Okhrin2014}), telling us that $\psi_{i}$ all being in a same family of Archimedean copulas for $i \in \{0,\dots,k\}$, say with parameter $\theta_{i}$, such that $\theta_{0} \leq \theta_{i}$ for $i \in \{1,\dots,k\}$ suffices to have complete monotonicity of these derivatives. \cite{Hofert2013} explicitly calculated the copula density of \eqref{eq: NAC} in several settings. Before looking at an example of the behaviour of a $\Phi$-dependence measure in case of a nested Archimedean copula, we make the following remark.
\newline \\ \noindent 
\textit{Remark 5.} The estimator \eqref{eq: parest} is quite a general one in terms of a certain copula family and function $\Phi$, motivated from the population version \eqref{eq: phidiv absc}. In some cases, it might be better to first simplify \eqref{eq: phidiv absc} and afterwards do the empirical mean approximation. For example, suppose the interest is in the Hellinger distance $\Phi(t) = (\sqrt{t}-1)^{2}$. A straightforward calculation shows that \eqref{eq: phidiv absc} can also be written as 
\begin{equation*}
    \mathcal{D}_{(\sqrt{t}-1)^{2}} = 2 - 2 \int_{\mathbb{I}^{q}} \sqrt{c(\mathbf{u}_{1},\dots,\mathbf{u}_{k})c_{1}(\mathbf{u}_{1}) \cdots c_{k}(\mathbf{u}_{k})} d\mathbf{u},
\end{equation*}
suggesting 
\begin{equation}\label{eq: helapprox}
    \widetilde{\mathcal{D}}_{(\sqrt{t}-1)^{2}} = 2 - 2 \frac{1}{M} \sum_{\ell = 1}^{M} \sqrt{\frac{c_{1} \left (\mathbf{U}_{1}^{(\ell)} \right ) \cdots c_{k} \left (\mathbf{U}_{k}^{(\ell)} \right )}{c \left (\mathbf{U}_{1}^{(\ell)},\dots,\mathbf{U}_{k}^{(\ell)} \right )}},
\end{equation}
where $(\mathbf{U}_{1}^{(\ell)},\dots,\mathbf{U}_{k}^{(\ell)})$ for $\ell = 1,\dots,M$ is a sample drawn from $c(\hspace{0.01cm} \cdot \hspace{0.01cm} ; \boldsymbol{\theta}_{C})$, as numerical approximation, assuming here that no estimation is done (if $\boldsymbol{\theta}_{C}$ is estimated, we use the notation $\widetilde{\mathcal{D}}_{(\sqrt{t}-1)^{2},n,M_{n}}$). The benefit of $\widetilde{\mathcal{D}}_{(\sqrt{t}-1)^{2}}$ is that 
\begin{equation*}
    \mathbb{E} \left [ \left ( \sqrt{\frac{c_{1} \left (\mathbf{U}_{1}^{(\ell)} \right ) \cdots c_{k} \left (\mathbf{U}_{k}^{(\ell)} \right )}{c \left (\mathbf{U}_{1}^{(\ell)},\dots,\mathbf{U}_{k}^{(\ell)} \right )}} \right )^{2} \right ] = \int_{\mathbb{I}^{q}} c_{1}(\mathbf{u}_{1}) \cdots c_{k}(\mathbf{u}_{k}) d\mathbf{u}_{1} \cdots d\mathbf{u}_{k}  = 1,
\end{equation*}
guaranteeing fast convergence of the law of large numbers, while the variance of the summand in \eqref{eq: parest} might be infinite. We further illustrate this in the following example.
\newline \\ \noindent
\textbf{Example 6.}
Consider $(X_{1},X_{2},X_{3},X_{4})$ having a four dimensional partially nested Archimedean copula given by 
\begin{equation}\label{eq: xmplNAC}
        \text{\scalebox{0.89}{$C(u_{1},u_{2},u_{3},u_{4}) = \psi_{0} \left ( \psi_{0}^{-1} \left ( \psi_{1} \left (\psi_{1}^{-1}(u_{1}) + \psi_{1}^{-1}(u_{2}) \right) \right) + \psi_{0}^{-1} \left ( \psi_{2} \left ( \psi_{2}^{-1}(u_{3}) + \psi_{2}^{-1}(u_{4}) \right ) \right ) \right ),$}}
\end{equation}
where $\psi_{i}(t) = \exp(-t^{1/\theta_{i}})$ is the generator of a Gumbel copula with parameter $\theta_{i} \in [1,\infty)$ for $i = 0,1,2$ satisfying $\theta_{0} \leq \theta_{1}$ and $\theta_{0} \leq \theta_{2}$ (sufficient nesting condition). We numerically approximate half the Hellinger distance $\mathcal{D}_{(\sqrt{t}-1)^{2}}/2$ using \eqref{eq: helapprox} with $M = 10\hspace{0.05cm}000$. Note that no estimation of the marginals is involved here. Figure \ref{fig: HelNAC} shows half the Hellinger distance as a function of $\theta_{1} = \theta_{2}$ for different values of $\theta_{0}$, and as a function of $\theta_{0}$ for different values of $\theta_{1} = \theta_{2}$.
Note that if $\theta_{0} = 1$, we have $\psi_{0}(t) = \exp(-t)$ such that 
\begin{equation*}
    \begin{split}
        C(u_{1},u_{2},u_{3},u_{4}) = C_{1}(u_{1},u_{2}) C_{2}(u_{3},u_{4}),
    \end{split}
\end{equation*}
meaning that $(X_{1},X_{2})$ and $(X_{3},X_{4})$ are independent.
In general, we observe that the strength of dependence between $(X_{1},X_{2})$ and $(X_{3},X_{4})$ is predominantly determined by the parameter $\theta_{0}$. Nested Archimedean copulas allow us to on the one hand control the dependence within each random vector (by parameters $\theta_{1}$ and $\theta_{2}$ here) and on the other hand control what remains to be specified between the random vectors (parameter $\theta_{0}$ here). Notice that if $\theta_{1} \to \infty$ and $\theta_{2} \to \infty$ (i.e. $(X_{1},X_{2})$ and $(X_{3},X_{4})$ both tend to have a comonotonicity copula), the half Hellinger distance $\mathcal{D}_{(\sqrt{t}-1)^{2}}/2$ tends to the half Hellinger distance of a bivariate Gumbel copula with parameter $\theta_{0}$ (red line in the right plot).

Coming back to Remark 5, suppose we would compute the Hellinger distance of this bivariate Gumbel copula with parameter $\theta_{0}$, say $c(u_{1},u_{2};\theta_{0})$, using the non-simplified empirical version \eqref{eq: parest} (but now with known $\theta_{0}$ and a true sample from $c$ with fixed $M$) instead of $\eqref{eq: helapprox}$, i.e. using
\begin{equation}\label{eq: hatD}
    \widehat{\mathcal{D}}_{(\sqrt{t}-1)^{2}} = \frac{1}{M} \sum_{\ell = 1}^{M} \frac{1}{c \left (U_{1}^{(\ell)},U_{2}^{(\ell)} \right )} \left ( \sqrt{c \left (U_{1}^{(\ell)},U_{2}^{(\ell)} \right )} - 1 \right )^{2},
\end{equation}
where $(U_{1}^{(\ell)},U_{2}^{(\ell)})$ for $\ell = 1,\dots,M$ is a sample drawn from $c$. Then, the pitfall is that
\begin{equation*}
   \text{\scalebox{0.95}{$\mathbb{E} \left [ \left (\frac{1}{c \left (U_{1}^{(\ell)},U_{2}^{(\ell)} \right )} \left ( \sqrt{c \left (U_{1}^{(\ell)},U_{2}^{(\ell)} \right )} - 1 \right )^{2} \right )^{2} \right ] = \int_{\mathbb{I}^{2}} \frac{1}{c(u_{1},u_{2})} \left (\sqrt{c(u_{1},u_{2})} - 1\right )^{4} du_{1}du_{2}$}}
\end{equation*}
\begin{equation*}
    \hspace{2cm} = \infty 
\end{equation*}
for $\theta_{0} \geq 2$ since 
\begin{equation*}
    \frac{1}{c(u_{1},u_{2})} \left (\sqrt{c(u_{1},u_{2})} - 1\right )^{4} = \mathcal{O} \left (\frac{1}{(-\log(u_{2}))^{\theta_{0}-1}} \right ) \hspace{0.2cm} \text{as} \hspace{0.2cm} u_{2} \to 1
\end{equation*}
for a fixed $u_{1} \in (0,1)$. The above can be seen from the fact that for the Gumbel generator $\psi(t) = \exp(-t^{1/\theta_{0}})$, it holds that $\psi^{\prime}(\psi^{-1}(u_{2})) = \mathcal{O}((-\log(u_{2}))^{1-\theta_{0}})$ as $u_{2} \to 1$. Hence, the law of large numbers still holds, but the convergence will be slower due to infinite variance, see also Section 6.2.
\newline \\ \noindent
\textbf{6. Simulation experiments} 
\\

In this section, we perform some simulations that ought to complement the theoretical results that we obtained. First, we focus on Theorem 1 and investigate how well the asymptotic normal distribution approximates the finite-sample distribution of the plug-in estimator for the mutual information and Hellinger distance in case of a Gaussian copula model. Second, we numerically assess how well the estimator \eqref{eq: parest} performs in the setting of Example 6, and compare the numerical quality of $\widehat{\mathcal{D}}$ given in \eqref{eq: hatD} with the numerical quality of $\widetilde{\mathcal{D}}$ given in \eqref{eq: helapprox} within a bivariate Gumbel copula model. 
\newline \\ \noindent
\textit{6.1. Asymptotic normality under Gaussian copula model} 

\begin{figure}[h!] \centering
\includegraphics[scale = 0.75]{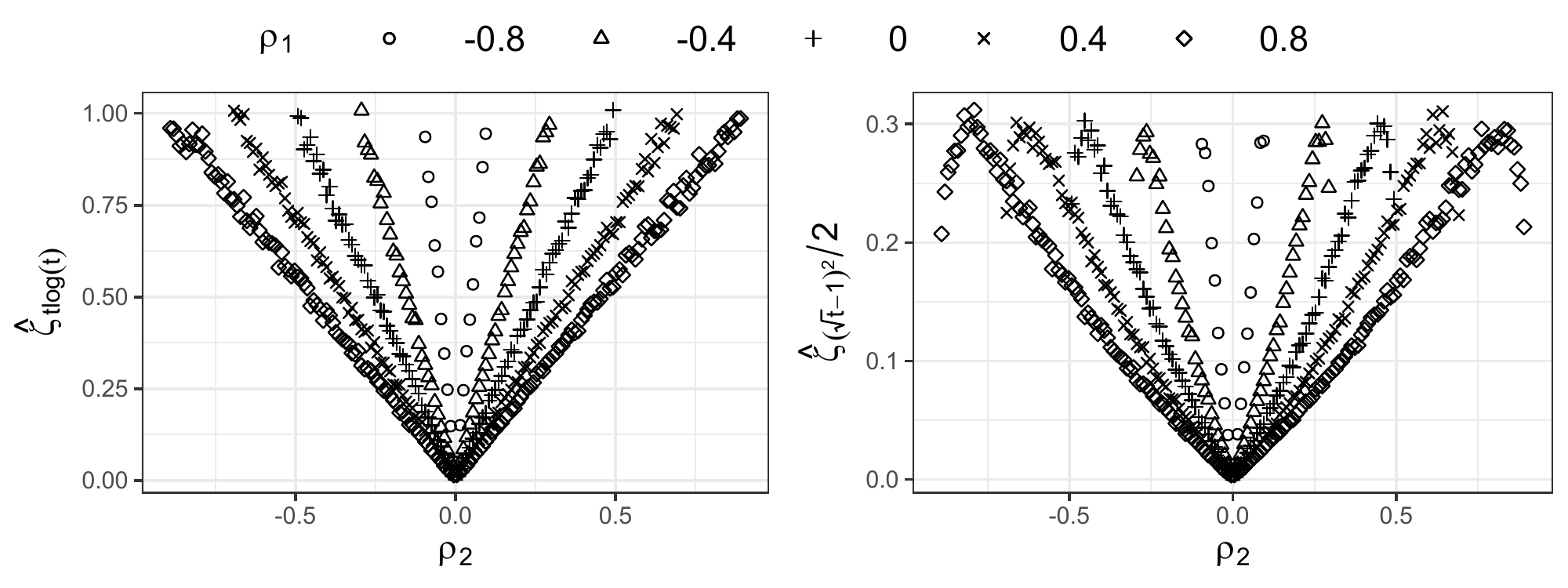}
\hspace*{0.4cm}
\includegraphics[scale = 0.55]{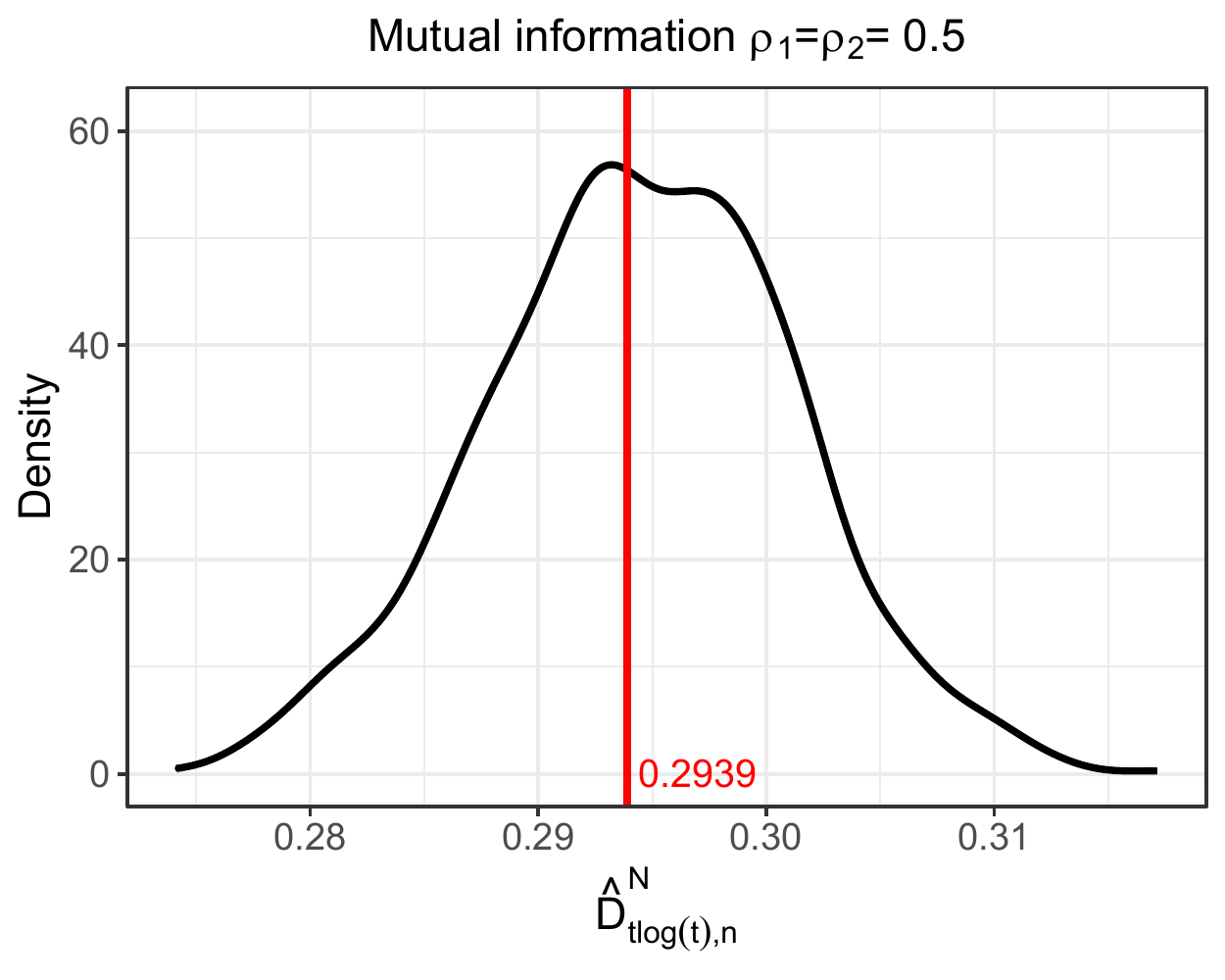}
\hspace{0.1cm}
\includegraphics[scale = 0.55]{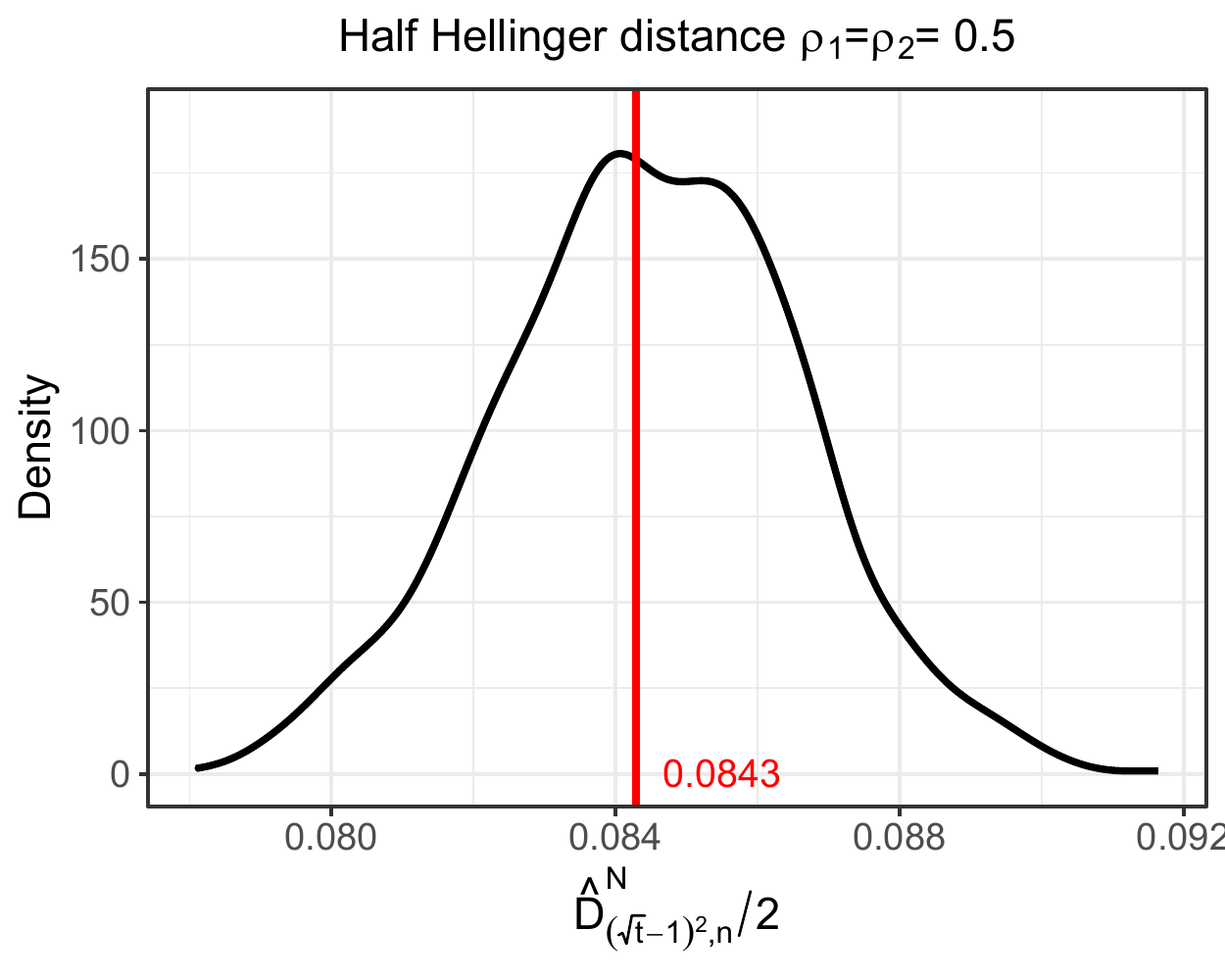}
\caption{Empirical standard deviation (sample size $N = 1000$) $\widehat{\zeta}_{t \log(t)}$ and $\widehat{\zeta}_{(\sqrt{t}-1)^{2}}/2$ in the setting of Example 4, and kernel density estimates for $\widehat{\mathcal{D}}^{\mathcal{N}}_{t \log(t),n}$ and $\widehat{\mathcal{D}}^{\mathcal{N}}_{(\sqrt{t}-1)^{2},n}/2$ when $\rho_{1} = \rho_{2} = 0.5$ and $n = 10 \hspace{0.05cm} 000$. The red vertical lines indicate the true value of the dependence coefficient.}
\label{fig: zetaxmpl}
\end{figure}

Recall the asymptotic normality result of Theorem 1 for the estimator $\widehat{\mathcal{D}}^{\mathcal{N}}_{\Phi,n} = \mathcal{D}^{\mathcal{N}}_{\Phi}(\widehat{\mathbf{R}}_{n})$, with $\widehat{\mathbf{R}}_{n}$ the matrix of sample normal scores rank correlation coefficients \eqref{eq: est cor matrix} and $\mathcal{D}^{\mathcal{N}}_{\Phi}$ given in \eqref{eq: mutN} and \eqref{eq: helN} for $\Phi(t) = t \log(t)$ and $\Phi(t) = (\sqrt{t}-1)^{2}$ respectively. We first turn our attention to Example 4 once more. In Figure \ref{fig: xmplNASD}, we depicted the asymptotic standard deviation 
$\zeta_{t \log(t)}$ and $\zeta_{(\sqrt{t}-1)^{2}}/2$. If we generate $N$ samples from e.g. a four dimensional multivariate Gaussian distribution with mean zero and covariance matrix $\mathbf{R}$ as in Example 4, we obtain $N$ estimates of $\mathcal{D}_{\Phi}^{\mathcal{N}}$, say $\widehat{\mathcal{D}}^{\mathcal{N},(1)}_{\Phi,n},\dots,\widehat{\mathcal{D}}^{\mathcal{N},(N)}_{\Phi,n}$. An empirical Monte Carlo version of $\sqrt{n}\text{SD}(\widehat{\mathcal{D}}_{\Phi,n}^{\mathcal{N}})$, with SD being the standard deviation, is then given by
\begin{equation*}
    \widehat{\zeta}_{\Phi} = \sqrt{\frac{n}{N-1} \sum_{\ell = 1}^{N} \left (\widehat{\mathcal{D}}^{\mathcal{N},(\ell)}_{\Phi,n} - \frac{1}{N}\sum_{t = 1}^{N} \widehat{\mathcal{D}}^{\mathcal{N},(t)}_{\Phi,n}  \right )^{2}},
\end{equation*}
which we can compute for different values of $\rho_{1}$ and $\rho_{2}$, see Figure \ref{fig: zetaxmpl} for some plots in case $n = 10\hspace{0.05cm}000$ and $N = 1000$. Comparing with Figure 3 of Example 4, this empirically verifies the formula for the asymptotic variance in Theorem 1 in this particular setting. Kernel density estimates for the density of $\widehat{\mathcal{D}}^{\mathcal{N}}_{\Phi,n}$ when $\rho_{1} = \rho_{2} = 0.5$ are included as well.
 \newline

Let now $\widehat{\zeta}_{\Phi,n}$ be the plug-in estimator of the asymptotic standard deviation $\zeta_{\Phi}$ obtained by using $\widehat{\mathbf{R}}_{n}$ instead of the true $\mathbf{R}$. Based on a sample $\mathbf{X}^{(\ell)}$ for $\ell = 1,\dots,n$ from a certain multivariate distribution having a Gaussian copula, we are able to compute one realization of the actual sampling distribution of the studentized estimator
$\sqrt{n} \left (\widehat{\mathcal{D}}^{\mathcal{N}}_{\Phi,n} - \mathcal{D}^{\mathcal{N}}_{\Phi}(\mathbf{R}) \right )/\widehat{\zeta}_{\Phi,n}$,
and several replications will give an idea about the entire distribution, which should, according to Theorem 1 approximately be a standard normal one for larger values of $n$. We consider four settings which we can generate samples from:
\begin{itemize}
    \item Setting 1: $k = 2, d_{1} = d_{2} = 2$, with standard normal marginals and a Gaussian copula having an autoregressive AR(1) correlation matrix with $\rho = 0.25$.
    \item Setting 2: as Setting 1, but now with marginals
    \begin{itemize}
        \item[{$\star$}] a $t$ distribution with $3$ degrees of freedom for $X_{11}$
        \item[{$\star$}] an exponential distribution with mean $1$ for $X_{12}$
        \item[{$\star$}] a beta distribution with parameters $2$ and $2$ for $X_{21}$
        \item[{$\star$}] an $F$-distribution with degrees of freedom $2$ and $6$ for $X_{22}$.
    \end{itemize}
    \item Setting 3: similar to Setting 1, but with $\rho = 0.8$.
    \item Setting 4: $k = 5, d_{1} = 4, d_{2} = 5, d_{3} = 3, d_{4} = 1, d_{5} = 2$, with standard normal marginals and a Gaussian copula having an equicorrelated correlation matrix with $\rho = 0.5$.
\end{itemize}
\begin{figure}[h!] \centering
\textbf{Mutual information}
\hspace*{-0.6cm}
\includegraphics[scale = 0.95]{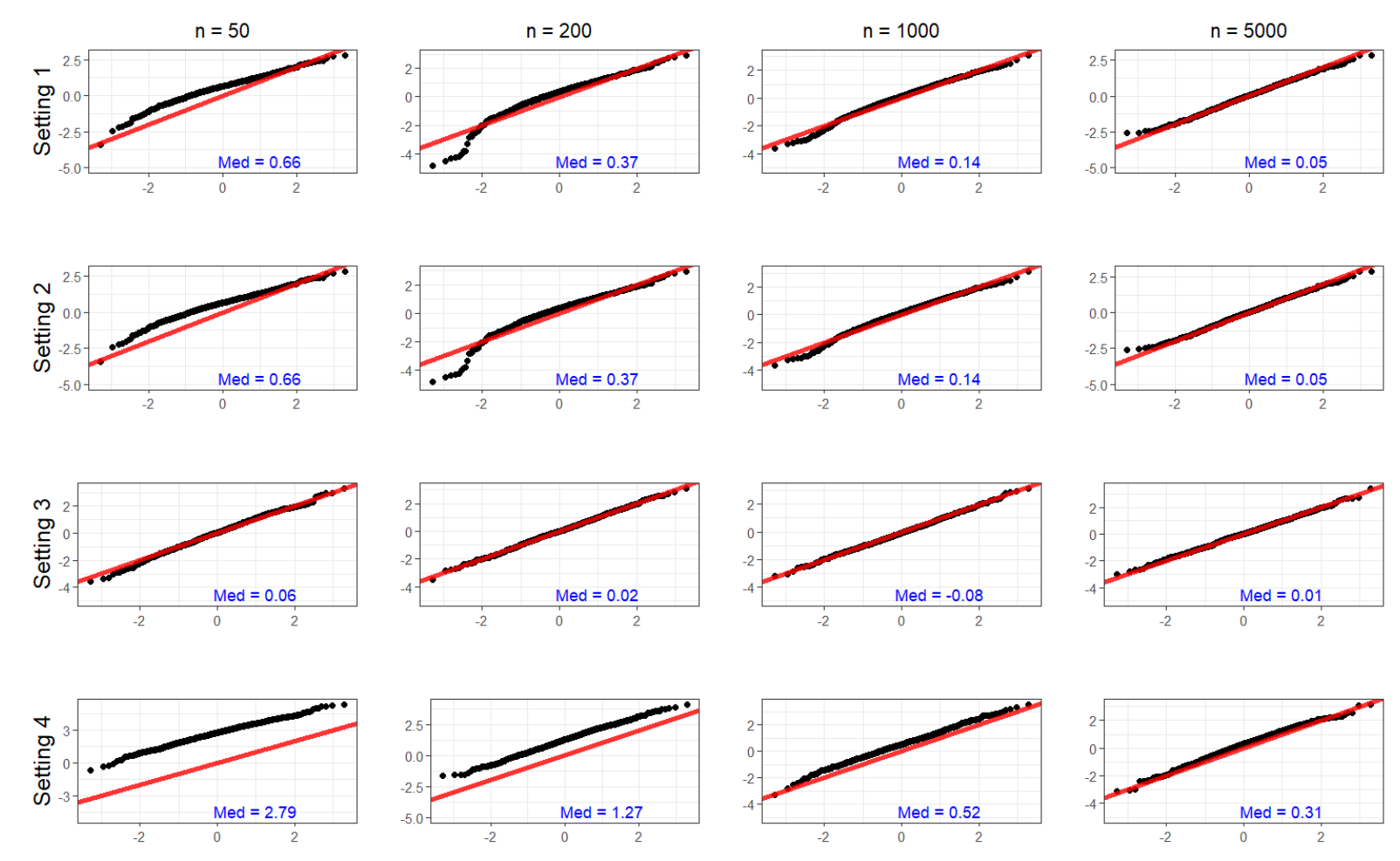}
\caption{Normal Q-Q plots for $1000$ Monte Carlo runs of the studentized plug-in estimator for the mutual information under four different settings with sample sizes $n = 50,200,1000,5000$. The median (``Med") of the studentized estimates is indicated in blue.}
\label{fig: figure2}
\end{figure} 

\begin{figure}[h!] \centering
\textbf{Half Hellinger distance}
\hspace*{-0.5cm}
\includegraphics[scale = 0.95]{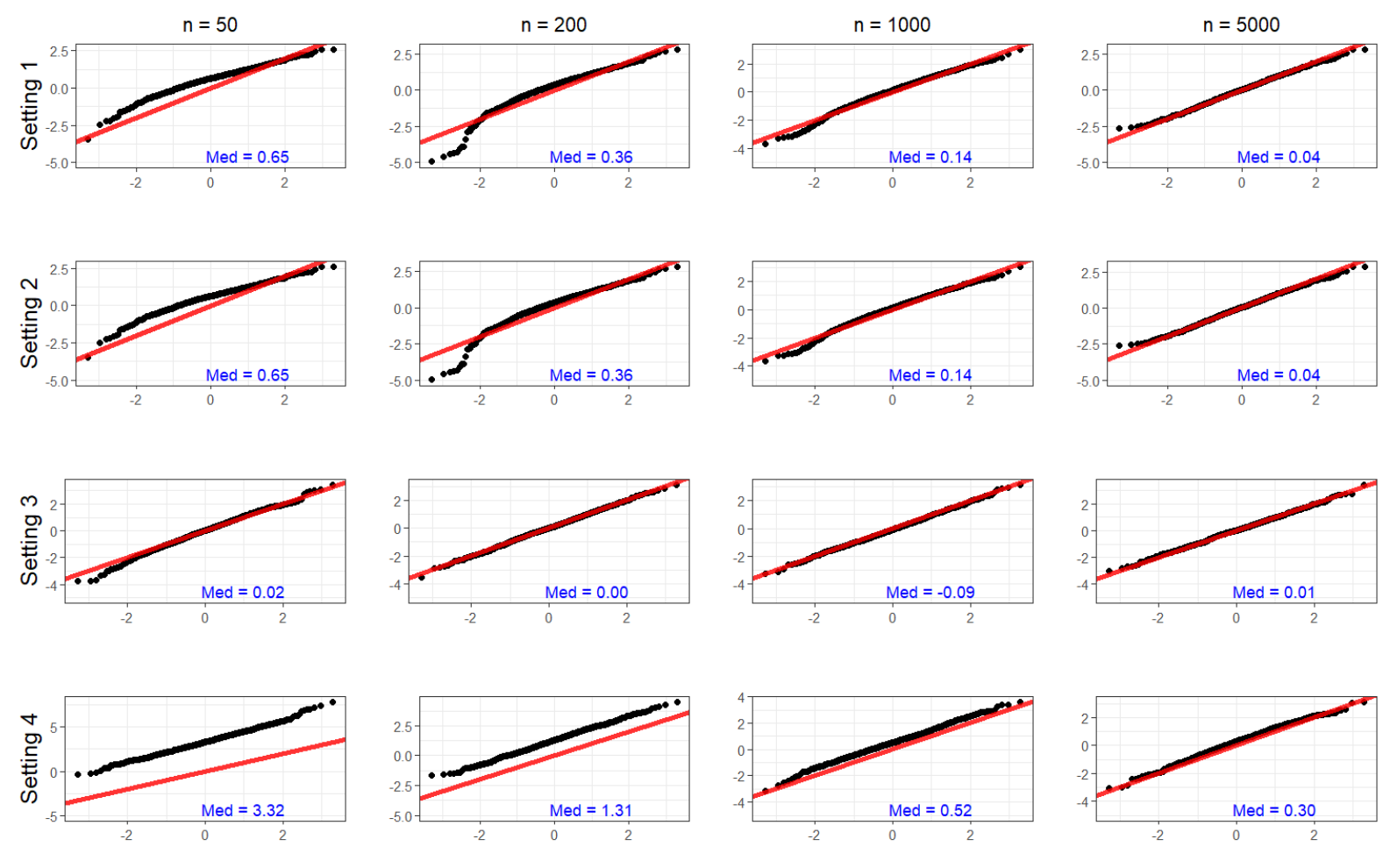}
\caption{Normal Q-Q plots for $1000$ Monte Carlo runs of the studentized plug-in estimator for the half Hellinger distance under four different settings with sample sizes $n = 50,200,1000,5000$. The median (``Med") of the studentized estimates is indicated in blue.}
\label{fig: figure3}
\end{figure}

Each time, we draw $1000$ samples of sizes $n = 50,200,1000,5000$ and make normal Q-Q plots to assess the goodness-of-fit with a standard normal distribution. See Figure \ref{fig: figure2} for the results of the mutual information, and Figure \ref{fig: figure3} for the half Hellinger distance.

In each setting, we have a qualitative normal approximation for larger sample sizes. Sampling from a multivariate normal distribution or from a multivariate normal copula with various marginals does not give a significant difference (Setting 1 versus Setting 2). For rather small correlations (Setting 1 and Setting 2), we observe a more pronounced lack-of-fit than for higher correlations (Setting 3). Increasing the total dimension to $q = 15$ (Setting 4) results in a large positive bias for small sample size, which is no shock since empirical covariance matrices tend to be more biased when the number of parameters to estimate increases. 
\newline \\ \noindent
\textit{6.2. Nested Archimedean copula model} 

In the context of a general parametric copula family, the estimator \eqref{eq: parest} relies on an estimator $\widehat{\boldsymbol{\theta}}_{n}$ for the copula parameter $\boldsymbol{\theta}_{C}$ on the one hand, and on a numerical integral approximation on the other hand. In Example 6, we mathematically illustrated that, in a bivariate Gumbel copula with parameter $\boldsymbol{\theta}_{C} = \theta_{0} = 3$, the estimator $\widehat{\mathcal{D}}_{(\sqrt{t}-1)^{2},n,M_{n}}$ is doomed to have a slower convergence than the estimator $\widetilde{\mathcal{D}}_{(\sqrt{t}-1)^{2},n,M_{n}}$ relying on a simplified integral approximation as in \eqref{eq: helapprox}.

Via a small simulation, we now compare the performance of these two estimators when $\theta_{0} = 3$ and $n = 200$, based on $1000$ samples drawn from this bivariate Gumbel copula. We consider $M_{200} = 100,1000,10\hspace{0.05cm}000$ and compare in Table \ref{tab: table2} the sample bias, variance and mean squared error of $\widehat{\mathcal{D}}_{(\sqrt{t}-1)^{2},n,M_{n}}/2$ and $\widetilde{\mathcal{D}}_{(\sqrt{t}-1)^{2},n,M_{n}}/2$.
The estimator $\widehat{\boldsymbol{\theta}}_{n}$ is based on maximizing the pseudo likelihood with non-parametric marginals, i.e. $\widehat{\boldsymbol{\theta}}_{C,n}^{\text{NP}}$, using a starting value of $1$. The true value of the dependence coefficient equals $\mathcal{D}_{(\sqrt{t}-1)^{2}}/2 = 0.20528$, and was computed, not by doing an empirical mean approximation of the two dimensional integral, but using numerical integration. As expected, the performance of $\widehat{\mathcal{D}}_{(\sqrt{t}-1)^{2},n,M_{n}}/2$ is poor due to a large variance that needs very large values of $M_{n}$ to go down. 

Next, we investigate how well the estimator $\widetilde{\mathcal{D}}_{(\sqrt{t}-1)^{2},n,M_{n}}/2$ performs in terms of increasing $n$ for the half Hellinger distance between $(X_{1},X_{2})$ and $(X_{3},X_{4})$ having the four dimensional copula given in \eqref{eq: xmplNAC} in Example 6. The estimator $\widehat{\mathcal{D}}_{t \log(t),n,M_{n}}$ for the mutual information is also considered. We look at two settings:
\begin{itemize}
    \item Setting 1: $\boldsymbol{\theta}_{C} = (\theta_{0},\theta_{1},\theta_{2}) = (1,3,4)$, such that $\mathcal{D}_{(\sqrt{t}-1)^{2}}/2 = \mathcal{D}_{t \log(t)} = 0$, since $(X_{1},X_{2})$ and $(X_{3},X_{4})$ are independent. 
    \item Setting 2: $\boldsymbol{\theta}_{C} = (\theta_{0},\theta_{1},\theta_{2}) = (3,3,4)$, such that $\mathcal{D}_{(\sqrt{t}-1)^{2}}/2 = 0.29007$ and $\mathcal{D}_{t \log(t)} = 0.99935$, computed via the true $\boldsymbol{\theta}_{C}$ and $10\hspace{0.05cm} 000 \hspace{0.05cm} 000$ samples to numerically approximate the integral. 
\end{itemize}
\begin{table}[h]
\begin{tabularx}{\textwidth}{||YYYYSYYYS||}
\hhline{|=========|}
\footnotesize $\mathbf{M_{200}}$ & \multicolumn{3}{c}{\footnotesize $\mathbf{\widehat{\mathcal{D}}_{(\sqrt{t}-1)^{2},200,M_{200}}/2}$} & & \multicolumn{3}{c}{\footnotesize $\mathbf{\widetilde{\mathcal{D}}_{(\sqrt{t}-1)^{2},200,M_{200}}/2}$} & \\ \cline{2-4} \cline{6-8}
& bias & var & mse & & bias & var & mse & \\ 
$100$ & $-0.0109$ & $0.1313$ & $0.1314$ & & $0.0031$ & $0.0039$ & $0.0039$ & \\
$1000$ & $-0.0078$ & $0.0249$ & $0.025$ & & $0.001$ & $0.0008$ & $0.0008$ & \\
$10 \hspace{0.05cm} 000$ & $0.0007$ & $0.0139$ & $0.0139$ & & $0.0008$ & $0.0005$ & $0.0005$ & \\
\hhline{|=========|}
\end{tabularx}
\caption{Sample bias, variance and mean squared error of two estimators for the half Hellinger distance in a bivariate Gumbel copula with parameter $\theta_{0} = 3$, based on $1000$ replications, a sample size $n = 200$, and $M_{200} = 100,1000,10\hspace{0.05cm}000$. The true value equals $\mathcal{D}_{(\sqrt{t}-1)^{2}}/2 = 0.20528$.}
\label{tab: table2}
\end{table} 

We take $\widehat{\boldsymbol{\theta}}_{n} = (\widehat{\theta}_{0,n},\widehat{\theta}_{1,n},\widehat{\theta}_{2,n}) =  \widehat{\boldsymbol{\theta}}_{C,n}^{\text{NP}}$, with $2$ as starting value for $\theta_{0}$ for maximizing the likelihood. The starting values for $\theta_{1}$ and $\theta_{2}$ are taken as the maximizers of the individual pseudo likelihoods (also with non-parametric marginals) corresponding to the marginal samples of $(X_{1},X_{2})$ and $(X_{3},X_{4})$ respectively, with both starting values equal to $2$. In each setting, we take $1000$ Monte Carlo runs and sample sizes $n = 50,200,1000,5000$ and fix $M_{n} = 10 \hspace{0.05cm} 000$ for every $n$. 

\begin{figure}[h!] \centering
\includegraphics[scale = 0.75]{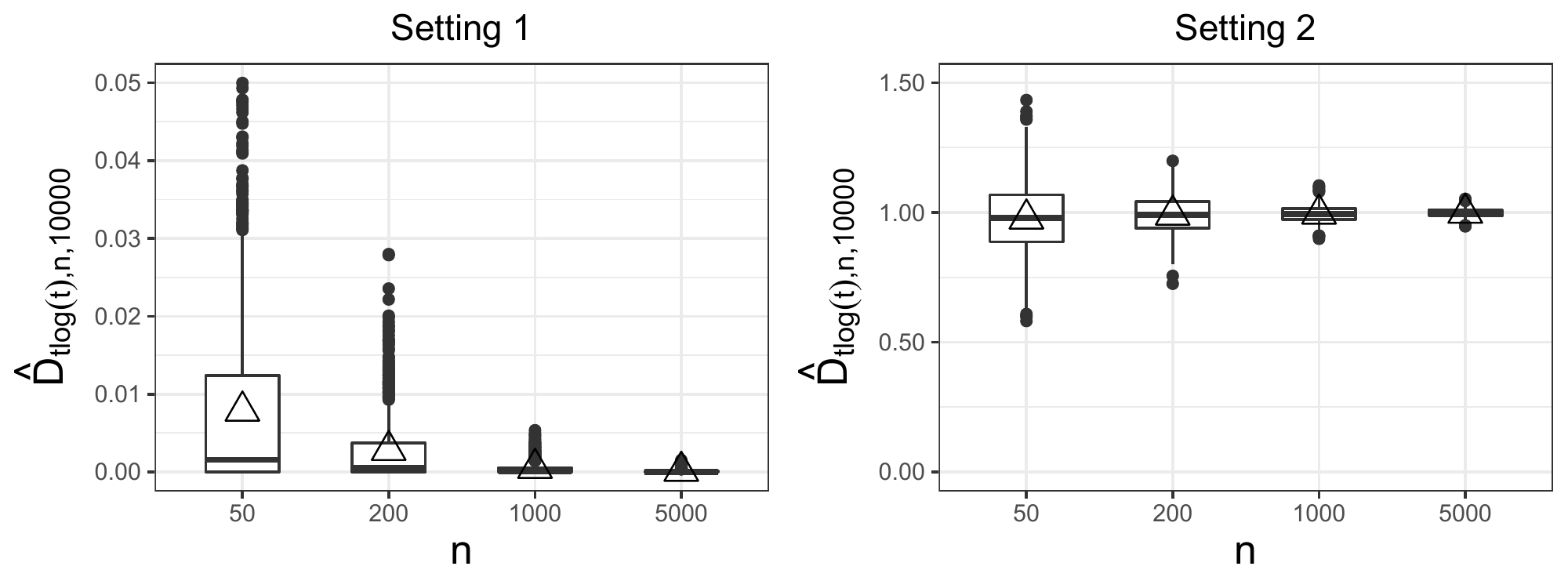}
\includegraphics[scale = 0.75]{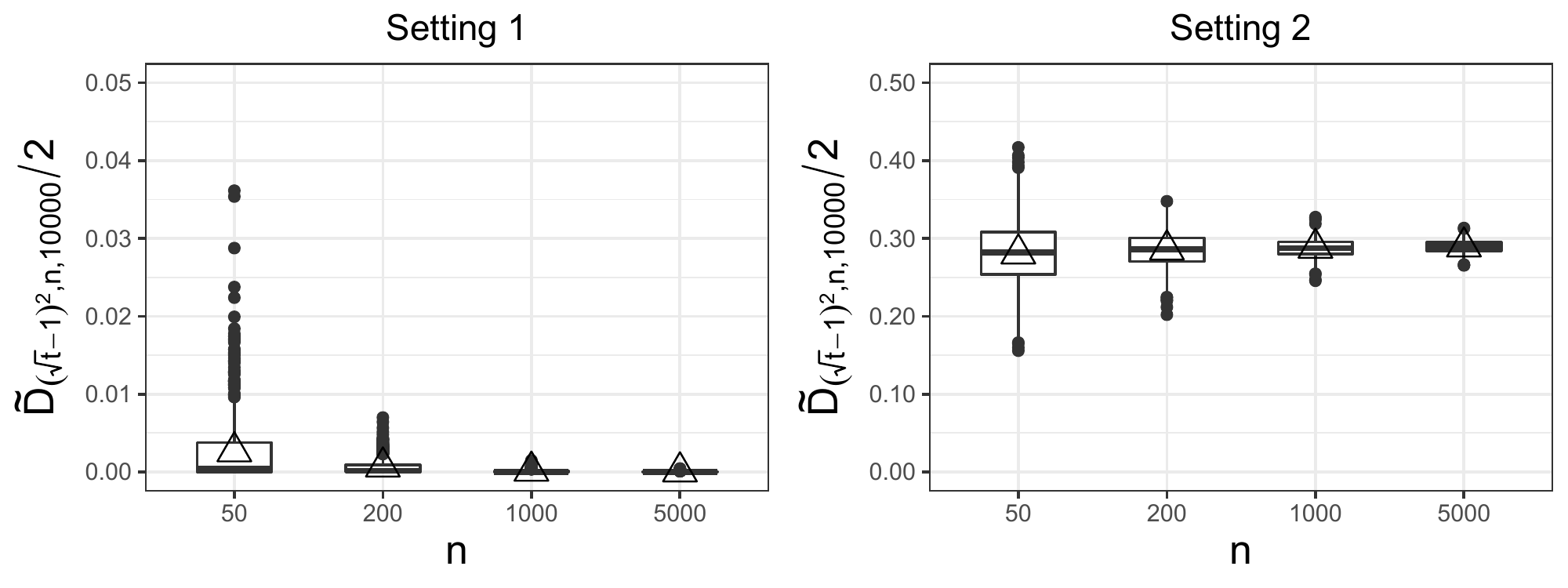}
\caption{Boxplots of estimated dependence coefficients for different sample sizes in different settings. The triangles indicate the mean values.}
\label{fig: NACsim}
\end{figure} 
\begin{table}[h]
\begin{tabularx}{\textwidth}{||YYYSYYS||}
\hhline{|=======|}
\footnotesize $\mathbf{n}$ & \multicolumn{2}{c}{\footnotesize \textbf{Mutual information}} & & \multicolumn{2}{c}{\footnotesize \textbf{Half Hellinger distance}} & \\ \cline{2-3} \cline{5-6}
& Setting 1 & Setting 2 & & Setting 1 & Setting 2 & \\ 
$50$ & $0.002$ & $0.1311$ & & $0.0001$ & $0.0124$ & \\
$200$ & $0.0002$ & $0.0734$ & & $1.19 \cdot 10^{-5}$ &  $0.0073$ & \\
$1000$ & $1.62 \cdot 10^{-5}$ & $0.0322$ & & $1.10 \cdot 10^{-6}$ & $0.0042$ & \\
$5000$ & $2.01 \cdot 10^{-6}$ & $0.0217$ & & $2.35 \cdot 10^{-7}$ & $0.0045$ & \\
\hhline{|=======|}
\end{tabularx}
\caption{Empirical Monte Carlo variances $n \widehat{\text{Var}}$ of estimated dependence coefficients for different sample sizes $n$ in different settings, based on $1000$ replications.}
\label{tab: table3}
\end{table} 

Boxplots are shown in Figure \ref{fig: NACsim}. For each dependence coefficient and in each setting, the bias and variance tend to zero when the sample size increases. Table \ref{tab: table3} shows empirical Monte Carlo variances $n \widehat{\text{Var}}$ of each estimator. They indicate that cases of stronger dependence (Setting 2) are harder to estimate than cases of weak dependence (Setting 1) where the asymptotic variance is smaller (as we have seen for some Gaussian examples too).
\newline \\ \\ \\ \noindent
\textbf{7. A financial application} 
\\

Quantifying the strength of relationship between variables is fundamental in finance, e.g. in portfolio management. Individual constituents are often (positively) related to one another because of contingency on macro-economic factors, known as systematic risk. For instance, market downturns can have detrimental consequences on the portfolio as association between assets can significantly increase. This phenomenon is known as asymmetric dependence, see e.g. \cite{Alcock2018}. Closely related is the concept of financial contagion, e.g. \cite{Gallegati2012}, \cite{Celik2012}, \cite{Wang2017} and \cite{Akhtaruzzaman2021}, among others, evidencing stronger linkages across markets in times of recession. 

These markets can be considered within one and the same region, or can be spread over multiple different regions. We might for instance have a random vector $\mathbf{X}_{1}$ describing equity indexes in North America and look at the intra-dependence, or investigate the inter-dependence with European indexes $\mathbf{X}_{2}$, neglecting the within region dependence. Here, we will use the $\Phi$-dependence measures between $\mathbf{X}_{1}$ and $\mathbf{X}_{2}$, intending to illustrate cross-regional financial contagion during the COVID-19 pandemic. 

In particular, we analyse historical daily logarithmic returns of stock indexes from North-America (US S\&P500, Canadian S\&P/TSX Composite Index and Mexican IPC Index), South-America (Brazilian IBOVESPA and Argentina Merval Index), Europe (Euronext 100, German GDAXI, Spanish IBEX 35 and Norwegian OMX Index) and Asia (Japanese Nikkei 225, Chinese SSE Composite Index, Indian S\&P BSE 500 and Hong Kong HSI Index), over a time span of Dec 07, 2016 to Dec 06, 2022. The data can freely be accessed and downloaded at \url{https://finance.yahoo.com/}. Notice that logarithmic returns are i.i.d. when assuming a random walk market.  Each set of index returns in each continent is considered as a random vector $\mathbf{X}_{i}$, and inter-regional financial dependence is assessed as $\mathcal{D}_{\Phi}(\mathbf{X}_{i},\mathbf{X}_{j})$ for $i \neq j$.

We first make use of the method discussed in Section 4, relying on the assumption of a Gaussian copula model. A primary impression of inter-regional financial contagion during the COVID-19 period is obtained by looking at the dependence $\mathcal{D}^{\mathcal{N}}_{\Phi}(\mathbf{X}_{i},\mathbf{X}_{j})$ over time. In total, we have $1099$ non-missing log-returns in the given period that are common for each stock index. We divide this data into $99$ windows of size $101$ with slide step $10$, i.e. $[1,101], [11,111], \dots, [981,1081]$ and a final $100$'th window $[991,1099]$. Next, we compute $\widehat{\mathcal{D}}^{\mathcal{N}}_{\Phi,n}(\mathbf{X}_{i},\mathbf{X}_{j})$ for the returns in each window and assign this to the date corresponding to the left bound of that window. As such, we get an idea of how the dependence between groups of stock indexes across different continents evolved over time, keeping in mind that a certain date reflects the dependence calculated from the $100$ future days available in the dataset. Thanks to Theorem 1, we can also add approximated confidence bounds $\widehat{\mathcal{D}}^{\mathcal{N}}_{\Phi,n} \pm z_{1-\frac{\alpha}{2}} \widehat{\zeta}_{\Phi,n}/\sqrt{n}$ for each window, with $z_{1-\frac{\alpha}{2}}$ the $1-\frac{\alpha}{2}$ lower quantile of a standard normal distribution. 
\begin{figure}[h!] \centering
\includegraphics[scale = 0.32]{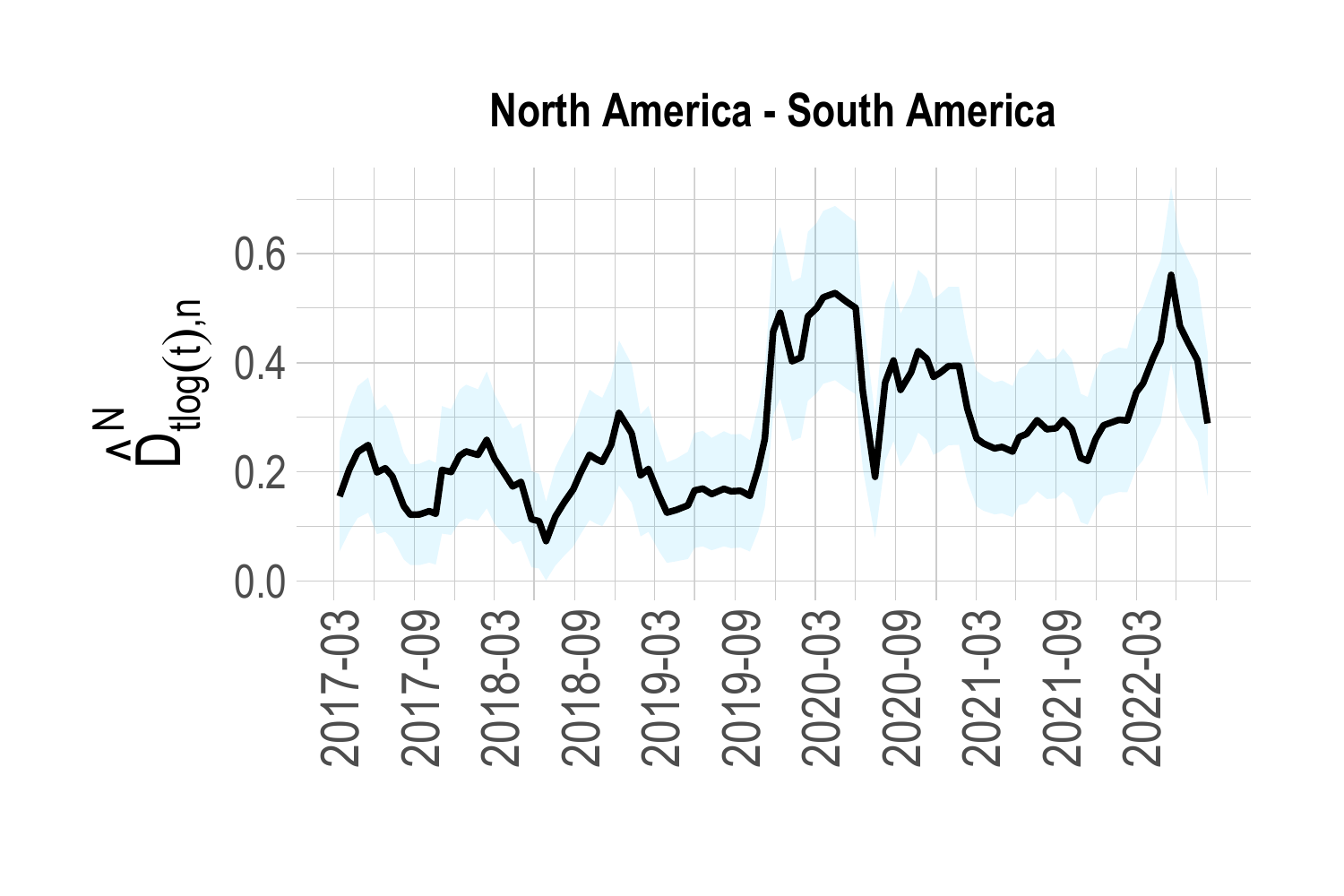}
\includegraphics[scale = 0.32]{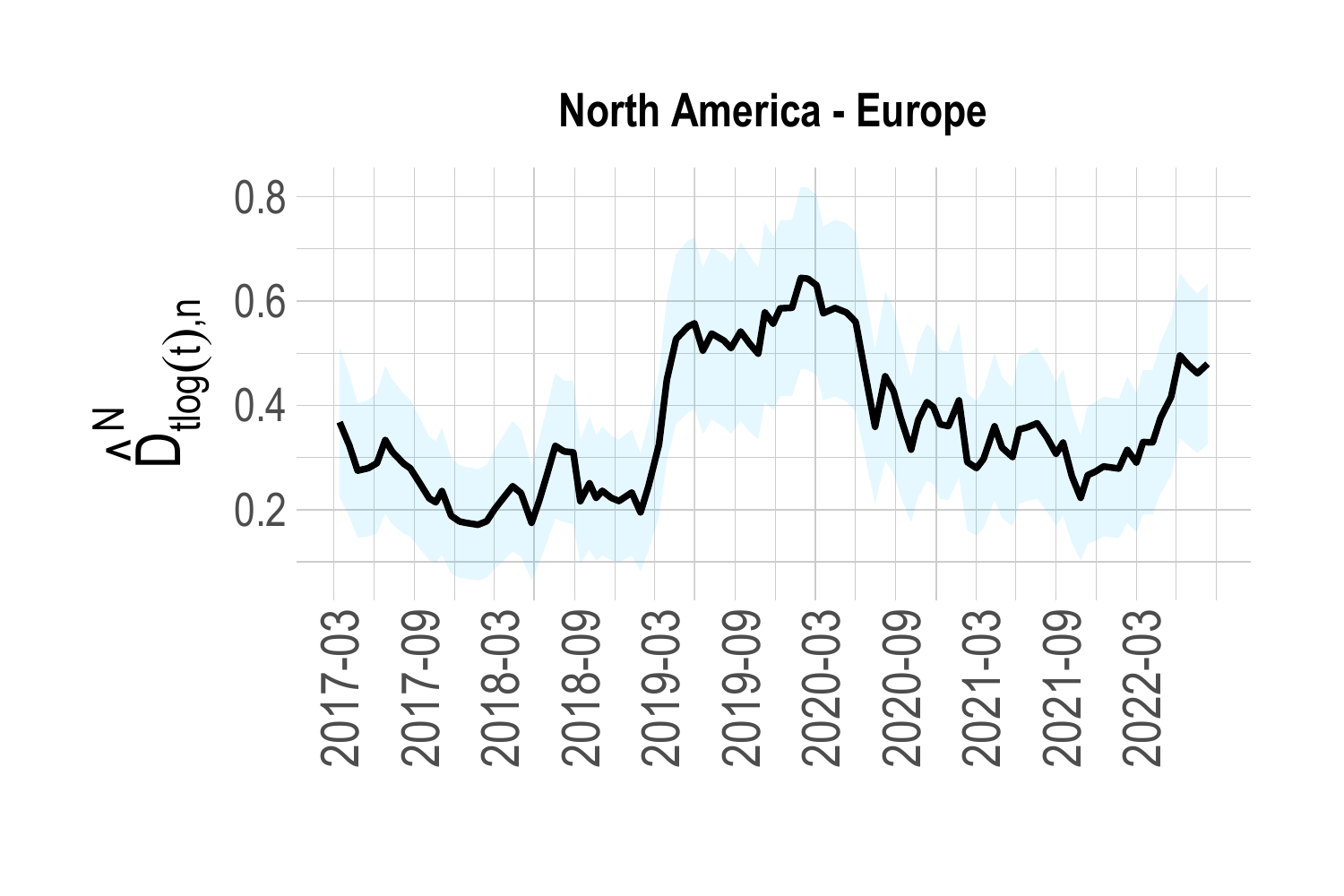}
\includegraphics[scale = 0.32]{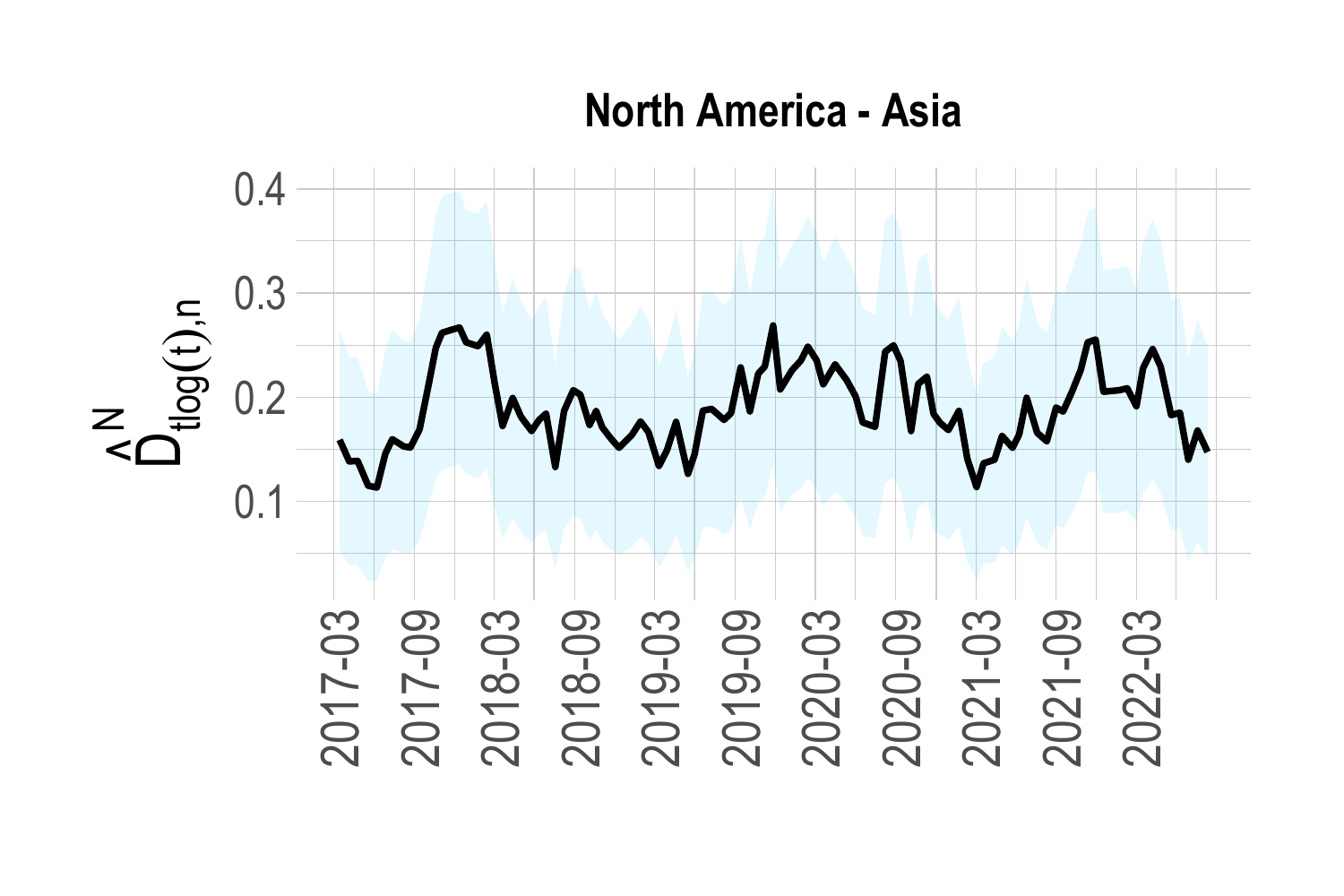}
\includegraphics[scale = 0.32]{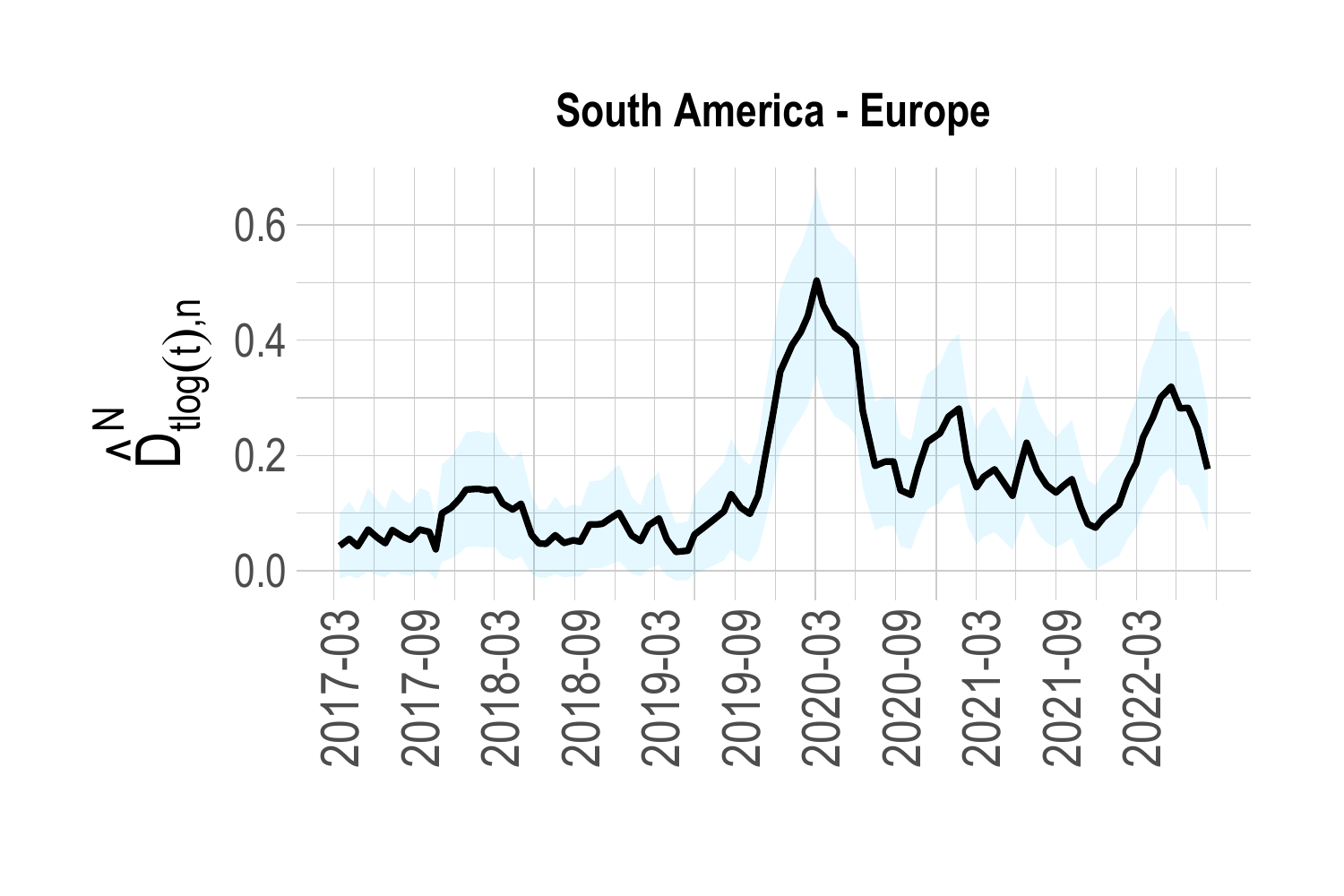}
\includegraphics[scale = 0.32]{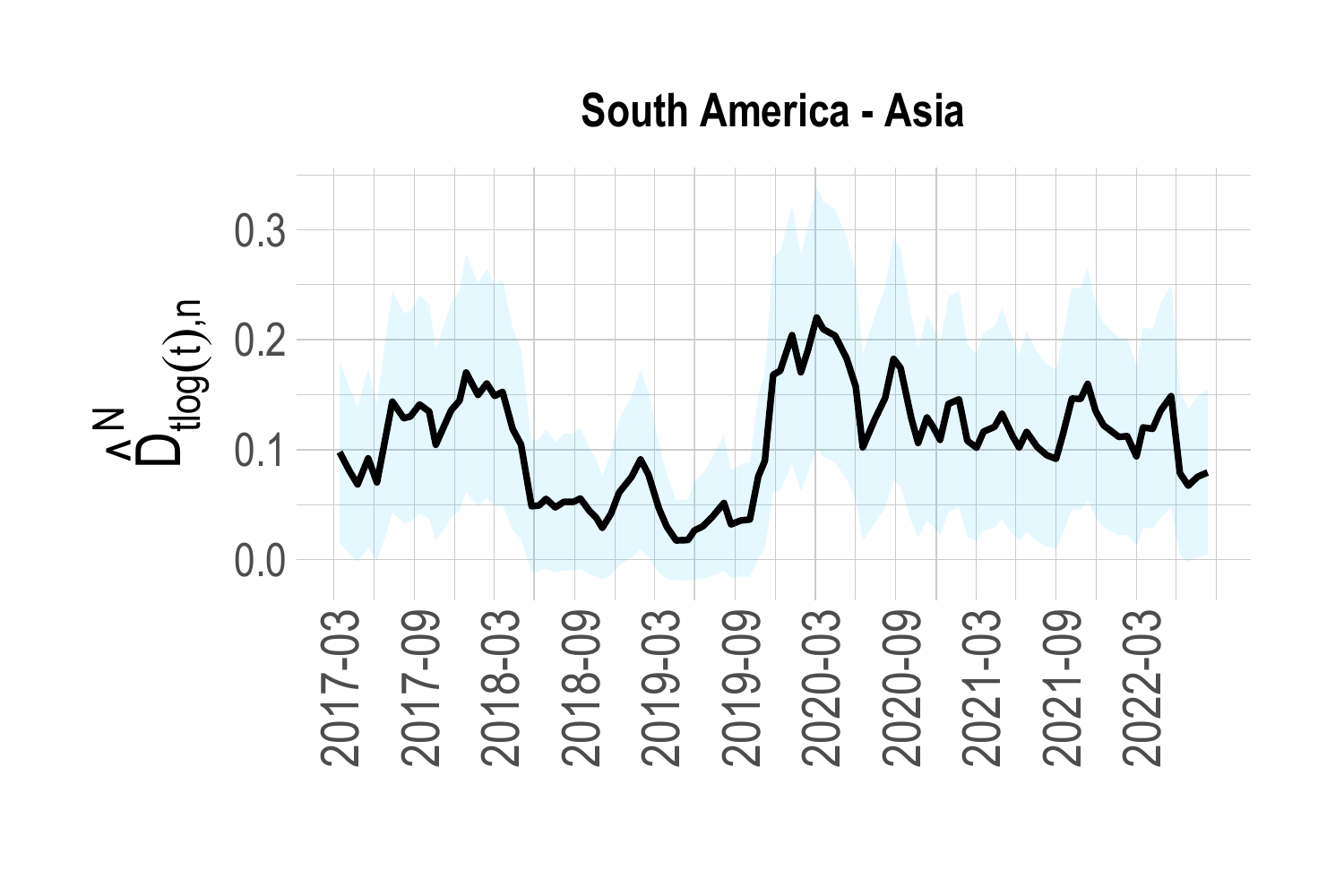}
\includegraphics[scale = 0.32]{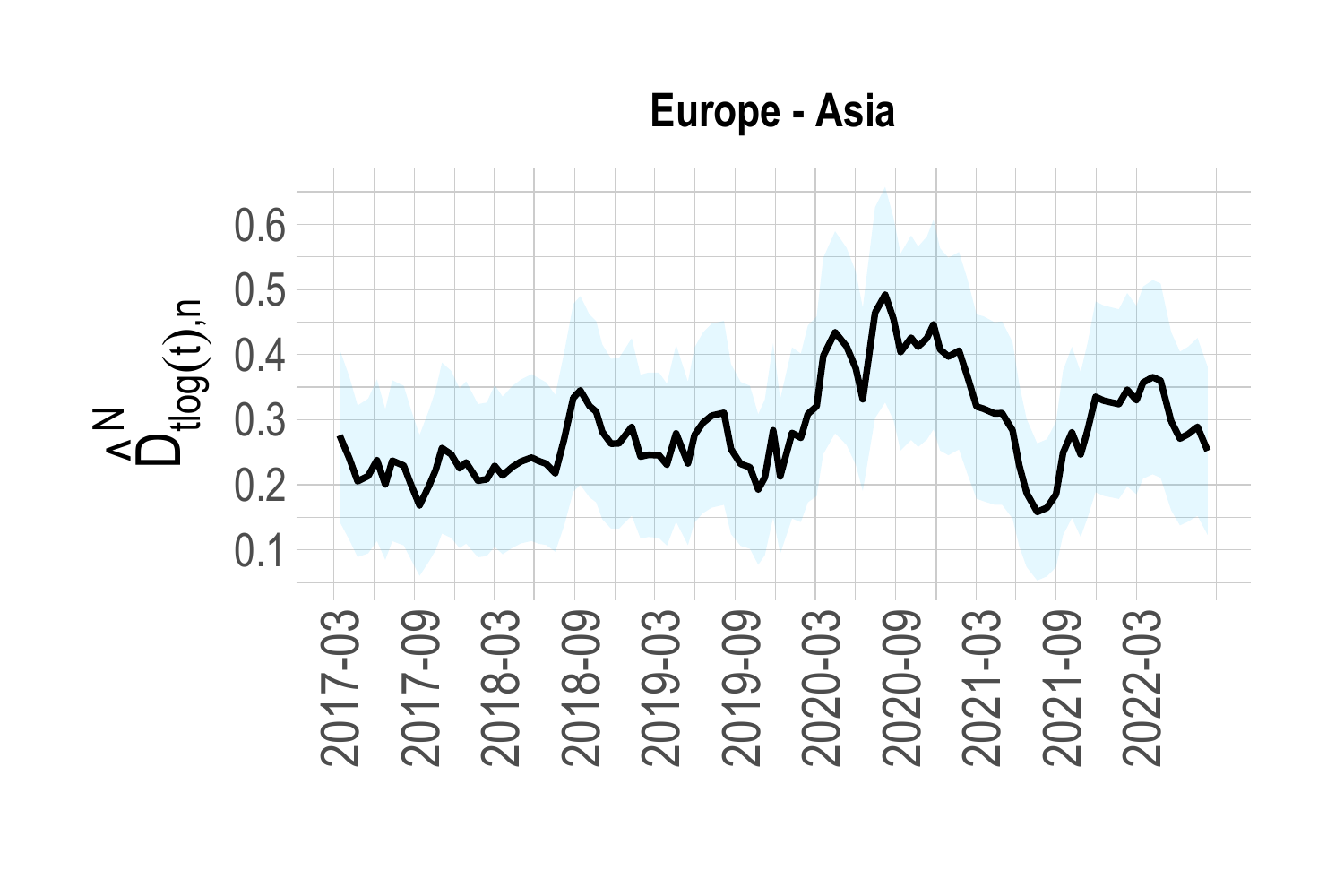}
\includegraphics[scale = 0.32]{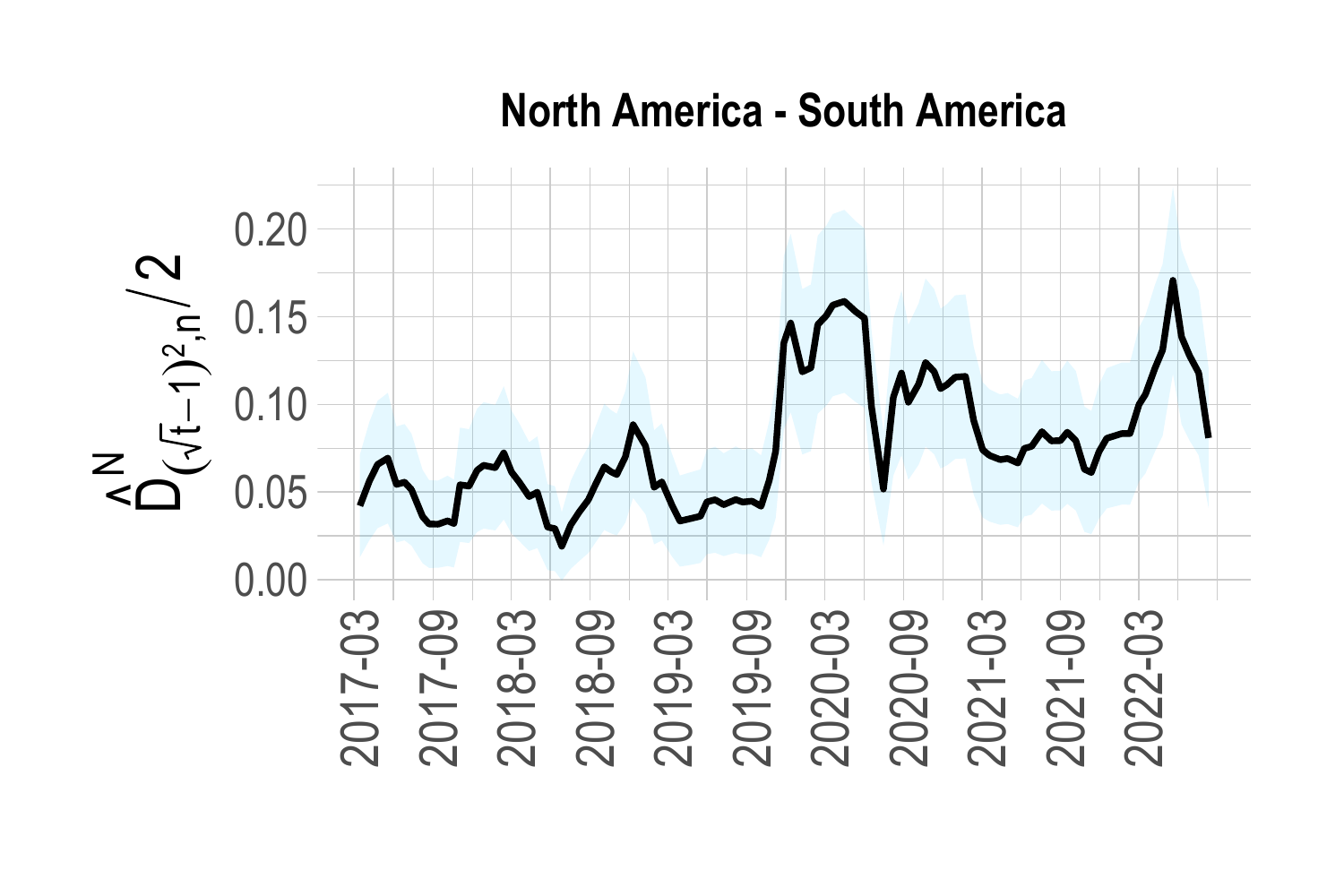}
\includegraphics[scale = 0.32]{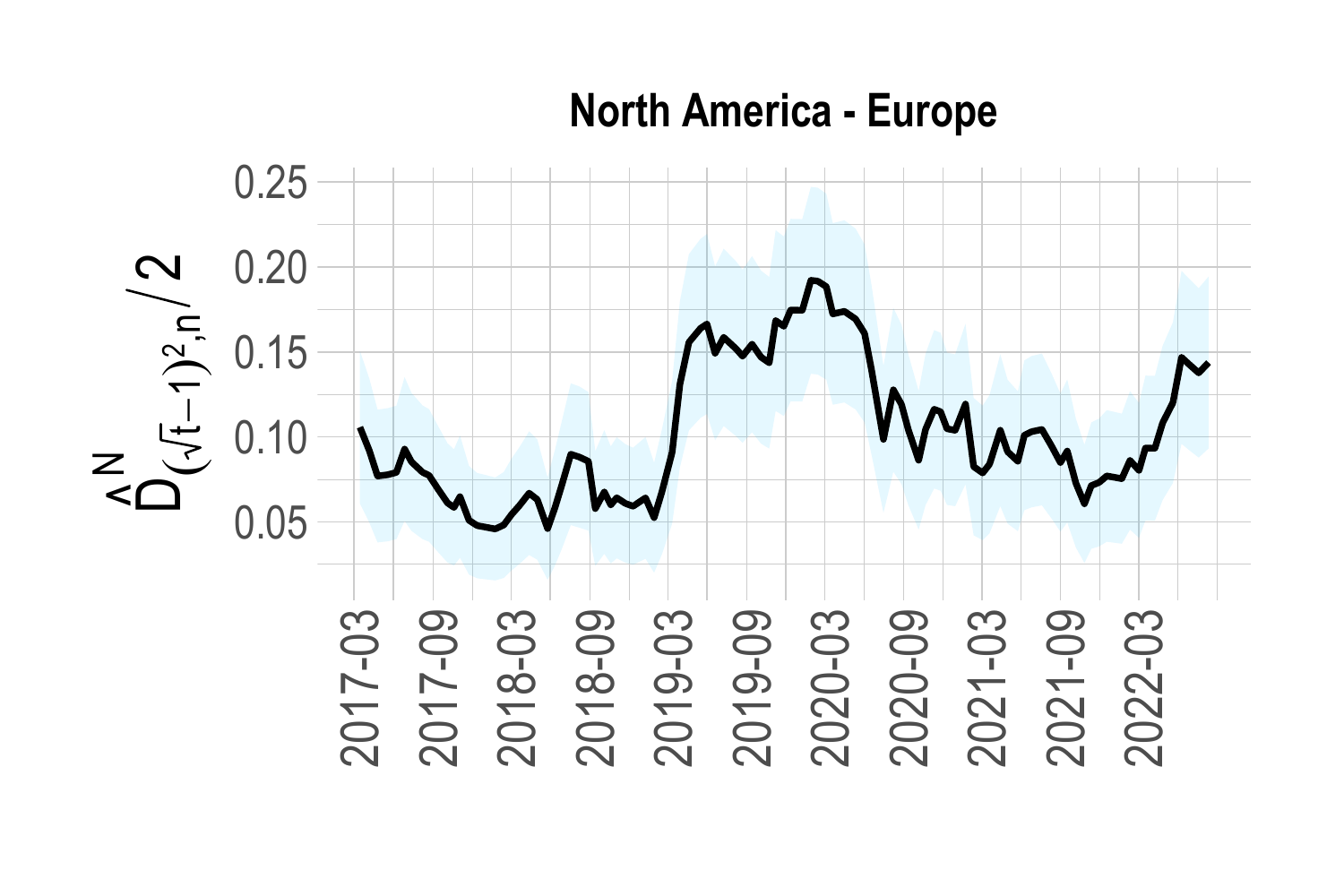}
\includegraphics[scale = 0.32]{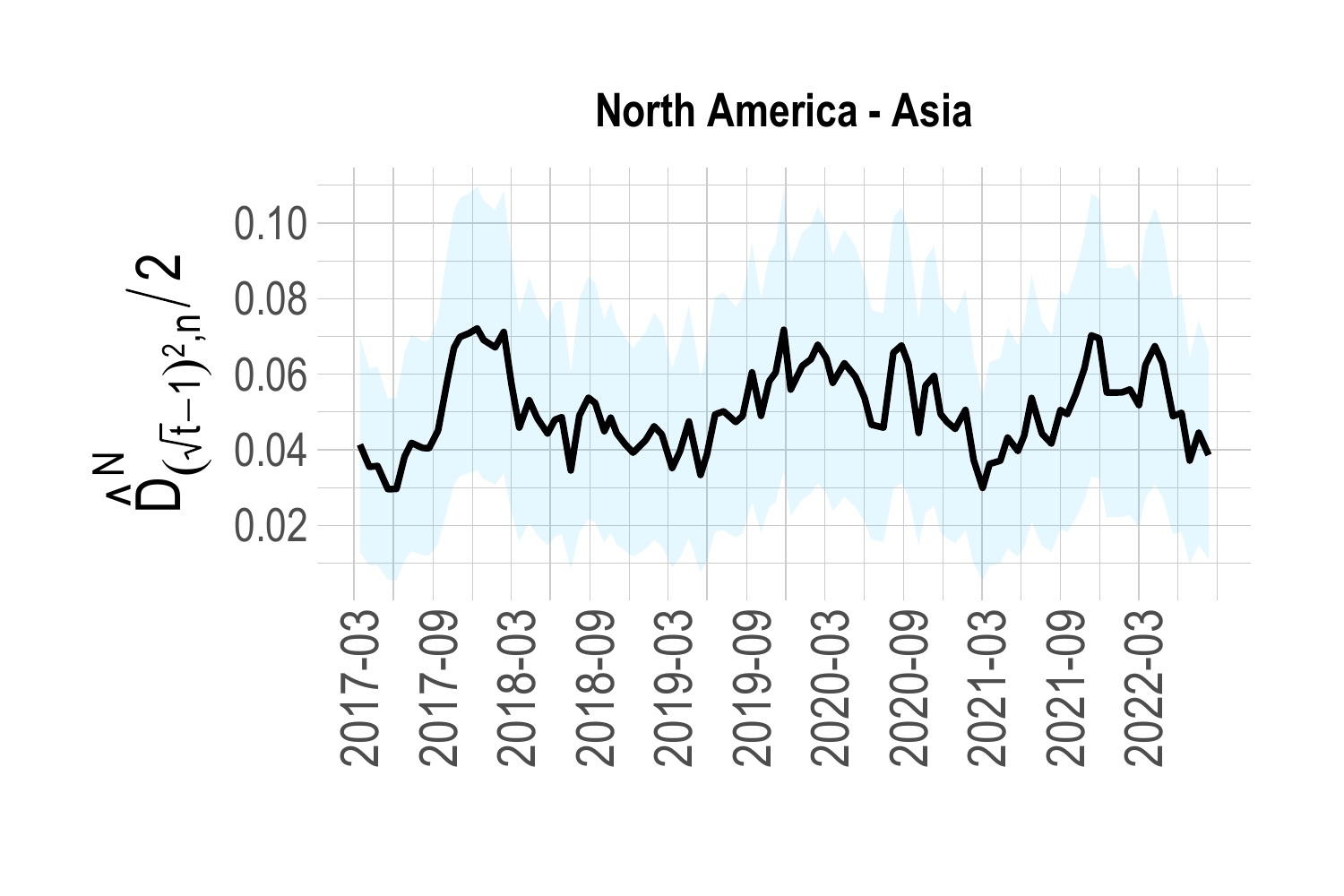}
\includegraphics[scale = 0.32]{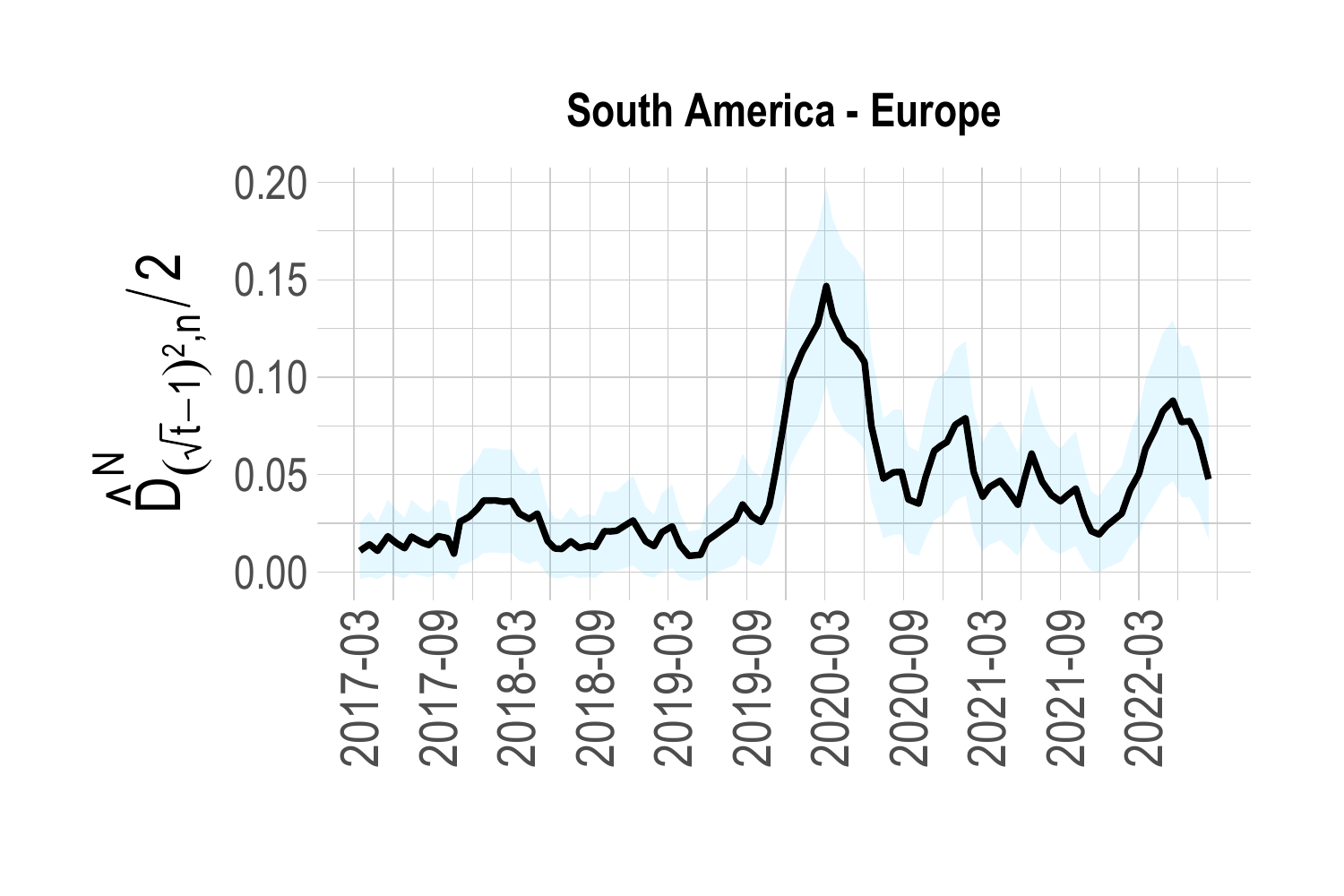}
\includegraphics[scale = 0.32]{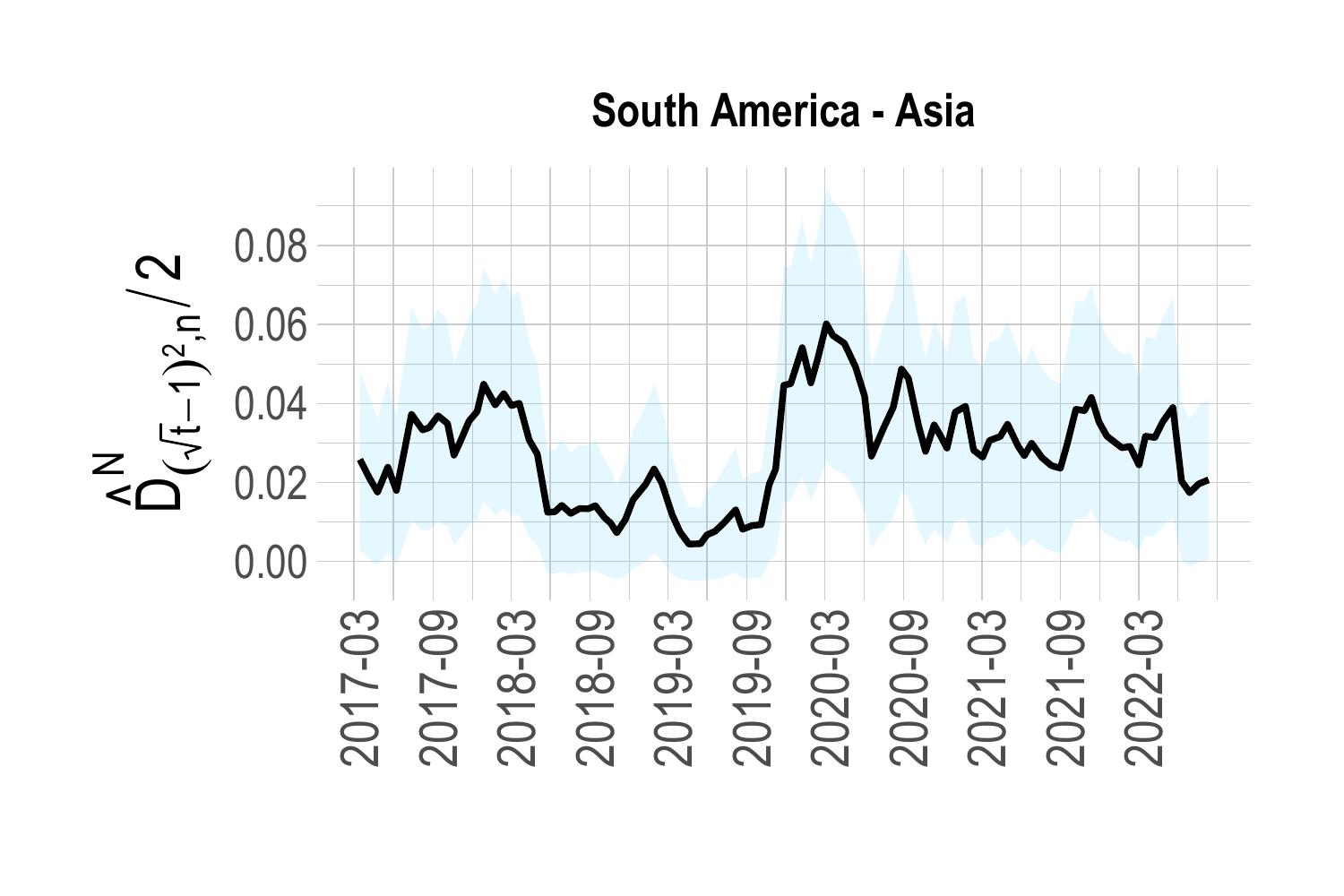}
\includegraphics[scale = 0.32]{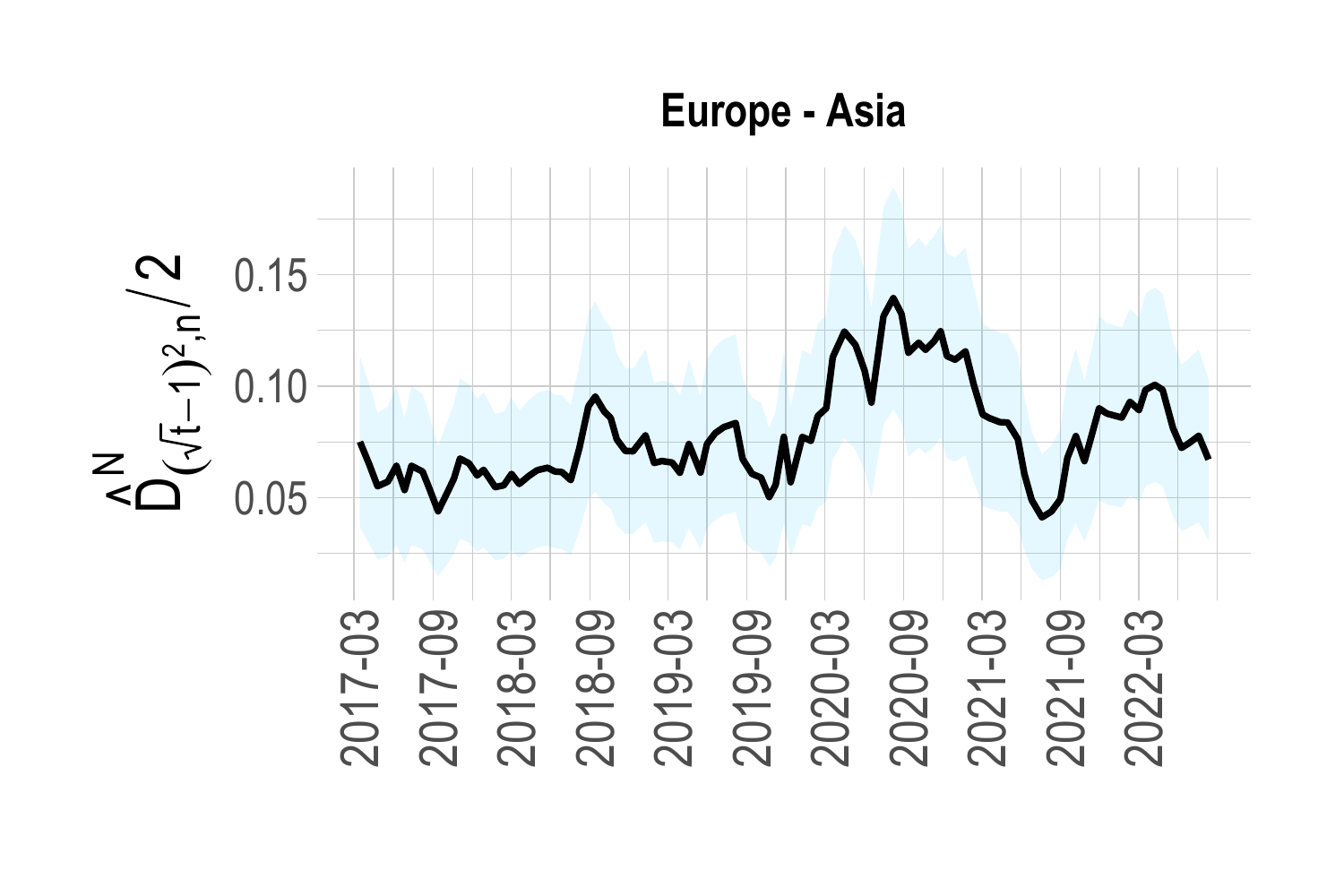}
\caption{Estimated mutual information $\widehat{\mathcal{D}}_{t \log(t),n}^{\mathcal{N}}$ and half Hellinger distance $\widehat{\mathcal{D}}_{(\sqrt{t}-1)^{2},n}^{\mathcal{N}}/2$ between groups of equity indexes of different continents over time. Regions with $95 \%$ confidence are shown in blue.}
\label{fig: FC1}
\end{figure} 

Figure \ref{fig: FC1} shows the results for the mutual information and half Hellinger distance with $\alpha = 0.05$. In all cases, we observe a (for some more pronounced than others) hump during the first year of the pandemic, attracting our attention.  In the first half of $2020$, the Corona pandemic gave rise to a stock market crash. Many indexes around the world were recovered by the end of $2020$. See \cite{Akhtaruzzaman2021} and references therein for more detailed information. We define the pre-crisis (period $1$) to be the period Dec 07, 2016 to Nov 29, 2019, the crisis period (period $2$) Dec 02, 2019 to Dec 30, 2020, and the post-crisis period (period $3$) as Jan 04, 2021 to Dec 06, 2022. For these respective periods, there are $n_{1} = 548, n_{2} = 198$ and $n_{3} = 352$ observed log-returns available. Denote $\mathcal{D}^{\mathcal{N}}_{\Phi,m}(\mathbf{X}_{i},\mathbf{X}_{j})$ with $m = 1,2,3$ for the Gaussian copula $\Phi$-dependence between $\mathbf{X}_{i}$ and $\mathbf{X}_{j}$ in period $m$, with asymptotic standard deviation $\zeta^{\mathcal{N}}_{\Phi,m}(\mathbf{X}_{i},\mathbf{X}_{j})$, and corresponding sample versions $\widehat{\mathcal{D}}^{\mathcal{N}}_{\Phi,n_{m}}(\mathbf{X}_{i},\mathbf{X}_{j})$ and $\widehat{\zeta}_{\Phi,n_{m}}(\mathbf{X}_{i},\mathbf{X}_{j})$. A test for financial contagion is
\begin{equation*}
    H_{0} : \mathcal{D}^{\mathcal{N}}_{\Phi,1}(\mathbf{X}_{i},\mathbf{X}_{j}) - \mathcal{D}^{\mathcal{N}}_{\Phi,2}(\mathbf{X}_{i},\mathbf{X}_{j}) \geq 0 \hspace{0.2cm} \text{versus} \hspace{0.2cm} H_{1} :  \mathcal{D}^{\mathcal{N}}_{\Phi,1}(\mathbf{X}_{i},\mathbf{X}_{j}) - \mathcal{D}^{\mathcal{N}}_{\Phi,2}(\mathbf{X}_{i},\mathbf{X}_{j}) < 0,
\end{equation*}
whose rejection provides statistical evidence for stronger linkages across markets during the crisis than before, and 
\begin{equation*}
    H_{0} : \mathcal{D}^{\mathcal{N}}_{\Phi,2}(\mathbf{X}_{i},\mathbf{X}_{j}) - \mathcal{D}^{\mathcal{N}}_{\Phi,3}(\mathbf{X}_{i},\mathbf{X}_{j}) \leq 0 \hspace{0.2cm} \text{versus} \hspace{0.2cm} H_{1} :  \mathcal{D}^{\mathcal{N}}_{\Phi,2}(\mathbf{X}_{i},\mathbf{X}_{j}) - \mathcal{D}^{\mathcal{N}}_{\Phi,3}(\mathbf{X}_{i},\mathbf{X}_{j}) > 0,
\end{equation*}
whose rejection illustrates weaker inter-regional dependence after the crisis than during. Asymptotic approximate $p$-values for these test are respectively given by 
\begin{equation}\label{eq: pvalue12}
    p_{12} = \mathbb{P} \left (Z_{12} \leq z_{12} \right ) \hspace{0.2cm} \text{with} \hspace{0.2cm} Z_{12} = \frac{\widehat{\mathcal{D}}^{\mathcal{N}}_{\Phi,n_{1}}(\mathbf{X}_{i},\mathbf{X}_{j}) - \widehat{\mathcal{D}}^{\mathcal{N}}_{\Phi,n_{2}}(\mathbf{X}_{i},\mathbf{X}_{j})}{\sqrt{\frac{\left ( \widehat{\zeta}_{\Phi,n_{1}}(\mathbf{X}_{i},\mathbf{X}_{j}) \right )^{2}}{n_{1}} + \frac{\left ( \widehat{\zeta}_{\Phi,n_{2}}(\mathbf{X}_{i},\mathbf{X}_{j})\right )^{2}}{n_{2}}}} \approx \mathcal{N}(0,1),
\end{equation}
and 
\begin{equation}\label{eq: pvalue23}
    p_{23} = \mathbb{P} \left (Z_{23} \geq z_{23} \right ) \hspace{0.2cm} \text{with} \hspace{0.2cm} Z_{23} = \frac{\widehat{\mathcal{D}}^{\mathcal{N}}_{\Phi,n_{2}}(\mathbf{X}_{i},\mathbf{X}_{j}) - \widehat{\mathcal{D}}^{\mathcal{N}}_{\Phi,n_{3}}(\mathbf{X}_{i},\mathbf{X}_{j})}{\sqrt{\frac{\left ( \widehat{\zeta}_{\Phi,n_{2}}(\mathbf{X}_{i},\mathbf{X}_{j}) \right )^{2}}{n_{2}} + \frac{\left ( \widehat{\zeta}_{\Phi,n_{3}}(\mathbf{X}_{i},\mathbf{X}_{j})\right )^{2}}{n_{3}}}} \approx \mathcal{N}(0,1),
\end{equation}
where $z_{12}$ and $z_{23}$ are the corresponding test values.
\begin{table}[h]
\begin{tabularx}{\textwidth}{||YYYSYYS||}
\hhline{|=======|}
\multicolumn{7}{||c||}{\normalsize \textbf{Financial contagion $p$-values}}  \\ 
& \multicolumn{2}{c}{\small Mutual information} & & \multicolumn{2}{c}{\small Half Hellinger distance} & \\ \cline{2-3} \cline{5-6}
& $p_{12}$ & $p_{23}$ & & $p_{12}$ & $p_{23}$ & \\ 
NA -  SA & $1.1657 \cdot 10^{-7}$ & $0.0008$ &  & $2.8958 \cdot 10^{-7}$ & $0.0008$ & \\ NA -  EU & $8.1439 \cdot 10^{-5}$ & $0.0039$ &  & $0.0001$ & $0.0052$ & \\ NA -  AS & $0.0151$ & $0.0672$ &  & $0.0157$ & $0.0655$ & \\ SA -  EU & $1.7442 \cdot 10^{-7}$ & $0.0012$ &  & $7.3208 \cdot 10^{-7}$ & $0.0016$ & \\ SA -  AS & $0.0004$ & $0.0114$ &  & $0.0006$ & $0.0126$ & \\ EU -  AS & $0.0016$ & $0.0388$ &  & $0.0018$ & $0.0317$ & \\
\hhline{|=======|}
\end{tabularx}
\caption{P-values \eqref{eq: pvalue12} for test of increased linkages between inter-regional equity indexes from pre-crisis to crisis and \eqref{eq: pvalue23} for test of decreased linkages from crisis to post-crisis, using the mutual information or half Hellinger distance as dependence measure (NA = North America, SA = South America, EU = Europe, AS = Asia).}
\label{tab: table4}
\end{table} 

For two continents neither of which is Asia, we have statistical evidence for inter-regional financial contagion at significance level $1\%$. When Asia is included, we find larger $p$-values, especially for $p_{23}$. 
\begin{figure}[h!] \centering
\includegraphics[width = 0.49\textwidth,height = 5cm]{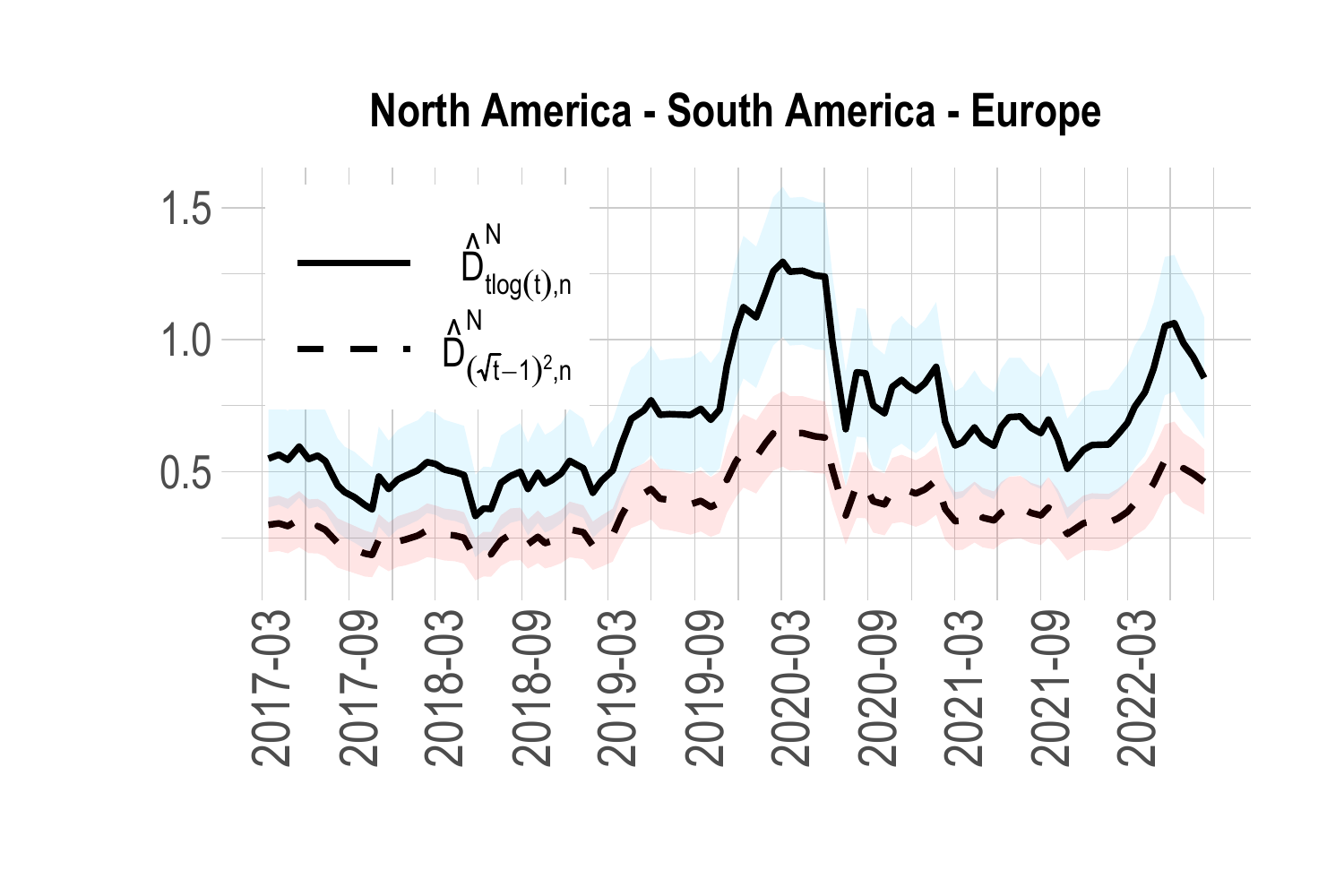}
\includegraphics[width = 0.49\textwidth,height = 5cm]{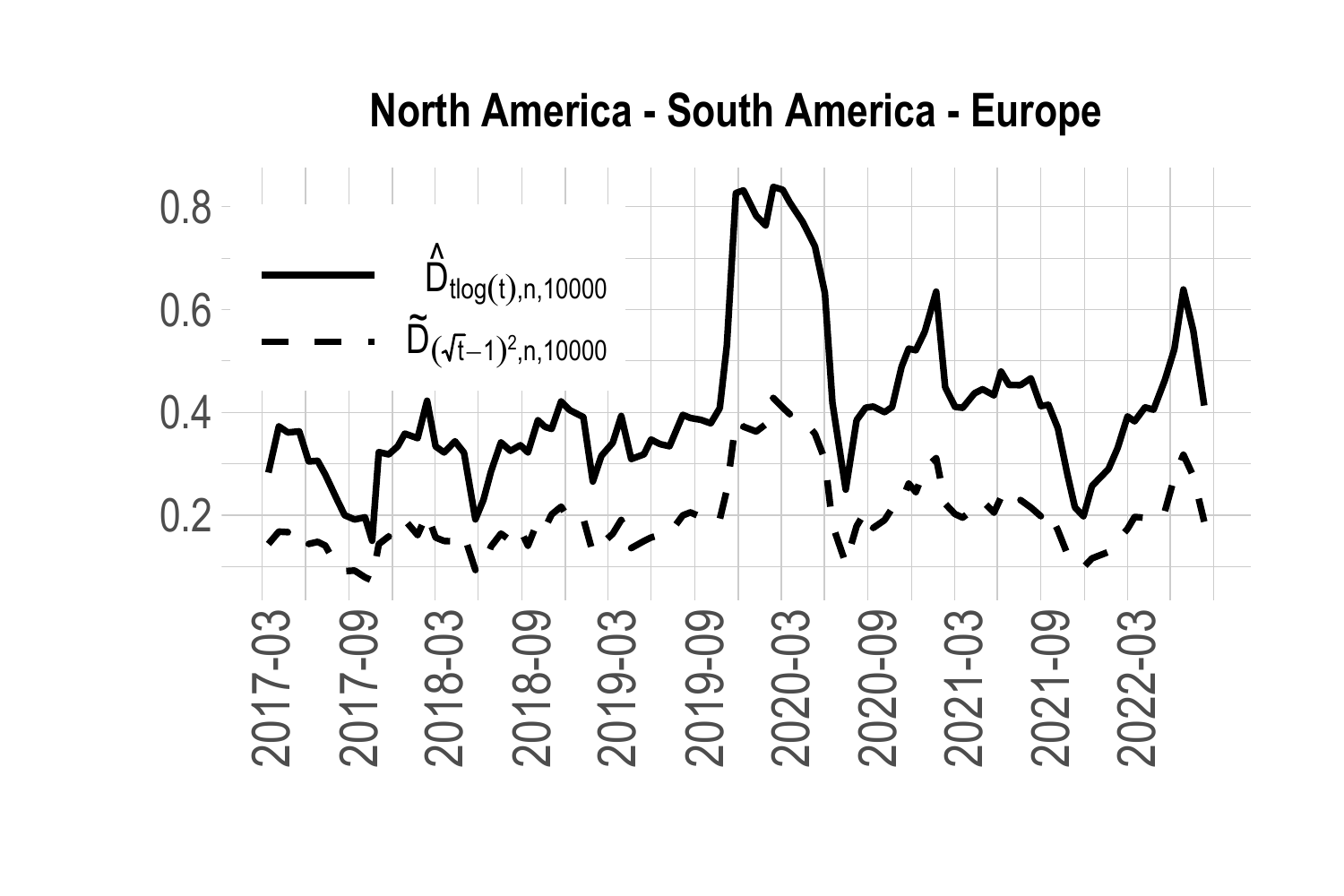}
\caption{Estimated mutual information and Hellinger distance between equity indexes of North America, South America and Europa over time, using a Gaussian copula (left) and a nested Clayton copula (right).}
\label{fig: FC2}
\end{figure} 
Over the past year, vaccination campaigns in Asia have been losing efficiency, healthcare systems have been struggling, and putting an end to the Chinese Zero-COVID policy has led to a new flare-up. 

We now look at the vector dependence between North America, South America and Europa simultaneously, not taking Asia into account. Figure \ref{fig: FC2} shows the estimated mutual information and Hellinger distance (we multiplied the half Hellinger distance with a factor $2$ in order to see the line more clearly) on the hand hand using the Gaussian copula approach (left) and on the other hand by fitting a nine dimensional nested Clayton copula (right) via the copula pseudo likelihood with non-parametrically estimated marginals using $0.1$ as starting value for the parameter of the root copula, and maximum likelihood estimates (also pseudo likelihood and non-parametric marginals) for the parameters of the child copulas, obtained with all three starting values equal to $2$. We opt for the Clayton family, as possible lower tail dependence might be incorporated as well in this copula family. Using this latter estimation approach, i.e. $\widehat{\mathcal{D}}_{t \log(t),n,M_{n}}$ or $\widetilde{\mathcal{D}}_{(\sqrt{t}-1)^{2},n,M_{n}}$, where we took $M_{n} = 10 \hspace{0.05cm} 000$, the hump of increased dependence during the COVID-19 financial recession is even more expressed. 

Finally, notice that returns can also be considered on a different periodicity, such as monthly returns. Longer horizons often allow to better capture the general trend as they are not affected so much by noise that might be present on a daily basis. For instance, as we are considering stock indexes from across the entire globe, trading hours depend on the geographic location and differences in opening and closing hours of the stock exchange might impact the dependence structure.
\newline \\ \noindent 
\textbf{8. Discussion} \newline \\ \noindent 
In this paper, we started from the general objective of quantifying dependence between a finite, arbitrary amount of random vectors. We postulated desirable properties giving preference to a copula based perspective, and proved conformity for our proposed family of $\Phi$-dependence measures. A notion of perfect dependence is not among the axioms, but a characterization of it is possible for certain choices of $\Phi$, and was illustrated by means of examples. 

Assuming a Gaussian copula model, an asymptotic normality result for the suggested plug-in estimator was established and could be interpreted further for specific $\Phi$-functions. Extensions to general parametric copula families were also obtained, focusing on maximum likelihood estimation and again a plug-in approach. Special attention went to nested Archimedean copulas, allowing for different parameters to control the intra- and inter-vector dependence. Simulations investigated finite-sample performances. In a real data analysis, the estimates indicate substantial levels of financial contagion during the COVID-19 recession. 

Among the interests for further research are high-dimensional settings ($q$ large), in which regularisation techniques are obliging. Also, semi- or non-parametric modelling of the copula can offer more flexible dependence capturing, and as such dig deeper than associations restricted to correlations (as in Gaussian copulas) and tail dependence (as in Archimedean copulas). These issues will be studied by the authors in future research.
\vspace*{0.6 cm}

\noindent
\textbf{Acknowledgement}. The authors gratefully acknowledge support from the Research Fund KU Leuven [C16/20/002 project].
\vspace*{0.24 cm}

\noindent
\bibliographystyle{myplainnat}
\bibliography{Biblio.bib}
\clearpage
\noindent \Large
\textbf{Appendix}
\newline \\ \noindent \normalsize 
\textbf{Proof of Theorem 1} \newline \\ \noindent
We first consider the specific cases $\Phi(t) = t \log(t)$ and $\Phi(t) = (\sqrt{t}-1)^{2}$. Recall that for these cases, expression \eqref{eq: phiN} reduces to \eqref{eq: mutN} and \eqref{eq: helN}, which are simpler expressions. Later in the proof, we look at the expression for $\mathbf{M}_{\Phi}$ for general $\Phi$. \newline \\ \noindent
\underline{Case $\Phi(t) = t \log(t)$} \newline \\ \noindent 
We start by showing the Fr\'{e}chet differentiability of the map
\begin{equation*}
   (\mathbb{S}^{s},||\cdot||_{\text{F}}) \rightarrow (\mathbb{R},|\cdot|) : \mathbf{\Sigma} \mapsto (\mathcal{D}_{t \log(t)}^{\mathcal{N}} \circ \varphi)(\mathbf{\Sigma}).
\end{equation*}
From Corollary 4.5 in \cite{Mordant2022}, we have that the Fr\'{e}chet derivative of $\varphi$ at $\mathbf{R}$ in the direction of a certain $\mathbf{H} \in \mathbb{S}^{q}$ equals 
\begin{equation*}
    \dot{\varphi}_{\mathbf{R}}(\mathbf{H}) = \mathbf{H} - \frac{1}{2} (\mathbf{D}_{\mathbf{H}} \mathbf{R} + \mathbf{R} \mathbf{D}_{\mathbf{H}}),
\end{equation*}
with $\mathbf{D}_{\mathbf{H}}$ the diagonal matrix containing the diagonal of $\mathbf{H}$. Let now $\mathbf{H}_{t},\mathbf{H} \in \mathbb{S}^{q}$ be such that $||\mathbf{H}_{t}-\mathbf{H}||_{\text{F}} \to 0$ as $t \to 0$. Already note that $\mathbf{R} + t \mathbf{H}_{t}$ will be in $\mathbb{S}^{q}_{>}$ for $t$ small enough since $\mathbf{R} \in \mathbb{S}^{q}_{>}$. Consider now the function 
\begin{equation*}
    f : (0,\infty)^{k+1} \rightarrow \mathbb{R} : \left (\overline{x}_{1},\dots,\overline{x}_{k},\overline{y} \right ) \mapsto - \frac{1}{2} \log \left (\frac{\overline{y}}{\prod_{i=1}^{k} \overline{x}_{i}} \right ).
\end{equation*}
Then, we see that $\mathcal{D}^{\mathcal{N}}_{t \log(t)}(\mathbf{R}) =$ $f(x_{1},\dots,x_{k},y)$ and $\mathcal{D}^{\mathcal{N}}_{t \log(t)}(\mathbf{R}+t\mathbf{H}_{t}) =$ $f(x_{1}^{t},\dots,x_{k}^{t},y^{t})$, where
\begin{equation*}
   \begin{split}
       y & = \left | \mathbf{R} \right | \hspace{1.7cm} x_{i} = \left |\mathbf{R}_{ii} \right | \\
       y^{t} & = \left |\mathbf{R} + t \mathbf{H}_{t} \right | \hspace{0.5cm} x_{i}^{t} = \left |(\mathbf{R} + t \mathbf{H}_{t})_{ii} \right |,
   \end{split} 
\end{equation*}
with $(\mathbf{R} + t \mathbf{H}_{t})_{ii}$ the $d_{i} \times d_{i}$ diagonal block of $\mathbf{R} + t \mathbf{H}_{t}$. Jacobi's formula in matrix calculus tells us that the Fr\'{e}chet derivative of $\mathbf{R} \mapsto \left |\mathbf{R} \right |$ in the direction of $\mathbf{H}$ is given by $\left |\mathbf{R} \right | \text{Tr}(\mathbf{R}^{-1}\mathbf{H})$. Hence, the Fr\'echet derivative of the map
\begin{equation*}
    g : \mathbb{S}^{q} \rightarrow (0,\infty)^{k+1} : \mathbf{R} \mapsto \left (\left |\mathbf{R}_{11} \right |,\dots,\left |\mathbf{R}_{kk} \right |,\left |\mathbf{R} \right | \right )
\end{equation*}
in the direction of $\mathbf{H}$ equals 
\begin{equation*}
    \boldsymbol{\Delta}(\mathbf{H}) = \left (\left |\mathbf{R}_{11} \right | \text{Tr} \left (\mathbf{R}_{11}^{-1} \mathbf{H}_{11} \right ),\dots,\left |\mathbf{R}_{kk} \right | \text{Tr} \left (\mathbf{R}_{kk}^{-1} \mathbf{H}_{kk} \right ), \left |\mathbf{R} \right | \text{Tr} \left (\mathbf{R}^{-1} \mathbf{H} \right ) \right ).
\end{equation*}
Furthermore, the Jacobian of $f$ is given by 
\begin{equation*}
    \mathbf{J}_{f} = \begin{pmatrix}\frac{\partial f}{\partial \overline{x}_{1}} & \cdots & \frac{\partial f}{\partial \overline{x}_{k}} & \frac{\partial f}{\partial \overline{y}}  \end{pmatrix} = \begin{pmatrix} \frac{1}{2\overline{x}_{1}} & \cdots & \frac{1}{2 \overline{x}_{k}} & -\frac{1}{2\overline{y}} \end{pmatrix},
\end{equation*}
such that the Fr\'{e}chet derivative of $f$ (being nothing more than a total derivative) at $g(\mathbf{R})$ in the direction of $\boldsymbol{\Delta}(\mathbf{H})$ is equal to
\begin{equation*}
    \mathbf{J}_{f}|_{g(\mathbf{R})}  \boldsymbol{\Delta}(\mathbf{H})^{\text{T}} = \frac{1}{2} \sum_{i=1}^{k} \text{Tr} \left (\mathbf{R}_{ii}^{-1} \mathbf{H}_{ii} \right ) - \frac{1}{2} \text{Tr} \left (\mathbf{R}^{-1} \mathbf{H} \right ) = \text{Tr} \left (\mathbf{M}_{t \log(t)} \mathbf{H} \right ).
\end{equation*}
Applying the chain rule, we have shown that the Fr\'{e}chet derivative of $\mathcal{D}_{t \log(t)}^{\mathcal{N}} \circ \varphi = f \circ g \circ \varphi$ at $\mathbf{R}$ evaluated in $\mathbf{H}$ equals 
\begin{equation*}
    \text{Tr} \left (\mathbf{M}_{t \log(t)} \dot{\varphi}_{\mathbf{R}}(\mathbf{H}) \right ) = \text{Tr} \left ( \left (\mathbf{M}_{t \log(t)} - \mathbf{D}_{\mathbf{M}_{t \log(t)}\mathbf{R}} \right ) \mathbf{H} \right ),
\end{equation*}
where the last equality follows from the exact same arguments as in Corollary 4.5 of \cite{Mordant2022}. Notice that this Fr\'{e}chet derivative is linear, as it should be. 
\newline \\ \noindent
\underline{Case $\Phi(t) = (\sqrt{t}-1)^{2}$} \newline \\ \noindent
We continue by showing the Fr\'{e}chet differentiability of $\mathbf{\Sigma} \mapsto (\mathcal{D}_{(\sqrt{t}-1)^{2}}^{\mathcal{N}} \circ \varphi)(\mathbf{\Sigma})$. Therefore, note that 
\begin{equation*}
    \mathcal{D}_{(\sqrt{t}-1)^{2}}^{\mathcal{N}}(\mathbf{R}) = 2 - 2 \frac{2^{q/2}}{\left |\mathbb{I}_{q}+\mathbf{R}_{0}^{-1}\mathbf{R} \right |^{1/2}} \exp \left (-\frac{1}{2} \mathcal{D}^{\mathcal{N}}_{t \log(t)}
    (\mathbf{R})\right ).
\end{equation*}
By the chain rule, the Fr\'{e}chet derivative of $\exp((-1/2)\mathcal{D}^{\mathcal{N}}_{t \log(t)})$ at $\mathbf{R}$ in the direction of $\mathbf{H}$ equals 
\begin{equation*}
    -\frac{1}{4} \exp \left (-\frac{1}{2} \mathcal{D}^{\mathcal{N}}_{t \log(t)}(\mathbf{R})  \right ) \text{Tr} \left (\left (\mathbf{R}_{0}^{-1} - \mathbf{R}^{-1} \right ) \mathbf{H} \right ).
\end{equation*}
From basic matrix calculus, the derivative of $\mathbf{R} \mapsto \mathbf{R}_{0}^{-1}$ at $\mathbf{R}$ evaluated in $\mathbf{H}$ is
\begin{equation*}
    \text{diag} \left (-\mathbf{R}_{11}^{-1} \mathbf{H}_{11} \mathbf{R}_{11}^{-1},\dots,-\mathbf{R}_{kk}^{-1}\mathbf{H}_{kk}\mathbf{R}_{kk}^{-1} \right ).
\end{equation*}
Obviously $\mathbf{R} \mapsto \mathbf{R}$ has the identity $\mathbf{H}$ as derivative, and the product rule gives
\begin{equation*}
    \widetilde{\boldsymbol{\Delta}}(\mathbf{H}) = \text{diag} \left (-\mathbf{R}_{11}^{-1} \mathbf{H}_{11} \mathbf{R}_{11}^{-1},\dots,-\mathbf{R}_{kk}^{-1}\mathbf{H}_{kk}\mathbf{R}_{kk}^{-1} \right )\mathbf{R} + \mathbf{R}_{0}^{-1} \mathbf{H}
\end{equation*}
as Fr\'{e}chet derivative of $\mathbf{R} \mapsto \mathbb{I}_{q} + \mathbf{R}_{0}^{-1} \mathbf{R}$. The derivative of $\mathbf{R} \mapsto |\mathbb{I}_{q} + \mathbf{R}_{0}^{-1} \mathbf{R}|$ at $\mathbf{R}$ in the direction of $\mathbf{H}$, is the derivative of $\mathbf{R} \mapsto |\mathbf{R}|$ at $\mathbb{I}_{q} + \mathbf{R}_{0}^{-1}\mathbf{R}$ in the direction of $\widetilde{\boldsymbol{\Delta}}(\mathbf{H})$, i.e.
\begin{equation*}
    \left |\mathbb{I}_{q} + \mathbf{R}_{0}^{-1} \mathbf{R} \right | \text{Tr} \left (\left (\mathbb{I}_{q} + \mathbf{R}_{0}^{-1} \mathbf{R} \right ) \widetilde{\boldsymbol{\Delta}}(\mathbf{H}) \right ).
\end{equation*}
From this, it is easily seen that 
\begin{equation*}
    -\frac{2^{q/2}}{2} \left |\mathbb{I}_{q} + \mathbf{R}_{0}^{-1} \mathbf{R} \right |^{-1/2} \text{Tr} \left (\left (\mathbb{I}_{q} + \mathbf{R}_{0}^{-1} \mathbf{R} \right ) \widetilde{\boldsymbol{\Delta}}(\mathbf{H}) \right )
 \end{equation*}
is the derivative of $2^{q/2} / \left |\mathbb{I}_{q} + \mathbf{R}_{0}^{-1} \mathbf{R} \right |^{1/2}$, and, by using the product rule again,
\begin{equation*}
    \frac{2^{q/2} \exp \left (-\frac{1}{2} \mathcal{D}^{\mathcal{N}}_{t \log(t)} \right )}{\left |\mathbb{I}_{q} + \mathbf{R}_{0}^{-1} \mathbf{R} \right |^{1/2}} \left [\frac{1}{2} \text{Tr} \left (\left (\mathbf{R}_{0}^{-1} - \mathbf{R}^{-1} \right ) \mathbf{H} \right ) + \text{Tr} \left (\left (\mathbb{I}_{q} + \mathbf{R}_{0}^{-1}\mathbf{R} \right )^{-1} \widetilde{\boldsymbol{\Delta}}(\mathbf{H}) \right ) \right ]
\end{equation*}
is the Fr\'{e}chet derivative of $\mathcal{D}^{\mathcal{N}}_{(\sqrt{t}-1)^{2}}$ at $\mathbf{R}$ evaluated in $\mathbf{H}$. Notice that the above derivative is not yet specifically of the from $\text{Tr}(\mathbf{M}\mathbf{H})$ for a certain matrix $\mathbf{M}$. However, since all linear maps are of that form, and Fr\'{e}chet derivatives are linear maps, we should be able to find such $\mathbf{M}$. Let 
\begin{equation*}
    \mathbf{P}_{i} = \big (\mathbf{0}_{d_{i} \times d_{1}} \cdots \hspace{0.1cm} \mathbf{0}_{d_{i} \times d_{i-1}} \hspace{0.1cm}  \pmb{\mathbb{I}}_{d_{i}} \hspace{0.1cm}  \mathbf{0}_{d_{i} \times d_{i+1}} \cdots \hspace{0.1cm}  \mathbf{0}_{d_{i} \times d_{k}} \big ) \in \mathbb{R}^{d_{i} \times q}
\end{equation*}
be the projection matrix onto the $d_{i}$ coordinates, satisfying $\mathbf{R}_{ii} = \mathbf{P}_{i} \mathbf{R} \mathbf{P}_{i}^{\text{T}}$, and
\begin{equation*}
    \mathbf{D} = \text{diag} \left (-\mathbf{R}_{11}^{-1} \mathbf{H}_{11} \mathbf{R}_{11}^{-1},\dots,-\mathbf{R}_{kk}^{-1}\mathbf{H}_{kk}\mathbf{R}_{kk}^{-1} \right ).
\end{equation*}
Observe that
\begin{equation*}
\begin{split}
    \text{Tr} \left (\left (\mathbb{I}_{q} + \mathbf{R}_{0}^{-1} \mathbf{R} \right )^{-1} \mathbf{D} \mathbf{R} \right ) = \text{Tr} \left (\mathbf{J} \mathbf{D} \right ) & = - \sum_{i=1}^{k} \text{Tr} \left (\mathbf{J}_{ii} \mathbf{R}_{ii}^{-1} \mathbf{H}_{ii} \mathbf{R}_{ii}^{-1} \right ) \\
    & = - \sum_{i=1}^{k} \text{Tr} \left (\mathbf{J}_{ii} \mathbf{R}_{ii}^{-1} \mathbf{P}_{i} \mathbf{H} \mathbf{P}_{i}^{\text{T}} \mathbf{R}_{ii}^{-1} \right ) \\
    & = - \sum_{i=1}^{k} \text{Tr} \left (\mathbf{P}_{i}^{\text{T}} \mathbf{R}_{ii}^{-1} \mathbf{J}_{ii} \mathbf{R}_{ii}^{-1} \mathbf{P}_{i} \mathbf{H} \right ) \\
    & = - \text{Tr} \left (\boldsymbol{\gamma \mathbf{H}} \right ),
\end{split}
\end{equation*}
using the cyclic trace property and the fact that 
\begin{equation*}
    \boldsymbol{\gamma} = \sum_{i=1}^{k} \mathbf{P}_{i}^{\text{T}} \mathbf{R}_{ii}^{-1} \mathbf{J}_{ii} \mathbf{R}_{ii}^{-1} \mathbf{P}_{i}.
\end{equation*}
Knowing this, it is quickly seen that the Fr\'{e}chet derivative of $\mathcal{D}^{\mathcal{N}}_{(\sqrt{t}-1)^{2}} \circ \varphi$ at $\mathbf{R}$ in the direction of $\mathbf{H}$ equals 
\begin{equation*}
    \text{Tr} \left (\left (\mathbf{M}_{(\sqrt{t}-1)^{2}} - \mathbf{D}_{\mathbf{M}_{(\sqrt{t}-1)^{2}}\mathbf{R}} \right )\mathbf{H} \right ).
\end{equation*}
\underline{General case} \newline \\ \noindent
We now look at the Fr\'{e}chet derivative of the integrand of $\mathcal{D}_{\Phi}^{\mathcal{N}}(\mathbf{R})$ in \eqref{eq: phiN} for general $\Phi$. The derivative of $\mathbf{R} \mapsto \mathbf{R}_{0}^{-1}$ in the direction of $\mathbf{H}$ equals
\begin{equation*}
    \mathbf{D} = \text{diag} \left (-\mathbf{R}_{11}^{-1}\mathbf{H}_{11}\mathbf{R}_{11}^{-1},\dots,-\mathbf{R}_{kk}^{-1}\mathbf{H}_{kk}\mathbf{R}_{kk}^{-1} \right ).
\end{equation*}
By the chain rule, the derivative of $\mathbf{R} \mapsto \exp((-1/2) \mathbf{x}^{\text{T}}\mathbf{R}_{0}^{-1}\mathbf{x})$ evaluated at $\mathbf{H}$ is
\begin{equation*}
    -\frac{1}{2} \exp \left (- \frac{1}{2} \mathbf{x}^{\text{T}} \mathbf{R}_{0}^{-1} \mathbf{x} \right ) \mathbf{x}^{\text{T}} \mathbf{D} \mathbf{x}. 
\end{equation*}
The derivative of $\mathbf{R} \mapsto |\mathbf{R}|^{1/2}$ in the direction of $\mathbf{H}$ equals
\begin{equation*}
    \frac{1}{2} \left |\mathbf{R} \right |^{1/2} \text{Tr} \left (\mathbf{R}^{-1}\mathbf{H} \right ).
\end{equation*}
Next, since the Jacobian of the function $f(\overline{x}_{1},\dots,\overline{x}_{k},\overline{y}) = \overline{y}/(\overline{x}_{1} \cdots \overline{x}_{k})$ is
\large 
\begin{equation*}
    \begin{pmatrix} \frac{-\overline{y}}{\overline{x}_{1}^{2} \cdot \overline{x}_{2} \cdots \overline{x}_{k}} & \cdots &  \frac{-\overline{y}}{\overline{x}_{1} \cdot \overline{x}_{2} \cdots \overline{x}_{k}^{2}} & \frac{1}{\overline{x}_{1} \cdots \overline{x}_{k}} \end{pmatrix},
\end{equation*} \normalsize 
it is quickly seen that the Fr\'{e}chet derivative of the factor in the integrand of \eqref{eq: phiN} in front of $\Phi$ in the direction of $\mathbf{H}$ equals
\begin{equation*}
    -\frac{\exp \left (-\frac{1}{2} \mathbf{x}^{\text{T}}\mathbf{R}_{0}^{-1} \mathbf{x} \right )}{(2\pi)^{q/2}2\prod_{i=1}^{k} \left |\mathbf{R}_{ii} \right |^{1/2}} \left [\text{Tr} \left (\mathbf{R}_{0}^{-1}\mathbf{H} \right ) + \mathbf{x}^{\text{T}} \mathbf{D} \mathbf{x}  \right ].
\end{equation*}
Regarding the term within the $\Phi$-function, note that the derivative of $\mathbf{R} \mapsto \mathbf{R}^{-1} - \mathbf{R}_{0}^{-1}$ evaluated in $\mathbf{H}$ equals 
\begin{equation*}
    -\mathbf{R}^{-1} \mathbf{H} \mathbf{R}^{-1} - \mathbf{D},
\end{equation*}
such that, using a similar reasoning as before
\begin{equation*}
    \text{\scalebox{0.95}{$\frac{\prod_{i=1}^{k}\left |\mathbf{R}_{ii} \right |^{1/2} \exp \left (-\frac{1}{2} \mathbf{x}^{\text{T}} \left (\mathbf{R}^{-1} - \mathbf{R}_{0}^{-1} \right ) \mathbf{x} \right )}{2 \left |\mathbf{R} \right |^{1/2}} \left [\text{Tr} \left (\left (\mathbf{R}_{0}^{-1} - \mathbf{R}^{-1} \right )\mathbf{H} \right ) + \mathbf{x}^{\text{T}} \left (\mathbf{R}^{-1} \mathbf{H} \mathbf{R}^{-1} + \mathbf{D} \right )\mathbf{x} \right ]$}}
\end{equation*}
is the Fr\'{e}chet derivative of the term inside the $\Phi$-function in the direction of $\mathbf{H}$. Denote now $f_{\mathbf{R}_{0}}(\mathbf{x})$ for the density function of a $\mathcal{N}(\mathbf{0},\mathbf{R}_{0})$ distribution, and similarly $f_{\mathbf{R}}(\mathbf{x})$. The product and chain rule tell us that the Fr\'{e}chet derivative of the integrand of $\mathcal{D}^{\mathcal{N}}_{\Phi}(\mathbf{R})$ in the direction of $\mathbf{H}$ is given by 
\begin{equation*}
\begin{split}
    \mathcal{I}(\mathbf{x}) = & -\frac{1}{2} f_{\mathbf{R}_{0}}(\mathbf{x}) \text{Tr} \left (\mathbf{R}_{0}^{-1} \mathbf{H} \right ) \alpha(\mathbf{x})  - \frac{1}{2} f_{\mathbf{R}_{0}}(\mathbf{x})\mathbf{x}^{\text{T}}\mathbf{D}\mathbf{x} \alpha(\mathbf{x}) \\ & + \frac{1}{2} f_{\mathbf{R}}(\mathbf{x}) \text{Tr} \left (\left (\mathbf{R}_{0}^{-1}-\mathbf{R}^{-1} \right ) \mathbf{H} \right ) \alpha^{\prime}(\mathbf{x})  + \frac{1}{2} f_{\mathbf{R}}(\mathbf{x}) \mathbf{x}^{\text{T}} \left (\mathbf{R}^{-1}\mathbf{H}\mathbf{R}^{-1} + \mathbf{D} \right )\mathbf{x} \alpha^{\prime}(\mathbf{x}).
\end{split}
\end{equation*}
Hence,
\begin{equation*}
\begin{split}
    & \int_{\mathbb{R}^{q}} \mathcal{I}(\mathbf{x}) d\mathbf{x} = -\frac{1}{2} \mathbb{E}_{\mathcal{N}(\mathbf{0},\mathbf{R}_{0})} \left [\alpha(\mathbf{X}) \right ] \text{Tr} \left (\mathbf{R}_{0}^{-1} \mathbf{H} \right ) - \frac{1}{2} \mathbb{E}_{\mathcal{N}(\mathbf{0},\mathbf{R}_{0})} \left [\mathbf{X}^{\text{T}} \mathbf{D} \mathbf{X} \alpha \left (\mathbf{X} \right ) \right] \\ & + \frac{1}{2} \mathbb{E}_{\mathcal{N}(\mathbf{0},\mathbf{R})} \left [\alpha^{\prime}(\mathbf{X}) \right ] \text{Tr} \left (\left (\mathbf{R}_{0}^{-1}-\mathbf{R}^{-1} \right ) \mathbf{H} \right ) + \frac{1}{2} \mathbb{E}_{\mathcal{N}(\mathbf{0},\mathbf{R})} \left [\mathbf{X}^{\text{T}} \left (\mathbf{R}^{-1}\mathbf{H}\mathbf{R}^{-1} + \mathbf{D} \right )\mathbf{X} \alpha^{\prime}(\mathbf{X}) \right ].
\end{split}
\end{equation*}
First of all, by definition 
\begin{equation*}
    \mathbb{E}_{\mathcal{N}(\mathbf{0},\mathbf{R}_{0})} \left [\alpha(\mathbf{X}) \right ] = \mathcal{D}_{\Phi}^{\mathcal{N}}(\mathbf{R}).
\end{equation*}
Furthermore, since a quadratic form is just a number, we also have that
\begin{equation}
\begin{split}
    \mathbb{E}_{\mathcal{N}(\mathbf{0},\mathbf{R}_{0})} \left [\mathbf{X}^{\text{T}} \mathbf{D} \mathbf{X} \alpha(\mathbf{X}) \right ]  & = \mathbb{E}_{\mathcal{N}(\mathbf{0},\mathbf{R}_{0})} \left [\text{Tr} \left (\mathbf{X}^{\text{T}} \mathbf{D} \mathbf{X} \alpha(\mathbf{X}) \right )\right ] \\ & = \mathbb{E}_{\mathcal{N}(\mathbf{0},\mathbf{R}_{0})} \left [\text{Tr} \left (\alpha(\mathbf{X}) \mathbf{X} \mathbf{X}^{\text{T}} \mathbf{D} \right )\right ] \\ & = \text{Tr} \left (\mathbb{E}_{\mathcal{N}(\mathbf{0},\mathbf{R}_{0})} \left [\alpha(\mathbf{X}) \mathbf{X}\mathbf{X}^{\text{T}} \right ] \mathbf{D} \right ) \\
    & = -\text{Tr} \left (\mathbf{F}_{1} \mathbf{H}\right ),
\end{split}
\tag{A1}
\end{equation}
where we used the linearity and cyclic property of the trace operator, and a similar reasoning with the projection matrices $\mathbf{P}_{i}$ as in the case $\Phi(t) = (\sqrt{t}-1)^{2}$. The other terms can be handled in a similar way, resulting in 
\begin{equation*}
    \int_{\mathbb{R}^{q}} \mathcal{I}(\mathbf{x}) d\mathbf{x} = \text{Tr} \left (\mathbf{M}_{\Phi}\mathbf{H} \right ).
\end{equation*}
\\ \noindent
\underline{Applying the delta method} \newline \\ \noindent 
Next, we consider the estimator $\widehat{\mathbf{R}}_{n}$. Theorem 3.1 in \cite{Klaassen1997} tells us that
\begin{equation*} 
    \text{\scalebox{0.86}{$\sqrt{n}(\widehat{\mathbf{R}}_{n} - \mathbf{R}) - \frac{1}{\sqrt{n}} \sum_{\ell=1}^{n} \left [\mathbf{Z}^{(\ell)} \left (\mathbf{Z}^{(\ell)} \right )^{\text{T}} - \frac{1}{2} \left (\mbox{diag} \left (\mathbf{Z}^{(\ell)} \left (\mathbf{Z}^{(\ell)} \right )^{\text{T}} \right )\mathbf{R} + \mathbf{R} \hspace{0.05cm} \text{diag}\left (\mathbf{Z}^{(\ell)} \left (\mathbf{Z}^{(\ell)} \right )^{\text{T}} \right )  \right )\right ] \xrightarrow{p} \mathbf{0}_{q \times q}$}}
\end{equation*}
as $n \to \infty$, where $\mathbf{Z}^{(\ell)} = (\mathbf{Z}_{1}^{(\ell)},\dots,\mathbf{Z}_{k}^{(\ell)})^{\text{T}}$, with $\mathbf{Z}_{i}^{(\ell)} = (Z_{i1}^{(\ell)},\dots,Z_{id_{i}}^{(\ell)})$ for $i = 1,\dots,k$, for $\ell = 1,\dots,n$ is a sample from the $\mathcal{N}_{q}(\mathbf{0}_{q},\mathbf{R})$ distribution. The same expansion holds when $\widehat{\mathbf{R}}_{n}$ is the empirical correlation matrix of $\mathbf{Z}^{(1)},\dots,\mathbf{Z}^{(n)}$, see e.g. Lemma 4.17 in \cite{Mordant2022}. Hence 
\begin{equation*}
    \sqrt{n}(\widehat{\mathbf{R}}_{n} - \mathbf{R}) - \sqrt{n} \left (\varphi \left (\frac{1}{n} \sum_{\ell=1}^{n} \mathbf{Z}^{(\ell)} \left (\mathbf{Z}^{(\ell)} \right )^{\text{T}} \right ) - \mathbf{R} \right  ) \xrightarrow{p} \mathbf{0}_{q \times q}
\end{equation*}
as $n \to \infty$, i.e. making use of the empirical correlation matrix based on a true Gaussian sample or based on a pseudo Gaussian sample, results in the same asymptotic expansion. Suppose further that $\mathbf{R} = \mathbf{U} \boldsymbol{\Lambda} \mathbf{U}^{\text{T}}$ is the eigendecomposition of $\mathbf{R}$. Then $\mathbf{Z}^{(\ell)} = \mathbf{U}\boldsymbol{\Lambda}^{1/2} \boldsymbol{\epsilon}^{(\ell)}$ for $\ell = 1,\dots, n$ and $\boldsymbol{\epsilon}^{(1)},\dots,\boldsymbol{\epsilon}^{(n)}$ a sample from $\mathcal{N}_{s}(\mathbf{0}_{q},\pmb{\mathbb{I}}_{q})$. From Lemma 4.16 of \cite{Mordant2022}, we have
\begin{equation*}
    \mathbf{W}_{n} = \frac{1}{\sqrt{n}} \sum_{\ell=1}^{n} (\boldsymbol{\epsilon}^{(\ell)} \left (\boldsymbol{\epsilon}^{(\ell)} \right )^{\text{T}} - \pmb{\mathbb{I}}_{q}) \xrightarrow{d} \mathbf{W},
\end{equation*}
as $n \to \infty$, where $\mathbf{W}$ is a random symmetric matrix with $\mathbf{W}_{jk} \sim \mathcal{N}(0,2)$ if $j = k \in \{1,\dots,q\}$ and $\mathbf{W}_{jk} \sim \mathcal{N}(0,1)$ if $1 \leq j < k \leq q$ independently (and similarly for $k < j$). Moreover, for $\mathbf{A},\mathbf{B} \in \mathbb{S}^{q}$, it holds that
\begin{equation*}
    \mathbb{E} \left (\text{Tr}(\mathbf{A}\mathbf{W})\text{Tr}(\mathbf{B}\mathbf{W}) \right ) = 2 \text{Tr}(\mathbf{A}\mathbf{B}).
\end{equation*}
We find 
\begin{equation*}
    \mathbf{U}\boldsymbol{\Lambda}^{1/2} \mathbf{W}_{n} \boldsymbol{\Lambda}^{1/2} \mathbf{U}^{\text{T}} = \sqrt{n} \left ( \frac{1}{n} \sum_{\ell=1}^{n} \left (\mathbf{Z}^{(\ell)}\left (\mathbf{Z}^{(\ell)} \right )^{\text{T}} \right ) - \mathbf{R} \right ) \xrightarrow{d} \mathbf{U}\boldsymbol{\Lambda}^{1/2} \mathbf{W} \boldsymbol{\Lambda}^{1/2} \mathbf{U}^{\text{T}},
\end{equation*}
as $n \to \infty$.
Applying the delta method (and using that $\varphi(\mathbf{R}) = \mathbf{D}_{\mathbf{R}}^{-1/2}\mathbf{R}\mathbf{D}_{\mathbf{R}}^{-1/2} = \mathbf{R}$), we obtain
\begin{equation*}
\begin{split}
    \sqrt{n} \left (\mathcal{D}_{\Phi}^{\mathcal{N}}(\widehat{\mathbf{R}}_{n}) - \mathcal{D}_{\Phi}^{\mathcal{N}}(\mathbf{R}) \right ) & \xrightarrow{d} \text{Tr} \left ((\mathbf{M}_{\Phi} - \mathbf{D}_{\mathbf{M}_{\Phi}\mathbf{R}})\mathbf{U}\boldsymbol{\Lambda}^{1/2} \mathbf{W} \boldsymbol{\Lambda}^{1/2} \mathbf{U}^{\text{T}} \right ) \\ & \hspace{0.2cm} = \text{Tr} \left (\boldsymbol{\Lambda}^{1/2} \mathbf{U}^{\text{T}} (\mathbf{M}_{\Phi} - \mathbf{D}_{\mathbf{M}_{\Phi} \mathbf{R}}) \mathbf{U} \boldsymbol{\Lambda}^{1/2} \mathbf{W} \right ),
\end{split}
\end{equation*}
as $n \to \infty$. The latter asymptotic expression is centered Gaussian with asymptotic variance 
\begin{equation*}
    2 \text{Tr} \left ( \left (\boldsymbol{\Lambda}^{1/2} \mathbf{U}^{\text{T}} (\mathbf{M}_{\Phi} - \mathbf{D}_{\mathbf{M}_{\Phi} \mathbf{R}}) \mathbf{U} \boldsymbol{\Lambda}^{1/2} \right )^{2} \right ) = 2 \text{Tr} \left ( \left (\mathbf{R}  (\mathbf{M}_{\Phi}-\mathbf{D}_{\mathbf{M}_{\Phi}\mathbf{R}}) \right )^{2} \right ),
\end{equation*}
using the trace cyclical property and finishing the proof. \hfill \qedsymbol{}
\newline \\ \noindent
\textbf{Proof of Theorem 2} \newline \\ \noindent
Note that because of the condition on $M_{n}$, both estimators $\widehat{\mathcal{D}}_{\Phi,n,M_{n}}$ and $\mathcal{D}_{\Phi}(\widehat{\boldsymbol{\theta}}_{n})$ have the same asymptotic normality result. The fact that $\mathcal{D}_{\Phi}(\widehat{\boldsymbol{\theta}}_{n})$ is asymptotically normal, follows from the delta method. If furthermore we can interchange differentiation and integration,  the Fr\'{e}chet derivative (total derivative) of $\boldsymbol{\theta}_{C} \mapsto \mathcal{D}_{\Phi}(\boldsymbol{\theta}_{C})$ in the direction of $\mathbf{h} = (h_{1},\dots,h_{D}) \in \mathbb{R}^{D}$ is given by 
\begin{equation*}
    \int_{\mathbb{I}^{q}} \sum_{i=1}^{D} h_{i} \frac{\partial f}{\partial \theta_{C,i}}(\mathbf{u};\boldsymbol{\theta}_{C}) d\mathbf{u} = \sum_{i=1}^{D} \left [\int_{\mathbb{I}^{q}} \frac{\partial f}{\partial \theta_{C,i}}(\mathbf{u};\boldsymbol{\theta}_{C}) d\mathbf{u} \right ] h_{i}.
\end{equation*}
Hence, an asymptotic normality result for $\widehat{\boldsymbol{\theta}}_{n}$,
\begin{equation*}
    \sqrt{n} \left (\widehat{\boldsymbol{\theta}}_{n} - \boldsymbol{\theta}_{C} \right ) \xrightarrow{d} \mathbf{Y} = (Y_{1},\dots,Y_{D})^{\text{T}} \sim \mathcal{N}(\mathbf{0},\mathbf{V}),
\end{equation*}
gives rise to an asymptotic normality result for $\mathcal{D}_{\Phi}(\widehat{\boldsymbol{\theta}}_{n})$, and thus for $\widehat{\mathcal{D}}_{\Phi,n,M_{n}}$,
\begin{equation*}
    \sqrt{n} \left (\widehat{\mathcal{D}}_{\Phi,n,M_{n}} - \mathcal{D}_{\Phi}(\boldsymbol{\theta}_{C}) \right ) \xrightarrow{d} \sum_{i=1}^{D} \left [\int_{\mathbb{I}^{q}} \frac{\partial f}{\partial \theta_{C,i}}(\mathbf{u},\boldsymbol{\theta}_{C}) d\mathbf{u} \right ] Y_{i} = \boldsymbol{\beta}^{\text{T}} \mathbf{Y}.
\end{equation*}
The result then follows from $\text{Var}(\boldsymbol{\beta}^{T}\mathbf{Y}) = \boldsymbol{\beta}^{\text{T}} \mathbf{V} \boldsymbol{\beta}$. \hfill \qedsymbol{}
\end{document}